\documentclass{article}
\usepackage{amsmath}[1996/11/01]
\usepackage{amssymb,amsthm,amsxtra}

\def\[#1\]{\begin{equation}#1\end{equation}}
\makeatletter
\def\beq{%
   \relax\ifmmode
      \@badmath
   \else
      \ifvmode
         \nointerlineskip
         \makebox[.6\linewidth]%
      \fi
      $$
   \fi
}
\def\eeq{%
   \relax\ifmmode
      \ifinner
         \@badmath
      \else
         $$
      \fi
   \else
      \@badmath
   \fi
   \ignorespaces
}

\def\enddisplaymath{\eeq\global\@ignoretrue}
\makeatother

\newtheorem{thm}{Theorem}
\newtheorem{cor}[thm]{Corollary}
\newtheorem{lem}[thm]{Lemma}
\newtheorem{prop}[thm]{Proposition}

\theoremstyle{remark}
\newtheorem*{rem}{Remark}
\newtheorem{rems}{Remark}[thm]

\theoremstyle{definition}
\newtheorem{defn}{Definition}
\newtheorem*{sq_defn}{``Definition''}

\numberwithin{equation}{section}
\numberwithin{thm}{section}
\numberwithin{eg}{section}

\newcommand{\C}{\mathbb C}
\newcommand{\CP}{\mathbb CP}

\newcommand{\Z}{\mathbb Z}

\DeclareMathOperator\Res{Res}

\DeclareMathOperator\sgn{sgn}

\DeclareMathOperator\cC{{\cal C}}
\DeclareMathOperator\cR{{\cal R}}

\newcommand{\II}{\mathord{I\!I}}
\newcommand{\blambda}{{\boldsymbol\lambda}}
\newcommand{\bkappa}{{\boldsymbol\kappa}}
\newcommand{\bmu}{{\boldsymbol\mu}}
\newcommand{\bnu}{{\boldsymbol\nu}}
\newcommand{\binomE}{\genfrac[]{0pt}{}}

\begin{document}

\title{Transformations of elliptic hypergeometric integrals} \author{Eric
  M. Rains\\Department of Mathematics, University of California, Davis}

\date{April 15, 2005}
\maketitle

\begin{abstract}
We prove a pair of transformations relating elliptic hypergeometric
integrals of different dimensions, corresponding to the root systems $BC_n$
and $A_n$; as a special case, we recover some integral identities
conjectured by van Diejen and Spiridonov.  For $BC_n$, we also consider
their ``Type II'' integral.  Their proof of that integral, together with
our transformation, gives rise to pairs of adjoint integral operators; a
different proof gives rise to pairs of adjoint difference operators.  These
allow us to construct a family of biorthogonal abelian functions
generalizing the Koornwinder polynomials, and satisfying the analogues of
the Macdonald conjectures.  Finally, we discuss some
transformations of Type II-style integrals.  In particular, we find
that adding two parameters to the Type II integral gives an integral
invariant under an appropriate action of the Weyl group $E_7$.

\end{abstract}

\tableofcontents

\section{Introduction}

In recent work, van Diejen and Spiridonov
\cite{vanDiejenJF/SpiridonovVP:2000,vanDiejenJF/SpiridonovVP:2001,SpiridonovVP:2003}
have produced a number of conjectural elliptic hypergeometric integration
formulae, common generalizations of Spiridonov's elliptic beta integral
\cite{SpiridonovVP:2001} and $q$-hypergeometric integration identities due
to Gustafson \cite{GustafsonRA:1992}.  In particular, for the $BC_n$ root
system, they gave two conjectures, ``Type I'' and ``Type II'' (proved as
Corollary \ref{cor:big_jackson} and Theorem \ref{thm:normalization} below),
and showed that the Type I integral would imply the Type II integral.
In an appropriate limit, their Type II integral transforms via residue
calculus into a sum originally conjectured by Warnaar \cite{WarnaarSO:2002}
(and proved by Rosengren \cite{RosengrenH:2001}).  In fact, Warnaar also
conjectured a more general formula, a Bailey-type transformation identity,
rather than a Jackson-type summation identity.  This suggests that there
should be transformation formulae on the integral level as well; this is
the topic of the present work.

Their Type I integral can be thought of as the ultimate generalization of
an integral identity used by Anderson \cite{AndersonGW:1991} in his proof
of the Selberg integral (which the Type II integral generalizes).  While
Anderson's proof of this identity (based on a clever change of variables)
does not appear to generalize any further, some recent investigations of
Forrester and the author \cite{ForresterPJ/RainsEM:2005} of a random matrix
interpretation of the Anderson integral suggested a different argument,
which as we will see does indeed generalize to the elliptic level.  While
the argument was not powerful enough to directly prove the Type I integral,
it {\em was} able to prove it for a countably infinite union of
submanifolds of parameter space.  This suggested that this argument should
at least suffice to produce the correct conjecture for a transformation
law; in the event, it turned out that it produced not only a conjecture but
a proof.  We thus obtain an identity relating an $n$-dimensional integral
with $2n+2m+4$ parameters to an $m$-dimensional integral with transformed
parameters; when $m=0$, this gives the van Diejen-Spiridonov integral, but
the proof requires this degree of freedom.  A similar identity for the
$A_n$ root system follows by a slight modification of the argument; this
gives a transformation generalization of a conjecture of Spiridonov
\cite{SpiridonovVP:2003}.  The basic idea for both proofs is that, in an
appropriate special case, the transformations can be written as
determinants of relatively simple one-dimensional transformations.  This
``determinantal'' case is thus easy to prove; moreover, by taking limits of
some of the remaining degrees of freedom, we can transform the
$n$-dimensional determinantal identity into a lower-dimensional, but
non-determinantal instance of the transformation.  Indeed, by repeating
this process, starting with a sufficiently large instance of the
determinantal case, we can obtain a dense set of special cases of the
desired transformation, thus proving the theorem.

As we mentioned, the Type II integral follows as a corollary of the Type I
integral.  In many ways, this integral is of greater interest, most notably
because it generalizes the inner product density for the Koornwinder
polynomials \cite{KoornwinderTH:1992}.  Since the inner product density
generalizes, it would be natural to suppose that the orthogonal polynomials
themselves should generalize.  It would be too much to expect them to
generalize to orthogonal functions, however; indeed, even in the univariate
case, the elliptic analogues of the Askey-Wilson polynomials are merely
{\em bi}orthogonal (these analogues are due to Spiridonov and Zhedanov
\cite{SpiridonovVP/ZhedanovAS:2000a,SpiridonovVP/ZhedanovAS:2000b} in the
discrete case (generalizing work of Wilson 
\cite{WilsonJA:1991}), and Spiridonov \cite{SpiridonovVP:2003} in the
continuous case (generalizing work of Rahman \cite{RahmanM:1986})).  With
this in mind, we will construct in the sequel a family of functions
satisfying biorthogonalilty with respect to the Type II integral.

There are two main ingredients in this construction.  The first is a family
of difference operators, generalizing some difference operators known to
act nicely on the Koornwinder polynomials \cite{bcpoly}, and satisfying
adjointness relations with respect to the elliptic inner product.  As a
special case, we obtain a difference-operator-based proof of the Type II
integral.  This suggests that the proof based on the Type I integral should
be related to a pair of adjoint {\em integral} operators, which form the
other main ingredient in our construction.  It turns out that the
$BC_n\leftrightarrow BC_m$ transformation plays an important role in
understanding these integral operators; indeed, by taking limits of the
transformation so that one side becomes a finite sum, we obtain formulas
for the images under the integral operator of a spanning set of its domain.
The biorthogonal functions are then constructed as the images of suitable
sequences of difference and integral operators.  (This construction is new
even at the level of Koornwinder polynomials.)

As these functions are biorthogonal with respect to a generalization of the
Koornwinder density (and indeed contain the Koornwinder polynomials as a
special case, although this turns out to be somewhat subtle to prove), one
of course expects that they satisfy analogous properties.  While the
Hecke-algebraic aspects of the Koornwinder theory (see, for instance,
\cite{SahiS:1999}) are still quite mysterious at the elliptic level, the
main properties, i.e., the ``Macdonald conjectures'', do indeed carry over.
Two of these properties, namely the closed forms for the principal
specialization and the nonzero values of the inner product, follow quite
easily from the construction and adjointness; the third (evaluation
symmetry) will be proved in a follow-up paper \cite{bctheta}.  The
arguments there are along the same lines as those given in \cite{bcpoly}
for the Koornwinder case, which were based on Okounkov's $BC_n$-symmetric
interpolation polynomials \cite{OkounkovA:1998a}.  Unlike in the Koornwinder
case, however, at the elliptic level the required interpolation functions
are actually special cases of the biorthogonal functions.

Just as in \cite{bcpoly}, these interpolation functions satisfy a number of
generalized hypergeometric identities, having Warnaar's multivariate
identities and conjectures \cite{WarnaarSO:2002} as special cases.  To be
precise, they satisfy multivariate analogues of Jackson's summation and
Bailey's transformation.  The former identity has an integral analogue,
extending the Type II integral.  In fact, the transformation also has an
analogue, stated as Theorem \ref{thm:int_bailey_lm} below; as a special
case, this gives our desired transformation analogue of the Type II
integral.  The simplest version of this transformation states that an
eight-parameter version of the Type II integral is invariant with respect
to an action of the Weyl group $E_7$; in fact, this action extends formally
to an action of $E_8$, acting on the parameters in the most natural way.

The plan of the paper is as follows.  After defining some notation at the
end of this introduction, we proceed in section 2 to discuss Anderson's
integral, as motivation for our proof of the $BC_n$ and $A_n$ integral
transformations.  These transformations are then stated and proved in
sections 3 and 4, respectively; we also briefly discuss in section 5 some
hybrid transformations arising from the fact that the $BC_1$ and $A_1$
integrals are the same, but the transformations are not.  Section 6 then
begins our discussion of the biorthogonal functions by describing the three
kinds of difference operators, as well as the spaces of functions on which
they act, and filtrations of those spaces with respect to which the
operators are triangular.  Section 7 discusses the corresponding integral
operators (all of which are special cases of a single operator defined from
the Type I integral), showing that they are triangular with respect to the
same filtrations as the difference operators.  Then, in section 8, we
combine these ingredients to construct the biorthogonal functions, and
describe their main properties.  In section 9, we discuss our Type II
transformations.  Finally, in an appendix, we give a general result
regarding when an integral of a meromorphic function is meromorphic, and
apply it to obtain precise information about the singularities of our
integrals.

The author would like to thank P. Forrester, A. Okounkov, H. Rosengren, and
V. Spiridonov for helpful comments on various drafts; also S. Ruijsenaars
for helpful conversations related to the appendix.  This work was supported
in part by NSF Grant No. DMS-0401387; in addition, some work was performed
while the author was employed by the Center for Communications Research,
Princeton.

\noindent{\bf Notation}

We will need a number of generalized $q$-symbols in the sequel.  First,
define the theta function and elliptic Gamma function
\cite{RuijsenaarsSNM:1997}:
\begin{align}
\theta(x;p)&:=\prod_{0\le k} (1-p^k x)(1-p^{k+1}/x)\\
\Gamma(x;p,q)&:=\prod_{0\le j,k} (1-p^{j+1}q^{k+1}/x)(1-p^j q^k x)^{-1}
\end{align}
In each case, the presence of multiple arguments before the semicolon
indicates a product; thus, for instance,
\[
\Gamma(z_i^{\pm 1}z_j^{\pm 1};p,q) =
\Gamma(z_iz_j;p,q)\Gamma(z_i/z_j;p,q)\Gamma(z_j/z_i;p,q)\Gamma(1/z_iz_j;p,q).
\]
These functions satisfy a number of identities, most notably
\[
\theta(x;p)=\theta(p/x;p)=(-x)\theta(1/x;p),
\]
and
\begin{align}
\Gamma(x;p,q) &= \Gamma(pq/x;p,q)^{-1}\\
\Gamma(px;p,q)&=\theta(x;q)\Gamma(x;p,q)\\
\Gamma(qx;p,q)&=\theta(x;p)\Gamma(x;p,q).
\end{align}

Using the theta function, one can define an elliptic analogue of the
$q$-symbol; in fact, just as the elliptic Gamma function is symmetric in
$p$ and $q$, we will want our elliptic $q$-symbol to also be symmetric.
Thus, we define
\[
\theta(x;p,q)_{l,m}:=\prod_{0\le k<l}\theta(p^k x;q)\prod_{0\le
  k<m}\theta(q^k x;p),
\]
so that
\[
\frac{\Gamma((p,q)^{l,m}x;p,q)}{\Gamma(x;p,q)}
=
(-x)^{-lm} p^{-lm(l-1)/2} q^{-lm(m-1)/2}
\theta(x;p,q)_{l,m},
\]
where
\[
(p,q)^{l,m}:=p^lq^m.
\]

We also need some multivariate symbols, indexed (as the biorthogonal
functions will be) by pairs of partitions.  By convention, we will use bold
greek letters to refer to such partition pairs, and extend transformations
and relations of partitions in the obvious way.  We then define, following
\cite{bcpoly},
\begin{align}
\cC^0_{\blambda}(x;t;p,q)
&:=
\prod_{1\le i}\theta(t^{1-i} x;p,q)_{\blambda_i}
\\
\cC^-_{\blambda}(x;t;p,q)
&:=
\prod_{1\le i\le j}
\frac{\theta(t^{j-i} x;p,q)_{\blambda_i-\blambda_{j+1}}}
     {\theta(t^{j-i} x;p,q)_{\blambda_i-\blambda_j}}
\\
\cC^+_{\blambda}(x;t;p,q)
&:=
\prod_{1\le i\le j}
\frac{\theta(t^{2-i-j} x;p,q)_{\blambda_i+\blambda_j}}
     {\theta(t^{2-i-j} x;p,q)_{\blambda_i+\blambda_{j+1}}}.
\end{align}
We note that each of the above $\cC$ symbols extends to a holomorphic
function on $x\in \C^*$.

Two particular combinations of $\cC$ symbols will occur frequently enough
to merit their own notation.  We define:
\begin{align}
\Delta^{0}_{\blambda}(a|\dots b_i\dots;t;p,q)
&:=
\frac{\cC^0_{\blambda}(\dots b_i\dots;t;p,q)}
     {\cC^0_{\blambda}(\dots pqa/b_i\dots;t;p,q)}\\
\Delta_{\blambda}(a|\dots b_i\dots;t;p,q)
&:=
\Delta^{0}_{\blambda}(a|\dots b_i\dots;t;p,q)
\frac{\cC^0_{2\blambda^2}(pqa;t;p,q)}
     {\cC^-_{\blambda}(pq,t;t;p,q)\cC^+_{\blambda}(a,pqa/t;t;p,q)}
\end{align}

We will also need the following notion, where $0<|p|<1$.

\begin{defn}
A $BC_n$-symmetric ($p$-)theta function of degree $m$ is a holomorphic function
$f(x_1,\dots x_n)$ on $(\C^*)^n$ such that
\begin{itemize}
\item[] $f(x_1,\dots x_n)$ is invariant under permutations of its
  arguments.
\item[] $f(x_1,\dots x_n)$ is invariant under $x_i\mapsto 1/x_i$ for each
  $i$.
\item[] $f(px_1,x_2,\dots x_n)=(1/pz_i^2)^m f(x_1,x_2,\dots x_n)$.
\end{itemize}
A $BC_n$-symmetric ($p$-)abelian function is a meromorphic function
satisfying the above conditions with $m=0$.
\end{defn}

In particular, a $BC_n$-symmetric theta function of degree $m$ is a
$BC_1$-symmetric theta function of degree $m$ in each of its
arguments.  Now, the space of $BC_1$-symmetric theta functions of
degree $m$ is $m+1$-dimensional, and moreover, any nonzero $BC_1$-symmetric
theta function vanishes at exactly $2m$ orbits of points (under
multiplication by $p$, and counting multiplicity).  Thus we can show that
a $BC_1$-symmetric theta function vanishes by finding $m+1$
independent points at which it vanishes.

The canonical example of a $BC_n$-symmetric theta function of degree $1$ is
\[
\prod_{1\le i\le n} \theta(u x_i^{\pm 1};p);
\]
indeed, the functions for any $n+1$ distinct choices of $u$ form a basis of
the space of $BC_n$-symmetric theta functions of degree $1$:
\[
f(\dots x_i\dots)
=
\sum_{0\le j\le n}
f(u_0,u_1,\dots u_{j-1},u_{j+1},\dots u_n)
\frac{\prod_{1\le i\le n} \theta(u_j x_i^{\pm 1};p)}
     {\prod_{k\ne j} \theta(u_j u_k^{\pm 1};p)},
\]
for any $BC_n$-symmetric theta function $f$ of degree $1$.  More generally,
the space of $BC_n$-symmetric theta functions of degree $m$ is spanned by
the set of products of $m$ such functions.

\section{The Anderson integral}

Many of our arguments in the sequel were inspired by considerations of an
extremely special (but quite important) case of Corollary
\ref{cor:big_jackson} below, a multivariate integral identity used in
Anderson's proof of the Selberg integral.

\begin{thm} \cite{AndersonGW:1991}
Let $x_1,\dots x_n$ and $s_1,\dots s_n$ be sequences of real
numbers such that
\[
x_1>x_2>\dots >x_n
\quad\text{and}\quad
0<s_1,s_2,\dots s_n
\]
Then
\[
\int_{x_n\le y_{n-1}\le \dots \le x_2\le y_1\le x_1}
\frac{
\prod_{1\le i<j\le n-1} |y_i-y_j|
\prod_{\substack{1\le i\le n-1\\1\le j\le n}}
|y_i-x_j|^{s_j-1}}
{\prod_{1\le i<j\le n} |x_i-x_j|^{s_i+s_j-1}}
\prod_{1\le i\le n-1} dy_i
=
\frac{\prod_{1\le i\le n} \Gamma(s_i)}{\Gamma(S)},
\]
where $S=\sum_{1\le i\le n}s_i$.
\end{thm}

\begin{rem}
In fact, although Anderson independently discovered the above integral, it
turns out that a more general identity (analogous to Theorem
\ref{thm:big_bailey} below) was discovered in 1905 by Dixon
\cite{DixonAL:1905}; see also the remark above Theorem
\ref{thm:anderson:schur_action} below.  However, Anderson was the first to
notice the significance of this special case in the theory of the Selberg
integral, so we will refer to it as the Anderson integral in the sequel.
\end{rem}

Since the integrand is nonnegative, we can normalize to obtain a
probability distribution.  It turns out that if the $s_i$ parameters are
all positive integers, then this probability distribution has a natural
random matrix interpretation.

\begin{thm} \cite{BaryshnikovY:2001} \cite[\S 4]{ForresterPJ/RainsEM:2005}
Let $A$ be a $S\times S$ complex Hermitian matrix, let $x_1>x_2>\dots>
x_n$ be the list of distinct eigenvalues of $A$, and let $s_1$, $s_2$,\dots
$s_n$ be the corresponding multiplicities.  Let $\Pi:\C^{S}\to
\C^{S-1}$ be the orthogonal projection onto a hyperplane chosen uniformly
at random.  Then the $S-1\times S-1$ Hermitian matrix $B=\Pi A \Pi^\dagger$
has eigenvalues $x_i$ with multiplicity $s_i-1$, together with $n-1$ more
eigenvalues $y_i$, distributed according to the Anderson distribution.
\end{thm}

\begin{rem}
A similar statement holds over the reals, except that now the
multiplicities correspond to $2s_1$, $2s_2$,\dots  Similarly, over the
quaternions, the multiplicities correspond to $s_1/2$, $s_2/2$, \dots
\end{rem}

Of particular interest is the generic case, in which the eigenvalues of $A$
are all distinct; that is $s_1=s_2=\dots s_n=1$ (this is the case
considered in \cite{BaryshnikovY:2001}).  In this case, the Anderson
integral is particularly simple to prove.  Indeed, the relevant integral
is:
\[
\frac{(n-1)!}{\prod_{1\le i<j\le n} (x_i-x_j)}
\int_{x_n\le y_{n-1}\le \dots \le x_2\le y_1\le x_1}
\prod_{1\le i<j\le n-1} (y_i-y_j)
\prod_{1\le i\le n-1} dy_i
\]
In particular, the integrand is simply a Vandermonde determinant,
\[
\prod_{1\le i<j\le n-1} (y_i-y_j)
=
(-1)^{n(n-1)/2}
\det_{1\le i,j\le n-1} (y_i-x_n)^{j-1}.
\]
Integrating this row-by-row gives
\begin{align}
\frac{(-1)^{n(n-1)/2}(n-1)!}{\prod_{1\le i<j\le n} (x_i-x_j)}
\det_{1\le i,j\le n-1} \int_{x_{i+1}}^{x_i} (y-x_n)^{j-1} dy
&=
\frac{(-1)^{n(n-1)/2}(n-1)!}{\prod_{1\le i<j\le n} (x_i-x_j)}
\det_{1\le i,j\le n-1} \int_{x_n}^{x_i} (y-x_n)^{j-1} dy\\
&=
\frac{(-1)^{n(n-1)/2}(n-1)!}{\prod_{1\le i<j\le n} (x_i-x_j)}
\det_{1\le i,j\le n-1} \frac{(x_i-x_n)^j}{j}\\
&=
1.
\end{align}

Now, in general, a Hermitian matrix with multiple eigenvalues can be
expressed as a limit of matrices with distinct eigenvalues; this suggests
that we should be able to obtain the general integer $s$ Anderson integral
as a limit of $s\equiv 1$ Anderson integrals.  Indeed, if we integrate over
$y_i$ and take a limit $x_{i+1}\to x_i$, the result is simply the Anderson
distribution with parameters
\[
x_1>\dots x_i>x_{i+2}>\dots x_n
\quad\text{and}\quad
s_1, s_2, \dots s_i+s_{i+1},\dots s_{n-1}, s_n.
\]
Combined with the determinantal proof for $s\equiv 1$, we thus obtain by
induction a proof of Anderson's integral for arbitrary positive integer
$s$.  We can then obtain the general case via analytic continuation (for
which we omit the argument, as it greatly simplifies in the cases of
interest below).  The resulting proof is less elegant than Anderson's
original proof; however, it has the distinct advantage for our purposes of
extending to much more general identities.  Indeed, our proofs of Theorems
\ref{thm:big_bailey} and \ref{thm:bb_An} below proceed by precisely this
sort of induction from large dimensional, but simple, cases.

\begin{rem}
Note that the key property of the ``determinantal'' case is not so much
that it is a determinant, but that it is a determinant of univariate
instances of the Anderson integral.  Indeed, the general Anderson integral
can be expressed as a determinant of univariate integrals; in fact, a
generalization of the resulting identity was proved by Varchenko
\cite{VarchenkoAN:1989} even before Anderson's work \cite{AndersonGW:1991},
but without noticing that it could be used to prove the Selberg integral.
See also \cite{RichardsD/ZhengQ:2002} (apparently the first article to
observe that Varchenko's identity could be expressed as a multivariate
integral).  It would be interesting to know if Varchenko's generalized
formula can be extended to the elliptic level.
\end{rem}

The random matrix interpretation also gives the following result.
Given a symmetric function $f$, we define $f(A)$ for a matrix $A$ to be
$f$ evaluated at the multiset of eigenvalues of $A$.

\begin{thm}\label{thm:anderson:schur_action}
Let $A$ be an $n$-dimensional Hermitian matrix, and let $\Pi$ be a random
orthogonal projection as before.  Then for any partition $\lambda$,
\[
{\bf E}_\Pi s_\lambda(\Pi A\Pi^\dagger)
=
\frac{s_\lambda(1_{n-1})s_\lambda(A)}{s_\lambda(1_n)}.
\]
\end{thm}

\begin{proof}
Since $\Pi$ was uniformly distributed, we have
\[
{\bf E}_\Pi s_\lambda(\Pi A\Pi^\dagger)
=
{\bf E}_\Pi s_\lambda(\Pi U A U^\dagger\Pi^\dagger)
\]
for any unitary matrix $U$.  In particular, we can fix $\Pi$ and take
expectations over $U$, thus obtaining
\[
{\bf E}_U s_\lambda(\Pi U A U^\dagger\Pi^\dagger)
=
{\bf E}_U s_\lambda(U A U^\dagger\Pi^\dagger \Pi)
=
\frac{s_\lambda(A)s_\lambda(\Pi^\dagger \Pi)}{s_\lambda(1_n)}
=
\frac{s_\lambda(A)s_\lambda(1_{n-1})}{s_\lambda(1_n)}.
\]
Here we have used the fact
\[
{\bf E}_U
s_\lambda(U A U^\dagger B)
=
\frac{s_\lambda(A)s_\lambda(B)}{s_\lambda(1_n)}
\]
from the theory of zonal polynomials, or equivalently from the fact that
Schur functions are irreducible characters of the unitary group.
\end{proof}

In other words, the Anderson distribution for $s\equiv 1$ acts as a raising
integral operator on Schur polynomials, taking an $n-1$-variable Schur
polynomial to the corresponding $n$-variable Schur polynomial.  Similarly,
the Anderson densities for $s\equiv \frac{1}{2}$ and $s\equiv 2$ act as
raising operators on the real and quaternionic zonal polynomials.  This
suggests that in general, an Anderson distribution with constant $s$ should
take polynomials to polynomials (mapping an appropriate Jack polynomial to
the corresponding $n$-variable Jack polynomial).

Indeed, we have the following fact, even for nonconstant $s$.

\begin{thm}
Let $y_1$, $y_2$,\dots $y_{n-1}$ be distributed according to the Anderson
distribution with parameters $s_1$, \dots $s_n>0$, $x_1>\dots >x_n>0$.
Then as a function of $a_1$, $a_2$,\dots $a_m$,
\[
{\bf E}_y(\prod_{\substack{1\le i\le n-1\\1\le j\le m}} (a_j-y_i))
=
\frac{\Gamma(S)}{\Gamma(S+m)}
\frac{\prod_{\substack{1\le i\le m\\1\le j\le n}} (a_i-x_j)^{1-s_j}}
     {\prod_{1\le i<j\le m} (a_j-a_i)}
\prod_{1\le i\le m}
\frac{\partial}{\partial a_i} \prod_{1\le j\le n} (a_i-x_j)^{s_j}
\prod_{1\le i<j\le m} (a_j-a_i).
\]
In particular, the left-hand side is a polynomial in the $x_j$.
\end{thm}

\begin{proof}
If $m=0$, this simply states that ${\bf E}_y(1)=1$; we may thus proceed by
induction on $m$.  Suppose the theorem holds for $m=m_0$, and consider what
becomes of that instance when $s_n=1$.  In that case, the density is
essentially independent of $x_n$, in that $x_n$ only affects the
normalization and the domain of integration.  Thus if we multiply both
sides by $\prod_{0\le i<j\le n} (x_i-x_j)^{1-s_i-s_j}$, we can
differentiate by $x_n$ to obtain an $n-1$-dimensional integral.  If we then
set $x_n=a_{m_0+1}$ and renormalize, the result is the $n-1$-dimensional case
of the theorem with $m=m_0+1$.

That the right-hand side is a polynomial in the $x_i$ is straightforward,
and thus the left-hand side is also polynomial.
\end{proof}

\begin{rems}
Compare the proof of Theorem \ref{thm:int_op_main} given in Remark
\ref{rems:int_op_main_alt_proof} following the theorem.
\end{rems}

\begin{rems}
The left-hand side of the above identity was studied by Barsky and
Carpentier \cite{BarskyD/CarpentierM:1996} using Anderson's change of
variables; they did not obtain as simple a right-hand side, however.
\end{rems}

\begin{cor}
As an integral operator, the Anderson distribution takes symmetric functions to
polynomials; if $s_1=s_2=\dots s_n=s$, it maps symmetric functions to
symmetric functions.
\end{cor}

\begin{rem}
A $q$-integral analogue of the Corollary was proved by Okounkov
\cite{OkounkovA:1998b}, who credits a private communication from Olshanski
for the Corollary itself.  
\end{rem}

We can also obtain integral operators on symmetric functions by fixing one
or two of the $x$ parameters and allowing their multiplicities to vary; the
result is then a symmetric function in the remaining $x$ parameters.  In
particular, Anderson's proof of the Selberg integral acquires an
interpretation in terms of pairs of adjoint integral operators.

\section{The $BC_n\leftrightarrow BC_m$ transformation}

For all nonnegative integers m,n, and parameters $p$, $q$, $t_0$\dots
$t_{2m+2n+3}$ satisfying
\[
|p|,|q|,|t_0|,\dots |t_{2m+2n+3}|<1,\quad \prod_{0\le r\le 2m+2n+3} t_r =
 (pq)^{m+1},
\]
define
\[
I^{(m)}_{BC_n}(t_0,t_1,\dots ;p,q)
:=
\frac{(p;p)^n(q;q)^n}{2^n n!}
\int_{T^n}
\frac{\prod_{1\le i\le n}\prod_{0\le r\le 2m+2n+3} \Gamma(t_r z_i^{\pm 1};p,q)}
{\prod_{1\le i<j\le n} \Gamma(z_i^{\pm 1}z_j^{\pm 1};p,q)
\prod_{1\le i\le n} \Gamma(z_i^{\pm 2};p,q)}
\prod_{1\le i\le n} \frac{dz_i}{2\pi\sqrt{-1}z_i},
\]
a contour integral over the unit torus.  We can extend this to a
meromorphic function on the set
\[
{\cal P}_{mn}:=\{(t_0,t_1,\dots t_{2m+2n+3},p,q)\mid \prod_{0\le r\le
 2m+2n+3} t_r=(pq)^{m+1},\ 0<|p|,|q|<1\}
\]
by replacing the unit torus with the $n$-th power of an arbitrary (possibly
disconnected) contour that contains the points of the form $p^i q^j t_r$,
$i,j\ge 0$ and excludes their reciprocals.  We thus find that the resulting
function is analytic away from points where $t_r t_s=p^{-i} q^{-j}$ for
some $0\le r,s\le 2m+2n+3$, $0\le i,j$.  (In fact, its singularities
consist precisely of simple poles along the hypersurfaces $t_r
t_s=p^{-i}q^{-j}$ with $0\le r<s\le 2m+2n+3$, $0\le i,j$; see the appendix.)

Note in particular that $I^{(m)}_{BC_0}(t_0,t_1,\dots t_{2m+3};p,q)=1$.

\begin{thm}\label{thm:big_bailey}
The following holds for $m,n\ge 0$ as an identity in meromorphic
functions on ${\cal P}_{mn}$.
\[
\label{eq:big_bailey}
I^{(m)}_{BC_n}(t_0,t_1,\dots t_{2m+2n+3};p,q)
=
\prod_{0\le r<s\le 2m+2n+3} \Gamma(t_r t_s;p,q)\,
I^{(n)}_{BC_m}(\frac{\sqrt{pq}}{t_0},\frac{\sqrt{pq}}{t_1},\dots 
\frac{\sqrt{pq}}{t_{2m+2n+3}};p,q).
\]
\end{thm}

\begin{rem}
If $\sqrt{pq}<|t_r|<1$ for all $r$, both contours may be taken to
be the unit torus.
\end{rem}

Taking $m=0$ gives the following:

\begin{cor}\label{cor:big_jackson}
\[
I^{(0)}_{BC_n}(t_0,t_1,\dots t_{2n+3};p,q)
=
\prod_{0\le r<s\le 2n+3} \Gamma(t_r t_s;p,q),
\]
\end{cor}

This is the ``Type I'' identity conjectured by van Diejen and Spiridonov
\cite{vanDiejenJF/SpiridonovVP:2001}, who showed that it would follow from
the fact that the integral vanishes of $t_0t_1=pq$, but were unable to
prove that fact.  The case $n=1$ of the $m=0$ integral is, however, known:
it is an elliptic beta integral due to Spiridonov \cite{SpiridonovVP:2001};
it also happens to agree with the case $n=1$ of Theorem
\ref{thm:normalization} below.  A direct proof of the corollary has since
been given by Spiridonov \cite{SpiridonovVP:2004}; see also the remark
following the second proof of Theorem \ref{thm:normalization} below.

Our strategy for proving Theorem \ref{thm:big_bailey} is as follows.  We
first observe that in a certain extremely special (``determinantal'') case,
each integrand can be expressed as a product of simple determinants, and
thus the integrals themselves can be expressed as determinants.  The
agreement of corresponding entries of the determinants then follows from
the $m=n=1$ instance of the determinantal case (Lemma
\ref{lem:ICn_basecase} below).

The next crucial observation is that if we take the limit $t_1\to pq/t_0$
in an instance of the general identity for given values of $m$ and $n$, the
result is an instance of the general identity with $m$ diminished by 1;
similarly, the limit $t_1\to 1/t_0$ decreases $n$ by 1.  It turns out,
however, that the determinantal case is {\em not} preserved by those
operations; thus by starting with ever larger determinantal cases and
dropping down to the desired $m$ and $n$, we obtain an ever increasing
collection of proved special cases of the identity.  The full set of
special cases obtained is in fact dense, and thus the theorem will follow.

\begin{rem}
A similar inductive argument based on a determinantal case was applied in
\cite{RosengrenH:2004} to prove the summation analogue of Corollary
\ref{cor:big_jackson} (see also Remark \ref{rems:int_op_main_rosengren}
following Theorem \ref{thm:int_op_main}); it is worth noting, therefore,
that the present argument is not in fact a generalization of Rosengren's.
This is not to say that the arguments are unrelated; indeed, in a sense,
the two arguments are dual.  In fact, Rosengren's determinantal case turns
out to be precisely Lemma \ref{lem:diffeq_norm} below, which is thus
related to the difference operators we will be considering in the sequel.
These, in turn, are related (by Theorem \ref{thm:int_op_main}, among other
things) to the integral operators we will define using Theorem
\ref{thm:big_bailey}.  The duality is most apparent on the series level; if
one interprets the sum as a sum over partitions, the two arguments are
precisely related by conjugation of partitions.  The main distinction for
our purposes is that Rosengren's argument, while superior in the series
case (as it does not require analytic continuation), does not appear to
extend to the integral case.
\end{rem}

The base case for the determinantal identity is the following:

\begin{lem}\label{lem:ICn_basecase}
The theorem holds for $m=n=1$, assuming the parameters have the form
\[
(t_0,t_1,t_2,\dots t_7) = (a_0,q/a_0,a_1,q/a_1,b_0,p/b_0,b_1,p/b_1)
\]
In other words, if we define
\[
F(a_0,a_1{:}b_0,b_1;p,q):=
\frac{(p;p)(q;q)}{2}
\int_C
\frac{\theta(z^2;q)\theta(z^{-2};p)}
{\theta(a_0 z^{\pm 1},a_1 z^{\pm 1};q)
\theta(b_0 z^{\pm 1},b_1 z^{\pm 1};p)}
\frac{dz}{2\pi\sqrt{-1}z}
,
\]
with contour as appropriate, then
\[
\label{eq:basecase_xform}
F(a_0,a_1{:}b_0,b_1;p,q)
=
\frac{\theta(a_0a_1/q,a_0/a_1;p)\theta(b_0b_1/p,b_0/b_1;q)}
{\theta(a_0a_1/q,a_0/a_1;q)\theta(b_0b_1/p,b_0/b_1;p)}\\
F(\frac{\sqrt{pq}}{b_0}, \frac{\sqrt{pq}}{b_1}{:}\frac{\sqrt{pq}}{a_0}, \frac{\sqrt{pq}}{a_1};p,q).
\]
\end{lem}

\begin{proof}
We first observe that it suffices to prove the Laurent series expansion
\[
b_1 \theta(b_0 b_1^{\pm 1};p)
\frac{z^{-1}\theta(z^2;p)}{\theta(b_0 z^{\pm 1},b_1 z^{\pm 1};p)}
=
(p;p)^{-2}
\sum_{k\ne 0} 
\frac{b_0^k+(p/b_0)^k-b_1^k-(p/b_1)^k}
{1-p^k} z^k
,
\]
valid for $|p|<|b_0|,|b_1|<1$, and $z$ in a neighborhood of the unit
circle.  Indeed, the desired integral is the constant term of the product
of two such expressions, and is thus expressed as an infinite sum, each
term of which already satisfies the desired transformation!

Consider the sum
\[
\sum_{k>0} 
\frac{b_0^k+(p/b_0)^k-b_1^k-(p/b_1)^k}
{1-p^k} \frac{z^k+z^{-k}}{k},
\]
which clearly differentiates (by $z\frac{d}{dz}$) to the stated sum.
If we expand $(1-p^k)^{-1}$ in a geometric series, each term can then be
summed over $k$ as a linear combination of logarithms.  We conclude that
\[
\sum_{k>0} 
\frac{b_0^k+(p/b_0)^k-b_1^k-(p/b_1)^k}
{1-p^k} \frac{z^k+z^{-k}}{k}
=
\log\left(\frac{\theta(b_1 z^{\pm 1};p)}{\theta(b_0 z^{\pm 1};p)}\right).
\]

Now, the derivative
\[
z\frac{d}{dz}
\log\left(\frac{\theta(b_1 z^{\pm 1};p)}{\theta(b_0 z^{\pm 1};p)}\right)
\]
is an elliptic function antisymmetric under $z\mapsto z^{-1}$, with only
simple poles, and those at points of the form $(p^k b_j)^{\pm 1}$.  It
follows that
\[
\frac{d}{d\log z}
\log\left(\frac{\theta(b_1 z^{\pm 1};p)}{\theta(b_0 z^{\pm 1};p)}\right)
=
C(b_0,b_1,p) \frac{z^{-1} \theta(z^2;p)}{\theta(b_0z^{\pm 1},b_1z^{\pm 1};p)},
\]
for some factor $C(b_0,b_1,p)$ independent of $z$.  Comparing asymptotics
at $z=b_0$ gives the desired result.
\end{proof}

\begin{rem}
If we take the limit $p\to 1$ in the above Laurent series expansion, we
obtain
\[
\lim_{p\to 1}
(1-p) (p;p)^2
b_1 \theta(b_0 b_1^{\pm 1};p)
\frac{z^{-1}\theta(z^2;p)}{\theta(b_0 z^{\pm 1},b_1 z^{\pm 1};p)}
=
\sum_{k\ne 0} 
\frac{b_0^k+b_0^{-k}-b_1^k-b_1^{-k}}
{k} z^k
.
\]
If $|b_0|=|b_1|=1$, $\Re(b_0)>\Re(b_1)$, then the limit is (up to a factor
of $2\pi\sqrt{-1}\sgn(\Im(z))$) the Fourier series expansion of the
indicator function for the arcs such that $\Re(b_0)\ge \Re(z)\ge \Re(b_1)$.
In particular, this explains how an integral like the Anderson integral,
with its relatively complicated domain of integration, can be a limiting
case of Corollary \ref{cor:big_jackson}.  The corresponding limit applied
to Theorem \ref{thm:big_bailey} gives an identity of Dixon \cite{DixonAL:1905}.
\end{rem}

\begin{lem}\label{lem:bb_detcase}
If $m=n$, then \eqref{eq:big_bailey} holds on the codimension $2n+2$
subset parametrized by:
\[
t_{2r} = a_r,
t_{2r+1} = q/a_r,
t_{2n+2+2r}=b_r,
t_{2n+2+2r+1}=p/b_r.
\]
\end{lem}

\begin{proof}
By taking a determinant of instances of \eqref{eq:basecase_xform}, we
obtain the identity:
\begin{align}
\det_{1\le i,j\le n}&\left(
\int_{C^n}
\frac{\theta(z^2;p)\theta(z^{-2};q)}
{\theta(a_0 z^{\pm 1},a_i z^{\pm 1};q)
\theta(b_0 z^{\pm 1},b_j z^{\pm 1};p)}
\frac{dz}{2\pi\sqrt{-1}z}
\right)
=\\
&\!\!\!
\det_{1\le i,j\le n}\left(
\frac{\theta(a_0a_i/q,a_0/a_i;p)\theta(b_0b_j/p,b_0/b_j;q)}
{\theta(a_0a_i/q,a_0/a_i;q)\theta(b_0b_j/p,b_0/b_j;p)}
\int_{C^{\prime n}}
\frac{\theta(z^2;p)\theta(z^{-2};q)}
{\theta(\frac{\sqrt{pq}}{a_0} z^{\pm 1},\frac{\sqrt{pq}}{a_i} z^{\pm 1};p)
\theta(\frac{\sqrt{pq}}{b_0} z^{\pm 1},\frac{\sqrt{pq}}{b_j} z^{\pm 1};q)}
\frac{dz}{2\pi\sqrt{-1}z}
\right)\notag
\end{align}
Consider the determinant on the left.  The $p$-theta functions in that
integral are independent of $j$, while the $q$-theta functions are
independent of $i$.  We may thus expand that determinant of integrals as an
integral of a product of two determinants:
\begin{align}
\det_{1\le i,j\le n}&\left(
\int_{C^n}
\frac{\theta(z^2;p)\theta(z^{-2};q)}
{\theta(a_0 z^{\pm 1},a_i z^{\pm 1};q)
\theta(b_0 z^{\pm 1},b_j z^{\pm 1};p)}
\frac{dz}{2\pi\sqrt{-1}z}
\right)\notag\\
&\qquad=
\frac{1}{n!}
\int_{C^n}
\det_{1\le i,j\le n}\left(
\frac{\theta(z_j^{-2};q)}
{\theta(a_0 z_j^{\pm 1},a_i z_j^{\pm 1};q)}
\right)
\det_{1\le j,i\le n}\left(
\frac{\theta(z_j^2;p)}
{\theta(b_0 z_j^{\pm 1},b_i z_j^{\pm 1};p)}
\right)
\prod_{1\le j\le n} \frac{dz_j}{2\pi\sqrt{-1}z_j}
\end{align}
These determinants can in turn be explicitly evaluated, using the following
identity:
\[
\det_{1\le i,j\le n}\left(
\frac{1}{a_i^{-1} \theta(a_i z_j^{\pm 1};q)}
\right)
=
(-1)^{n(n-1)/2}
\frac{
\prod_{1\le i<j\le n} a_i^{-1}\theta(a_i a_j^{\pm 1};q)
\prod_{1\le i<j\le n} z_i^{-1}\theta(z_i z_j^{\pm 1};q)}
{\prod_{1\le i,j\le n} a_i^{-1} \theta(a_i z_j^{\pm 1};q)}.
\]
(This is, for instance, a special case of a determinant identity of Warnaar
\cite{WarnaarSO:2002}, and can also be obtained as a special case of the
Cauchy determinant.)  The resulting identity is precisely the desired
result.
\end{proof}

As we mentioned above, the other key element to the proof is an
understanding of the limit of \eqref{eq:big_bailey} as $t_1\to pq/t_0$.  On
the left-hand side, the integral is perfectly well-defined when
$t_1=pq/t_0$, but the right-hand side ends up identifying two poles that
should be separated.  Thus we need to understand how
$I^{(m)}_{BC_n}(t_0,t_1,\dots)$ behaves as $t_1\to 1/t_0$.

\begin{lem}\label{lem:bb_reduction}
We have the limit:
\[
\lim_{t_1\to t_0^{-1}}
\frac{I^{(m)}_{BC_n}(t_0,t_1,\dots ;p,q)}
{
\Gamma(t_0t_1;p,q)
\prod_{2\le r<2m+2n+3} \Gamma(t_0 t_r,t_1 t_r;p,q)}
=
I^{(m-1)}_{BC_n}(t_2,t_3,\dots ;p,q)
\]
\end{lem}

\begin{proof}
If we deform the contour on the left through the points $t_1$ and $1/t_1$,
the resulting integral will have a finite limit, and will thus be
annihilated by the factor of $\Gamma(t_0t_1)$ in the denominator.  In other
words, the desired limit is precisely the limit of the sum of residues
corresponding to the change of contour.  By symmetry, each variable
contributes equally, as do $t_1$ and $1/t_1$; we thus find (using the
identity
\[
\lim_{y\to x} \Gamma(x/y)(1-x/y) = 1/(p;p)(q;q),
\]
which is easily verified):
\begin{align}
\lim_{t_1\to t_0^{-1}}&
\frac{I^{(m)}_{BC_n}(t_0,t_1,\dots ;p,q)}
{
\Gamma(t_0t_1;p,q)
\prod_{2\le r<2m+2n+3} \Gamma(t_r t_0,t_r t_1;p,q)}\notag\\
&=
\lim_{t_1\to t_0^{-1}}
\frac{(p;p)^{n-1}(q;q)^{n-1}}
     {2^{n-1}(n-1)!}
\frac{\Gamma(t_0/t_1;p,q)\prod_{2\le r<2m+2n+3} \Gamma(t_r/t_1;p,q)}
{\Gamma(1/t_1^2;p,q)\prod_{2\le r<2m+2n+3} \Gamma(t_0 t_r;p,q)}
\\
&\qquad\qquad
\int_{C^n}
\frac{\prod_{1\le i<n}\Gamma(t_0 z_i^{\pm 1};p,q)}
     {\prod_{1\le i<n} \Gamma(z_i^{\pm 1}/t_1;p,q)}
\frac{
\prod_{1\le i<n}\prod_{2\le r<2m+2n+3} \Gamma(t_r z_i^{\pm 1};p,q)}
{
\prod_{1\le i<j<n} \Gamma(z_i^{\pm 1}z_j^{\pm 1};p,q)
\prod_{1\le i<n} \Gamma(z_i^{\pm 2};p,q)}
\prod_{1\le i<n} \frac{dz_i}{2\pi\sqrt{-1}z_i}\notag
\\
&=
I^{(m-1)}_{BC_n}(t_2,t_3,\dots t_{2m+2n+3};p,q)
\end{align}
as required.
\end{proof}

We can now prove Theorem \ref{thm:big_bailey}.

\begin{proof}
For $m$, $n\ge 0$, let ${\cal C}_{mn}$ be the set of
parameters $c_0c_1\dots c_{m+n+1}=(pq)^{m+1}$ such that the theorem holds
on the manifold with
\[
t_{2i}t_{2i+1}=c_i, 0\le i\le m+n+1.
\]
Thus, for instance, Lemma \ref{lem:bb_detcase} states that the point
$(q,q,q,\dots q,p,p,p\dots p)$ is in ${\cal C}_{nn}$.

The key idea is that if $(c_0,c_1,c_2,c_3,\dots ,c_{m+n+1})\in {\cal C}_{mn}$, then
we also have:
\begin{align}
(c_0c_1,c_2,c_3,\dots c_{m+n+1})&\in {\cal C}_{m(n-1)}\\
(c_0c_1/pq,c_2,c_3,\dots c_{m+n+1})&\in {\cal C}_{(m-1)n},
\end{align}
so long as the generic point on the corresponding manifolds gives
well-defined integrals; in other words, so long as none of the $c_i$ are of
the form $p^iq^j$, $i,j>0$ or $p^{-i}q^{-j}$, $i,j\le 0$.  Indeed, if we
use Lemma \ref{lem:bb_reduction} to take the limit $t_2\to pq/t_0$ in the
generic identity corresponding to $(c_0,c_1,\dots c_{m+n+1})$, we find that
in the left-hand side, the gamma factors corresponding to $t_2$ and $t_0$
cancel, while on the right-hand side, the residue formula gives an
$n-1$-dimensional integral; the result is the generic identity
corresponding to $(c_0c_1/pq,\dots c_{m+n+1})$.  The other combination
follows symmetrically.

Thus, starting with the point $(q,q,\dots q,p,p,\dots p)\in C_{NN}$ for $N$
sufficiently large, we can combine the $q$'s with each other to obtain an
arbitrary collection of values of the form $q^{j+1} p^{-k}$ with $j,k\ge
0$, and similarly combine the $p$'s to values of the form $p^{j+1} q^{-k}$,
subject only to the global condition that their product is $(pq)^{m+1}$.
In other words (taking $N\to \infty$), the theorem holds for a dense set of
points, and thus holds in general.
\end{proof}

\section{The $A_n\leftrightarrow A_m$ transformation}

Consider the following family of $A_n$-type integrals:
\begin{align}
I^{(m)}_{A_n}(Z|t_0,\dots t_{m+n+1};&u_0,\dots u_{m+n+1};p,q)\\
&:=
\frac{(p;p)^n(q;q)^n}{(n+1)!}
\int_{\prod_{0\le i\le n} z_i=Z}
\frac{
\prod_{0\le i\le n}
\prod_{0\le r<m+n+2} \Gamma(t_r z_i,u_r/z_i;p,q)
}{
\prod_{0\le i<j\le n} \Gamma(z_i/z_j,z_j/z_i;p,q)
}
\prod_{1\le i\le n} \frac{dz_i}{2\pi\sqrt{-1}z_i}
.\notag
\end{align}
If $|u_r|<|Z|^{1/(n+1)}<1/|t_r|$, we may take the contour to be the torus of
radius $|Z|^{1/(n+1)}$; outside this range, we must choose the contour to
meromorphically continue the integral.  Such contour considerations can be
greatly simplified by multiplying by a test function $f(Z)$ holomorphic on 
$\C^*$ and integrating over $Z$.  In the resulting integral, the correct
contour has the form $C^{n+1}$, where $C$ contains all points of the form
$p^j q^k u_r$, $j,k\ge 0$, $0\le r\le m+n+1$ and excludes all points of
the form $p^{-j}q^{-k}/t_r$.

Note that unlike the $BC_n$ case, the $A_n$ integral is not equal to 1 for
$n=0$; instead, we pick up the value of the integrand at $Z$:
\[
I^{(m)}_{A_0}(Z|t_0,\dots t_{m+1};u_0,\dots u_{m+1};p,q)
=
\prod_{0\le r<m+2} \Gamma(t_r Z,u_r/Z;p,q)
\]
We also observe that the $Z$ parameter is not a true degree of freedom;
indeed:
\[
I^{(m)}_{A_n}(c^{n+1} Z|\dots t_i\dots;\dots u_i\dots;p,q)
=
I^{(m)}_{A_n}(Z|\dots c t_i\dots;\dots c^{-1} u_i\dots;p,q).
\]
In particular, we could in principle always take $Z=1$ (in which case it
will be omitted), although this is sometimes notationally inconvenient.

\begin{thm}\label{thm:bb_An}
For otherwise generic parameters satisfying $\prod_{0\le r<m+n+2}
t_ru_r=(pq)^{m+1}$,
\[
I^{(m)}_{A_n}(Z|\dots t_i\dots;\dots u_i\dots;p,q)
=
\prod_{0\le r,s<m+n+2} \Gamma(t_r u_s;p,q)
I^{(n)}_{A_m}(Z|\dots T^{\frac{1}{m+1}}/t_i\dots;
                \dots U^{\frac{1}{m+1}}/u_i\dots;p,q),
\]
where $T=\prod_{0\le r<m+n+2}t_r$,
$U=\prod_{0\le r<m+n+2}u_r$.
\end{thm}

\begin{rem}
It appears that this can be viewed as an integral analogue of a series
transformation of Rosengren \cite{RosengrenH:2003} and Kajihara and Noumi
\cite{KajiharaY/NoumiM:2003}, in that the latter should be derivable via
residue calculus from the former.
\end{rem}

For $m=0$, we obtain the following integral conjectured by Spiridonov
\cite{SpiridonovVP:2003}:

\begin{cor}\label{cor:jackson_An}
For otherwise generic parameters satisfying $\prod_{0\le r<n+2} t_ru_r=pq$,
\[
I^{(0)}_{A_n}(\dots t_i\dots;\dots u_i\dots;p,q)
=
\prod_{0\le r,s<n+2} \Gamma(t_r u_s;p,q)
\prod_{0\le r<n+2} \Gamma(U/u_r,T/t_r;p,q)
\]
\end{cor}

The main difficulty with applying the $BC_n$ approach in this case is the
fact that the variables are coupled by the condition $\prod_i z_i=Z$; in
general the integral over this domain of the usual sort of product of
determinants will not be expressible as a determinant of univariate
integrals.  Another difficulty is that, in any event, even in the
``right'' specialization, the integrand is not quite expressible as a
product of determinants.  As we shall see, it turns out that these
problems effectively cancel each other.

In particular, we note the extra factor in the following determinant
identity.

\begin{lem} \cite{FrobeniusG:1882}
\[
\det_{0\le i,j<n}(\frac{\theta(t x_i y_j;p)}{\theta(t,x_i y_j;p)})
=
\frac{\theta(t \prod_{0\le i<n} x_i y_i;p)}{\theta(t;p)}
\prod_{0\le i<j<n} x_jy_j \theta(x_i/x_j,y_i/y_j;p)
\prod_{0\le i,j<n} \theta(x_i y_j;p)^{-1}
\]
\end{lem}

\begin{proof}
Consider the function
\[
F(t;\dots x_i\dots;\dots y_i\dots)
:=
\prod_{0\le i<j<n} (x_jy_j)^{-1} \theta(x_i/x_j,y_i/y_j;p)^{-1}
\prod_{0\le i,j<n} \theta(x_i y_j;p)
\det_{0\le i,j<n}(\frac{\theta(t x_i y_j;p)}{\theta(t,x_i y_j;p)}).
\]
This is clearly holomorphic on $(\C^*)^{2n}$ for $t$ fixed; moreover, since
\[
F(t;p x_0,x_1,\dots x_{n-1};\dots y_i\dots)
=
-(t\prod_{0\le j<n} (x_jy_j))^{-1}
F(t;x_0,x_1,\dots x_{n-1};\dots y_i\dots)
\]
we conclude that $F$ vanishes if $t\prod_{0\le j<n} (x_jy_j)=1$; indeed,
$F(x_0)$ is a degree one theta function, and thus uniquely determined by
its multiplier.  Thus the function
\[
\theta(t \prod_{0\le i<n} x_i y_i;p)^{-1}
F(t;x_0,x_1,\dots x_{n-1};\dots y_i\dots)
\]
is still holomorphic, and indeed we verify that it is an abelian function
of all variables except $t$, so is in fact a function of $t$ alone.  
The remaining factors can thus be recovered from the limiting case:
\[
\lim_{\substack{y_i\to x_i^{-1}\\i=0\dots n-1}}
F(t;x_0,x_1,\dots x_{n-1};\dots y_i\dots)
=
1.
\]
\end{proof}

\begin{rem}
A presumably related application of this determinant to hypergeometric
series identities can be found in \cite{KajiharaY/NoumiM:2003}.
\end{rem}

\begin{lem}
The theorem holds for the special case
\[
I^{(n-1)}_{A_{n-1}}(Z|\dots x_i\dots,\dots \frac{q}{y_i}\dots;\dots
\frac{p}{x_i}\dots,\dots y_i\dots;p,q).
\]
\end{lem}

\begin{proof}
We first observe that the integral:
\[
\int
\frac{\theta(s x z;p)}{\theta(s,x z;p)}
\frac{\theta(t y/z;q)}{\theta(t,y/z;q)}
\frac{dz}{2\pi \sqrt{-1} z}
\]
is symmetric in $x$ and $y$, as follows from the change of variable
$z\mapsto yz/x$.  It thus follows that the determinant
\[
\det_{0\le i,j<n}\Bigl(
\int
\frac{\theta(s x_i z;p)}{\theta(s,x_i z;p)}
\frac{\theta(t y_j/z;q)}{\theta(t,y_j/z;q)}
\frac{dz}{2\pi\sqrt{-1} z}
\Bigr)
\]
is invariant under exchanging the roles of the $x$ and $y$ variables.  As
before, we can write this as a multiple integral of a product of two
determinants:
\begin{align}
n!\det_{0\le i,j<n}\Bigl(
\int
\frac{\theta(s x_i z;p)}{\theta(s,x_i z;p)}&
\frac{\theta(t y_j/z;q)}{\theta(t,y_j/z;q)}
\frac{dz}{2\pi\sqrt{-1} z}
\Bigr)\notag\\
&=
\int
\det_{0\le i,j<n}(\frac{\theta(s x_i z_j;p)}{\theta(s,x_i z_j;p)})
\det_{0\le i,j<n}(\frac{\theta(t y_i/z_j;q)}{\theta(t,y_i/z_j;q)})
\prod_{0\le i<n} \frac{dz_i}{2\pi\sqrt{-1} z_i}\\
&=
\prod_{0\le i<j<n} x_jy_j \theta(x_i/x_j;p)\theta(y_i/y_j;q)\notag\\
&\phantom{{}={}}\int
\frac{\theta(s X Z;p)\theta(t Y/Z;q)}{\theta(s;p)\theta(t;q)}
\frac{
\prod_{0\le i<j<n} \theta(z_i/z_j;p)\theta(z_j/z_i;q)}
{
\prod_{0\le i,j<n} \theta(x_i z_j;p)\theta(y_i/z_j;q)}
\prod_{0\le i<n} \frac{dz_i}{2\pi\sqrt{-1} z_i},
\end{align}
where $X = \prod_i x_i$, $Y = \prod_i y_i$, $Z=\prod_i z_i$.
We thus conclude:
\begin{align}
\int
&
\frac{\theta(s X Z;p)\theta(t Y/Z;q)}{\theta(s;p)\theta(t;q)}
\frac{\prod_{0\le i<j<n} \theta(z_i/z_j;p)\theta(z_j/z_i;q)}
{\prod_{0\le i,j<n} \theta(x_i z_j;p)\theta(y_i/z_j;q)}
\prod_{0\le i<n} \frac{dz_i}{2\pi\sqrt{-1} z_i}\\
&=
\prod_{0\le i<j<n} 
  \frac{\theta(x_i/x_j;q)\theta(y_i/y_j;p)}
       {\theta(x_i/x_j;p)\theta(y_i/y_j;q)}
\int
\frac{\theta(s Y Z;p)\theta(t X/Z;q)}{\theta(s;p)\theta(t;q)}
\frac{\prod_{0\le i<j<n} \theta(z_i/z_j;p)\theta(z_j/z_i;q)}
{\prod_{0\le i,j<n} \theta(y_i z_j;p)\theta(x_i/z_j;q)}
\prod_{0\le i<n} \frac{dz_i}{2\pi\sqrt{-1} z_i}.\notag
\end{align}
Now, if we replace $s$ in this identity by $p^k s$, we find:
\begin{align}
\int
&
\frac{\theta(s X Z;p)\theta(t Y/Z;q)}{(XZ)^k\theta(s;p)\theta(t;q)}
\frac{\prod_{0\le i<j<n} \theta(z_i/z_j;p)\theta(z_j/z_i;q)}
{\prod_{0\le i,j<n} \theta(x_i z_j;p)\theta(y_i/z_j;q)}
\prod_{0\le i<n} \frac{dz_i}{2\pi\sqrt{-1} z_i}\\
&=
\prod_{0\le i<j<n} 
  \frac{\theta(x_i/x_j;q)\theta(y_i/y_j;p)}
       {\theta(x_i/x_j;p)\theta(y_i/y_j;q)}
\int
\frac{\theta(s Y Z;p)\theta(t X/Z;q)}{(YZ)^k\theta(s;p)\theta(t;q)}
\frac{\prod_{0\le i<j<n} \theta(z_i/z_j;p)\theta(z_j/z_i;q)}
{\prod_{0\le i,j<n} \theta(y_i z_j;p)\theta(x_i/z_j;q)}
\prod_{0\le i<n} \frac{dz_i}{2\pi\sqrt{-1} z_i}.\notag
\end{align}
As this is true for all integers $k$, we find that
\begin{align}
\int
f(X&Z)
\frac{\prod_{0\le i<j<n} \theta(z_i/z_j;p)\theta(z_j/z_i;q)}
{\prod_{0\le i,j<n} \theta(x_i z_j;p)\theta(y_i/z_j;q)}
\prod_{0\le i<n} \frac{dz_i}{2\pi\sqrt{-1} z_i}\\
&=
\prod_{0\le i<j<n} 
  \frac{\theta(x_i/x_j;q)\theta(y_i/y_j;p)}
       {\theta(x_i/x_j;p)\theta(y_i/y_j;q)}
\int
f(YZ)
\frac{\prod_{0\le i<j<n} \theta(z_i/z_j;p)\theta(z_j/z_i;q)}
{\prod_{0\le i,j<n} \theta(y_i z_j;p)\theta(x_i/z_j;q)}
\prod_{0\le i<n} \frac{dz_i}{2\pi\sqrt{-1} z_i}.\notag
\end{align}
for any function $f$ holomorphic in a neighborhood of the contour (the
dependence on $s$ and $t$ having been absorbed in $f$).  
But this implies
\begin{align}
\int_{\prod_{0\le i<n} z_i = Z}
&\frac{\prod_{0\le i<j<n} \theta(z_i/z_j;p)\theta(z_j/z_i;q)}
{\prod_{0\le i,j<n} \theta(x_i z_j;p)\theta(y_i/z_j;q)}
\prod_{1\le i<n} \frac{dz_i}{2\pi\sqrt{-1} z_i}\\
&\!\!\!\!=
\prod_{0\le i<j<n} 
  \frac{\theta(x_i/x_j;q)\theta(y_i/y_j;p)}
       {\theta(x_i/x_j;p)\theta(y_i/y_j;q)}
\int_{\prod_{0\le i<n} z_i = ZX/Y}
\frac{
\prod_{0\le i<j<n} \theta(z_i/z_j;p)\theta(z_j/z_i;q)}
{\prod_{0\le i,j<n} \theta(y_i z_j;p)\theta(x_i/z_j;q)^{-1}}
\prod_{1\le i<n} \frac{dz_i}{2\pi\sqrt{-1} z_i}.\notag
\end{align}
Applying the change of variables $z_i\to (X/Y)^{1/n} z_i$ on the right gives
the desired result.
\end{proof}

We also have the following analogue of Lemma \ref{lem:bb_reduction}, with
essentially the same proof.

\begin{lem}
We have the limit:
\[
\lim_{u_0\to t_0^{-1}}
\frac{I^{(m)}_{A_n}(Z|t_0\dots t_{m+n+1};u_0\dots u_{m+n+1};p,q)}
{\Gamma(t_0u_0;p,q)\prod_{0<r<m+n+2} \Gamma(t_0 u_r,t_r u_0;p,q)}
=
I^{(m)}_{A_{n-1}}(t_0 Z|t_1\dots t_{m+n+1};u_1\dots u_{m+n+1};p,q).
\]
\end{lem}

Theorem \ref{thm:bb_An} follows as in the proof of Theorem
\ref{thm:big_bailey}, except that in the definition of ${\cal C}_{mn}$, we
take $t_i u_i = c_i$; we have
\[
(q,q,\dots q,p,p,\dots p)\in {\cal C}_{nn},
\]
and if $(c_0\dots c_{m+n+1})\in {\cal C}_{mn}$, then
\begin{align}
(c_0c_1,c_2,\dots c_{m+n+1})&\in {\cal C}_{m(n-1)}\\
(c_0c_1/pq,c_2,\dots c_{m+n+1})&\in {\cal C}_{m(n-1)}
\end{align}
as long as both sides of the corresponding identities are generically
well-defined.  As before, this shows that ${\cal C}_{mn}$ is dense, and
thus the Theorem \ref{thm:bb_An} holds in general.

\section{Mixed transformations}

Consider the integral associated to $A_1$.  If we eliminate $z_2$ from the
integral using the relation $z_1z_2=1$, we find that the result is
invariant under $z_1\mapsto z_1^{-1}$, and is thus an instance of the
$BC_1$ integral.  Indeed, if
\[
\prod_i t_i u_i = (pq)^{m+1},
\]
then
\begin{align}
I^{(m)}_{A_1}(t_0\dots t_{m+2};u_0\dots u_{m+2};p,q)
&=
\frac{(p;p)(q;q)}{2}
\int
\frac{\Gamma(t_r z^{\pm 1},u_r z^{\pm 1};p,q)}
     {\Gamma(z^{\pm 2};p,q)}
\frac{dz}{2\pi\sqrt{-1} z}\\
&=
I^{(m)}_{BC_1}(t_0\dots t_{m+2},u_0\dots u_{m+2};p,q)
\end{align}

As a consequence, we obtain an identity between the $m=1$ integrals of
types $A_n$ and $BC_n$.

\begin{cor}
If $\prod_{0\le i\le n+2} t_i u_i = (pq)^2$, then
\[
I^{(1)}_{A_n}(\dots t_i\dots;\dots u_i\dots;p,q)
=
\prod_{0\le i<j\le n+2} \Gamma(T/t_it_j,U/u_iu_j;p,q)
I^{(1)}_{BC_n}(\dots (U/T)^{1/4} t_i\dots,
              \dots (T/U)^{1/4} u_i\dots;p,q),
\]
where $T=\prod_{0\le i\le n+2} t_i$, $U=\prod_{0\le i\le n+2} u_i$.
\end{cor}

In particular, since the $BC_n$ integral is symmetric in its $2n+6$
parameters, we obtain an $S_{2n+6}$ symmetry of $I^{(1)}_{A_n}$.  We thus
obtain a total of $n+4$ essentially different transformations of the $A_n$
integral, corresponding to the $n+4$ double cosets of $S_{n+3}\times
S_{n+3}$ in $S_{2n+6}$.

\begin{cor}
Let $k$ be an integer $0\le k\le n+3$.  Then
\begin{align}
I^{(1)}_{A_n}(t_0,\dots,t_{n+2};&
              u_0,\dots,u_{n+2};p,q)\notag\\
&=
\prod_{\substack{0\le r<k\\k\le s\le n+2}} \Gamma(t_r u_s,t_s
              u_r,T/t_rt_s,U/u_ru_s;p,q)
I^{(1)}_{A_n}(t'_0,\dots,t'_{n+2};
              u'_0,\dots,u'_{n+2};p,q),
\end{align}
where
\begin{align}
t'_r &= \begin{cases}
(T/U)^{(n+1-k)/2(n+1)} (T_k/U_k)^{1/(n+1)} u_r,&0\le r<k\\
(U/T)^{k/2(n+1)} (T_k/U_k)^{1/(n+1)} t_r,&k\le r\le n+2
\end{cases}\\
u'_r &= \begin{cases}
(U/T)^{(n+1-k)/2(n+1)} (U_k/T_k)^{1/(n+1)} t_r,&0\le r<k\\
(T/U)^{k/2(n+1)} (U_k/T_k)^{1/(n+1)} u_r,&k\le r\le n+2
\end{cases}\\
T &= \prod_{0\le r\le n+2} t_r\\
U &= \prod_{0\le r\le n+2} t_r\\
T_k &= \prod_{0\le r<k} t_r\\
U_k &= \prod_{0\le r<k} t_r.
\end{align}
\end{cor}

For $k=0$, we obtain the identity transformation, while for $k=n+3$, we
simply switch the $t_i$ and $u_i$ parameters (corresponding to taking
$z\mapsto 1/z$ in the integral).  The case $k=1$ was stated as equation
(6.11) of \cite{SpiridonovVP:2003} (conditional on Corollary
\ref{cor:jackson_An}).  Again, apparently related series identities are
known; see \cite{RosengrenH:2003} and \cite{KajiharaY/NoumiM:2003}.

\section{Difference operators}

The following identity was originally conjectured by van Diejen and
Spiridonov \cite{vanDiejenJF/SpiridonovVP:2001} (their ``Type II'' integral):

\begin{thm}\label{thm:normalization}
For otherwise generic parameters satisfying $|p|,|q|,|t|<1$ and
$t^{2n-2}\prod_{0\le r\le 5} t_r = pq$,
\begin{align}
\frac{(p;p)^n(q;q)^n\Gamma(t;p,q)^n}{2^n n!}
\int_{C^n}
\prod_{1\le i<j\le n}
\frac{\Gamma(t z_i^{\pm 1}z_j^{\pm 1};p,q)}
     {\Gamma(z_i^{\pm 1}z_j^{\pm 1};p,q)}
&
\prod_{1\le i\le n}
\frac{\prod_{0\le r\le 5} \Gamma(t_r z_i^{\pm 1};p,q)
}
{\Gamma(z_i^{\pm 2};p,q)}
\frac{dz_i}{2\pi \sqrt{-1} z_i}\label{eq:normalization}\\
&\qquad=
\prod_{0\le j<n}
\Gamma(t^{j+1};p,q)
\prod_{0\le r<s\le 5}\Gamma(t^j t_rt_s;p,q),
\notag
\end{align}
where the contour $C=C^{-1}$ contains all points of the form $p^i q^j t_r$
for $i,j\ge 0$, excludes their reciprocals, and contains the contours $p^i
q^j t C$ for $i,j\ge 0$.  (In particular, if $|t_r|<1$ for $0\le r\le 5$,
$C$ may be taken to be the unit circle.)
\end{thm}

\begin{proof} \cite{vanDiejenJF/SpiridonovVP:2001}
Suppose $t^{2n} \prod_{0\le r\le 5} t_r = pq$, and consider the 
double integral
\begin{align}
\int_{C^{n+1}} \int_{C^{\prime n}}&
\frac{
\prod_{\substack{0\le i\le n\\1\le j\le n}} \Gamma(\sqrt{t} x_i^{\pm 1} y_j^{\pm 1};p,q)
}{\prod_{0\le i<j\le n} \Gamma(x_i^{\pm 1}x_j^{\pm 1};p,q)
  \prod_{1\le i<j\le n} \Gamma(y_i^{\pm 1}y_j^{\pm 1};p,q)}\\
&
\prod_{0\le i\le n}
\frac{\Gamma(t^n t_0 x_i^{\pm 1};p,q)\prod_{1\le r\le 5} \Gamma(t_r x_i^{\pm 1};p,q)}{\Gamma(x_i^{\pm 2};p,q)}
\prod_{1\le i\le n}
\frac{\Gamma(pq t^{-n-1/2} y_i^{\pm 1}/t_0,t^{-1/2} t_0 y_i^{\pm
    1};p,q)}{\Gamma(y_i^{\pm 2};p,q)}\notag\\
&\prod_{1\le i\le n} \frac{dy_i}{2\pi\sqrt{-1}y_i}
\prod_{0\le i\le n} \frac{dx_i}{2\pi\sqrt{-1}x_i}
\notag
\end{align}
Both the $x$ and $y$ integrals can be evaluated via Corollary
\ref{cor:big_jackson}; comparing both sides gives a recurrence for the
left-hand side of \eqref{eq:normalization}, the unique solution of which
is the right-hand side, as required.
\end{proof}

We will discuss this proof (of which Anderson's proof of the Selberg
integral is a limiting case) in greater detail in the sequel; for the
moment, however, it will be instructive to consider a different proof.  The
main ingredient in the alternate proof is the following identity:

\begin{lem}\label{lem:diffeq_norm}
Let $n$ be a nonnegative integer, and let $u_0$, $u_1$, $u_2$, $u_3$, $t$
satisfy $t^{n-1}u_0u_1u_2u_3=p$.  Then
\begin{align}\label{eq:diffeq_norm}
\sum_{\sigma\in \{\pm 1\}^n}
\prod_{1\le i\le n}
\frac{\prod_{0\le r\le 3} \theta(u_r z_i^{\sigma_i};p)}
     {\theta(z_i^{2\sigma_i};p)}
\prod_{1\le i<j\le n}
\frac{\theta(t z_i^{\sigma_i} z_j^{\sigma_j};p)}
     {\theta(z_i^{\sigma_i} z_j^{\sigma_j};p)}
&=
\prod_{0\le i<n} \theta(t^i u_0u_1,t^i u_0u_2,t^i u_0u_3;p)\\
&=
\prod_{0\le i<n} \theta(t^i u_0u_1,t^i u_0u_2,t^i u_1u_2;p)
\end{align}
\end{lem}

\begin{proof}
We first observe that the condition on the $u_r$ ensures that every term in
the above sum is invariant under all translations $z_i\to pz_i$, and thus
the same is true of their sum.  Moreover, the sum is manifestly invariant
under permutations of the $z_i$ as well as reflections $z_i\to 1/z_i$.
Thus if we multiply the sum by
\[
\prod_{1\le i\le n} z_i^{-1} \theta(z_i^2;p)
\prod_{1\le i<j\le n} z_i^{-1} \theta(z_i z_j,z_iz_j^{-1};p),
\]
the result is a (holomorphic) theta function anti-invariant under the same
group.  But any such theta function is a multiple of the above product; it
thus follows that the desired sum has no singularities in $z_i$, and must
therefore be independent of $z_i$.

To evaluate the sum, we may therefore specialize $z_i = u_0 t^{n-i}$, in
which case all but one of the terms in the sum vanish, so the sum is given
by the remaining term (that with $\sigma_i=1$ for all $i$):
\[
\frac{\prod_{0\le r\le 3} \theta(u_0 u_r t^{n-i};p)}{\theta(u_0^2 t^{2n-2i};p)}
\prod_{1\le i<j\le n}
\frac{\theta(u_0^2 t^{2n+1-i-j};p)}
     {\theta(u_0^2 t^{2n-i-j};p)}
\]
The factors involving $u_0^2$ cancel, and we are thus left with the
evaluation claimed above.
\end{proof}

\begin{proof} (of Theorem \ref{thm:normalization})
Divide the integral by the claimed evaluation, and consider the result as
a meromorphic function on the set $t^{2n-2}t_0t_1t_2t_3t_4t_5=pq$.  We
claim that this function is invariant under the translations
\begin{align}
(t_0,t_1,t_2,t_3,t_4,t_5)\to
(p^{1/2}t_0,p^{1/2}t_1,p^{1/2}t_2,p^{-1/2}t_3,p^{-1/2}t_4,p^{-1/2}t_5)\\
(t_0,t_1,t_2,t_3,t_4,t_5)\to
(q^{1/2}t_0,q^{1/2}t_1,q^{1/2}t_2,q^{-1/2}t_3,q^{-1/2}t_4,q^{-1/2}t_5),
\end{align}
and all permutations thereof.  It will then follow that the ratio is a
constant; to evalute the constant, we may then consider the limit $t_1\to
t^{1-n} t_0^{-1}$ as in Lemma \ref{lem:bb_reduction} above.  (In other
words, we apply the special case of the residue formula of van Diejen and
Spiridonov in which the resulting sum consists of precisely one term.)

Since both sides are symmetric in $p$ and $q$, it suffices to consider the
$q$ translation.  If we factor the integrand as
\[
\Delta^{(n)}(z_1,z_2,\dots z_n)
\Delta^{(n)}(z_1^{-1},z_2^{-1},\dots z_n^{-1})
\prod_{1\le i\le n} \frac{dz_i}{2\pi\sqrt{-1}z_i},
\]
where
\[
\Delta^{(n)}(z_1,z_2,\dots z_n)
=
\prod_{1\le i\le n}
\frac{\Gamma(t_0 z_i,t_1 z_i,t_2 z_i,t_3 z_i,t_4 z_i,t_5 z_i,
             pz_i/(t^{n-1}t_0t_1t_2);p,q)
      }
     {\Gamma(z_i^2,p/(z_it^{n-1}t_0t_1t_2);p,q)}
\prod_{1\le i<j\le n} \frac{\Gamma(t z_iz_j^{\pm 1};p,q)}{\Gamma(z_i
  z_j^{\pm 1};p,q)},
\]
and similarly let $\tilde\Delta^{(n)}$ be the corresponding product with
parameters
\[
(q^{-1/2}t_3,q^{-1/2}t_4,q^{-1/2}t_5,q^{1/2}t_0,q^{1/2}t_1,q^{1/2}t_2)
\]
(permuting the parameters to make the transformation an involution), then
we find that
\[
\frac{\tilde{\Delta}^{(n)}(\dots q^{1/2}z_i\dots)}
     {\Delta^{(n)}(\dots z_i \dots)}
=
\prod_{1\le i\le n}
\frac{\theta(t_0 z_i,t_1 z_i,t_2 z_i,p z_i/t^{n-1} t_0t_1t_2;p)}
     {\theta(z_i^2;p,q)}
\prod_{1\le i<j\le n} \frac{\theta(t z_i z_j;p)}
       {\theta(z_i z_j;p)},
\]
and thus
\[\label{eq:norm_a}
\sum_{\sigma_i\in \{\pm 1\}^n}
\frac{\tilde{\Delta}^{(n)}(\dots q^{1/2}z_i^{\sigma_i}\dots)}
     {\Delta^{(n)}(\dots z_i^{\sigma_i} \dots)}
=
\prod_{0\le i<n}
\theta(t^i t_0 t_1,t^i t_0 t_2,t^i t_1 t_2;p)
\]
by Lemma \ref{lem:diffeq_norm}.  Similarly,
\[\label{eq:norm_b}
\sum_{\sigma_i\in \{\pm 1\}^n}
\frac{\Delta^{(n)}(\dots q^{1/2}z_i^{\sigma_i}\dots)}
     {\tilde{\Delta}^{(n)}(\dots z_i^{\sigma_i} \dots)}
=
\prod_{0\le i<n}
\theta(t^i t_3 t_4/q,t^i t_3 t_5/q,t^i t_4 t_5/q;p).
\]

Now, consider the integral:
\[
\int_C
\tilde{\Delta}^{(n)}(\dots q^{1/2}z_i\dots)
\Delta^{(n)}(\dots z_i^{-1}\dots)
\prod_{1\le i\le n} \frac{dz_i}{2\pi\sqrt{-1}z_i},
\]
where the contour is chosen to contain the points $p^i q^j t_r$ for $i,j\ge
0$, exclude their reciprocals, and contain the contours $tC$ and $tC^{-1}$;
here we note that the poles of $\tilde{\Delta}^{(n)}(\dots
q^{1/2}z_i\dots)$ are a subset of the poles of $\Delta^{(n)}(\dots
z_i\dots)$, so this constraint on the contour is still reasonable.  If we
then perform the change of variable $z_i\mapsto q^{-1/2}/z_i$, we find that
the new contour is legal for the transformed parameters.  In other words,
we have
\begin{align}
\int_C
\tilde{\Delta}^{(n)}(\dots q^{1/2}z_i\dots)
\Delta^{(n)}(\dots z_i^{-1}\dots)
&\prod_{1\le i\le n} \frac{dz_i}{2\pi\sqrt{-1}z_i}
\\
&=\int_{C'}
\Delta^{(n)}(\dots q^{1/2}z_i\dots)
\tilde{\Delta}^{(n)}(\dots z_i^{-1}\dots)
\prod_{1\le i\le n} \frac{dz_i}{2\pi\sqrt{-1}z_i}.\notag
\end{align}
Since the constraints on the contours are symmetrical under $z_i\mapsto
1/z_i$, we may symmetrize the integrands, losing the same factor of $2^n$
on both sides.  The theorem follows upon applying equations
\eqref{eq:norm_a} and \eqref{eq:norm_b} to simplify the symmetrized
integrands.
\end{proof}

\begin{rem}
One can also prove Corollary \ref{cor:big_jackson} by a similar argument,
based on the straightforward identity
\[
\sum_{\sigma\in\{\pm 1\}^n}
\frac{\theta(\prod_{0\le r<n+2} t_r/\prod_{1\le i\le n} z_i^{\sigma_i};p)
      \prod_{1\le i\le n}
        \prod_{0\le r<n+2} \theta(t_r z_i^{\sigma_i};p)}
     {\prod_{1\le i\le j\le n} \theta(z_i^{\sigma_i}z_j^{\sigma_j};p)}
=
\prod_{0\le r<s<n+2} \theta(t_rt_s;p).
\]
\end{rem}

Define a $q$-difference operator $D^{(n)}_{q}(u_0,u_1,u_2,u_3;t,p)$ by
setting
\[
(D^{(n)}_{q}(u_0,u_1,u_2,u_3;t,p)f)(\dots z_i\dots)
:=
\!\!\!
\sum_{\sigma\in \{\pm 1\}^n}
\prod_{1\le i\le n}
\frac{\prod_{0\le r\le 3} \theta(u_r z_i^{\sigma_i};p)}
     {\theta(z_i^{2\sigma_i};p)}
\prod_{1\le i<j\le n}
\frac{\theta(t z_i^{\sigma_i} z_j^{\sigma_j};p)}
     {\theta(z_i^{\sigma_i} z_j^{\sigma_j};p)}
f(\dots q^{\sigma_i/2} z_i\dots).
\]
Thus Lemma \ref{lem:diffeq_norm} gives a formula for the image of $1$ under
$D^{(n)}_{q}(u_0,u_1,u_2,u_3;t,p)$ when $t^{n-1}u_0u_1u_2u_3=p$.  Moreover,
the resulting proof of Theorem \ref{thm:normalization} would appear to be
based on an adjointness relation between two such difference operators, as
we will confirm below.

To make this precise, we need some suitable spaces of functions on which to
act.  Let $A^{(n)}(u_0;p,q)$ be the space of $BC_n$-symmetric $p$-abelian
functions $f$ such that
\[
\prod_{1\le i\le n} \theta(pq z_i^{\pm 1}/u_0;p,q)_{0,m} f(\dots z_i\dots)
\]
is holomorphic for sufficiently large $m$; that is, $f$ is smooth except at
the points $p^k u_0/q^l$, $p^k q^l/u_0$ for $k\in \Z$, $1\le l\le m$, where
it has at most simple poles.  The canonical (multiplication) map from the
tensor product of $A^{(n)}(u_0;p,q)$ and $A^{(n)}(u_0;q,p)$ to the space of
meromorphic functions on $(\C^*)^n$ is generically injective; denote the
image by ${\cal A}^{(n)}(u_0;p,q)$.  In particular, we observe that if
$f\in {\cal A}^{(n)}(u_0;p,q)$, then
\[
\prod_{1\le i\le n} \theta(pq z_i^{\pm 1}/u_0;p,q)_{l,m}
f(\dots z_i\dots)
\propto
\prod_{1\le i\le n}
\frac{\Gamma(u_0 z_i^{\pm 1};p,q)}{\Gamma(p^{-l} q^{-m} u_0 z_i^{\pm 1};p,q)}
f(\dots z_i\dots)
\]
is holomorphic for sufficiently large $l$, $m$.

\begin{rem}
Our main motivation for considering the large space ${\cal
A}^{(n)}(u_0;p,q)$, rather than the smaller spaces in which the functions
are actually abelian, is that such product functions already appear in the
family of univariate biorthogonal functions considered by Spiridonov
\cite[Appendix A]{SpiridonovVP:2003}.
\end{rem}

We now define
\[
{\cal D}^{(n)}_{q}(u_0,u_1,u_2;t,p)f
:=
\frac{
D^{(n)}_{q}(u_0,u_1,u_2,t^{1-n}p/u_0u_1u_2;t,p)f}
{\prod_{1\le i\le n} \theta(t^{n-i} u_0 u_1,t^{n-i} u_0 u_2,t^{n-i} u_1u_2;p)}
\]
We will also need a shift operator $T^{(n)}_{\omega,q}$:
\[
(T^{(n)}_{\omega,q}f)(\dots z_i\dots)=f(\dots q^{1/2} z_i\dots).
\]
Note that this maps $BC_n$-symmetric $q$-abelian functions to
$BC_n$-symmetric $q$-abelian functions.

\begin{lem}\label{lem:diffeq}
The operator ${\cal D}^{(n)}_{q}(u_0,u_1,u_2;t,p)$ induces a linear
transformation
\[
{\cal D}^{(n)}_{q}(u_0,u_1,u_2;t,p):
{\cal A}^{(n)}(\sqrt{q}u_0;p,q)\to {\cal A}^{(n)}(u_0;p,q).
\]
Moreover, the corresponding map
\[
{\cal D}^{(n)}_{q}(u_0,u_1,u_2;t,p):
A^{(n)}(\sqrt{q}u_0;p,q)
\otimes
A^{(n)}(\sqrt{q}u_0;q,p)
\to
A^{(n)}(u_0;p,q)
\otimes
A^{(n)}(u_0;q,p)
\]
can be decomposed as
\[
{\cal D}^{(n)}_{q}(u_0,u_1,u_2;t,p)
=
{\cal D}^{(n)}_{q}(u_0,u_1,u_2;t,p)
\otimes
T^{(n)}_{\omega,q}.
\]
\end{lem}

\begin{proof}
Let
\begin{align}
g&\in A^{(n)}(\sqrt{q}u_0;p,q)\notag\\
h&\in A^{(n)}(\sqrt{q}u_0;q,p).\notag
\end{align}
A straightforward computation, using the fact that $h$ is $q$-abelian,
gives:
\[
{\cal D}^{(n)}_{q}(u_0,u_1,u_2;t,p)(gh)
=
({\cal D}^{(n)}_{q}(u_0,u_1,u_2;t,p)g)
(T^{(n)}_{\omega,q}h)
\]
as required.  That $T^{(n)}_{\omega,q}h\in A^{(n)}(u_0;q,p)$ is
straightforward; that
\[
{\cal D}^{(n)}_{q}(u_0,u_1,u_2;t,p)g
\]
is $p$-abelian follows as in the proof of Lemma \ref{lem:diffeq_norm}.
Finally, we observe that this function is holomorphic at $z_i=u_0$, as
required.
\end{proof}

The desired adjointness relation can then be stated as follows.  For
parameters satisfying $t^{2n-2} u_0u_1t_0t_1t_2t_3 = pq$, define a scalar
product between ${\cal A}^{(n)}(u_0;p,q)$ and ${\cal A}^{(n)}(u_1;p,q)$ as
follows:
\begin{align}
\langle f,g\rangle_{t_0,t_1,t_2,t_3;u_0,u_1;t,p,q}
&:=
\frac{1}{Z}
\int_{C^n}
f(\dots z_i\dots)g(\dots z_i\dots)
\prod_{1\le i<j\le n}
\frac{\Gamma(t z_i^{\pm 1}z_j^{\pm 1};p,q)}
     {\Gamma(z_i^{\pm 1}z_j^{\pm 1};p,q)}\\
&\phantom{{}:=\frac{1}{Z}\int_{C^n}{}}
\prod_{1\le i\le n}
\frac{\prod_{0\le r\le 5} \Gamma(t_r z_i^{\pm 1};p,q)
}
{\Gamma(z_i^{\pm 2};p,q)}
\frac{dz_i}{2\pi \sqrt{-1} z_i},\notag
\end{align}
where
\begin{align}
Z &= 
\frac{2^n n!}{(p;p)^n(q;q)^n\Gamma(t;p,q)^n}
\prod_{1\le i\le n}
\Gamma(t^i;p,q)
\prod_{0\le r<s\le 5} \Gamma(t^{n-i} t_r t_s;p,q),\\
t_4 &= u_0\\
t_5 &= u_1,
\end{align}
and the contour is chosen as in Theorem \ref{thm:normalization}, except
that we first absorb the singularities of $f$ and $g$ into the factors
$\Gamma(u_r z_i^{\pm 1};p,q)$ of the integrand.  In particular, we have
\[
\langle 1,1\rangle_{t_0,t_1,t_2,t_3;u_0,u_1;t,p,q}=1.
\]

\begin{thm}
If $f\in {\cal A}^{(n)}(q^{1/2} u_0;p,q)$,
   $g\in {\cal A}^{(n)}(u_1;p,q)$ and
$t^{2n-2}u_0u_1t_0t_1t_2t_3=pq$, then
\[
\langle {\cal D}^{(n)}_{q}(u_0,t_0,t_1;t,p) f,
g\rangle_{t_0,t_1,t_2,t_3;u_0,u_1;t,p,q}
=
\langle f,{\cal D}^{(n)}_{q}(u'_1,t'_2,t'_3;t,p)
g\rangle_{t'_0,t'_1,t'_2,t'_3;u'_0,u'_1;t,p,q},
\]
where
\[
(t'_0,t'_1,t'_2,t'_3,u'_0,u'_1)=(q^{1/2}t_0,q^{1/2}t_1,q^{-1/2}t_2,q^{-1/2}t_3,q^{1/2}u_0,q^{-1/2}u_1).
\]
\end{thm}

\begin{proof}
The second proof of Theorem \ref{thm:normalization} applies, essentially
without change.
\end{proof}

To understand the significance of this result, we need to introduce a
filtration of the space ${\cal A}^{(n)}(u_0;p,q)$.  Let $\Lambda_n$ be the
set of partitions of at most $n$ parts, and let $\subset$ denote the
inclusion partial order; we also let $\subset$ denote the product partial
order on $\Lambda_n\times \Lambda_n$.  Then for any pair of partitions
$\lambda,\mu\in \Lambda_n$, we define
\[
{\cal A}^{(n)}_{\lambda\mu}(u_0;t;p,q)
\]
to be the subspace of ${\cal A}^{(n)}(u_0;p,q)$ consisting of functions
$f$ such that whenever $(\kappa,\nu)\not\subset (\lambda,\mu)$, we have the
limit
\[
\lim_{\substack{z_i\to p^{-\kappa_i} q^{-\nu_i} t^{i-1} u_0\\i=1\dots n}}
\prod_{1\le i\le n} \theta(pq z_i^{\pm 1}/u_0;p,q)_{l,m}
f(\dots z_i\dots)
=
0
\]
whenever
\[
\prod_{1\le i\le n} \theta(pq z_i^{\pm 1}/u_0;p,q)_{l,m}
f(\dots z_i\dots)
\]
is holomorphic.  Note that enlarging $l$ or $m$ multiplies the equation by
a (possibly zero) scalar, so we really have only one equation for each pair
$(\kappa,\nu)$.

\begin{rem}
In the univariate case
(\cite{WilsonJA:1991,RahmanM:1986,RahmanM/SuslovSK:1993,SpiridonovVP/ZhedanovAS:2000a,SpiridonovVP/ZhedanovAS:2000b,SpiridonovVP:2003}),
this filtration simply corresponds to a sequence of allowed poles.  Given
the role played by vanishing conditions in the theory of Koornwinder
polynomials \cite{OkounkovA:1998a,bcpoly}, it would seem to be natural to
generalize the forbiddance of a pole to the vanishing (after clearing the
denominator) at an appropriate point, thus obtaining our filtration.
\end{rem}

\begin{lem}
For generic $u_0$, $p$, $q$, $t$, the filtration ${\cal
  A}^{(n)}_{\blambda}(u_0;t;p,q)$ is tight in the sense that
\[
\dim {\cal A}^{(n)}_{\blambda}(u_0;t;p,q)
=
1
+
\dim \sum_{\bkappa\subsetneq\blambda}
{\cal A}^{(n)}_{\bkappa}(u_0;t;p,q),
\]
for any partition pair $\blambda\in\Lambda_n^2$.
In particular, each space in the filtration is finite-dimensional.
\end{lem}

\begin{proof}
Let $\blambda=(\lambda,\mu)$, $\bkappa=(\kappa,\nu)$.  Since the spaces
\[
{\cal A}^{(n)}_{\lambda\mu}(u_0;t;p,q)
\]
and
\[
\sum_{(\kappa,\nu)\subsetneq(\lambda,\mu)}
{\cal A}^{(n)}_{\kappa\nu}(u_0;t;p,q)
\]
differ by a single equation, their dimensions differ by at most 1; it thus
suffices to construct a function in the former but not in the latter.

Define a function $F^{(n)}_{\lambda\mu}(u_0:\dots z_i\dots)$ by the following
product:
\[
F^{(n)}_{\lambda\mu}(u_0:\dots z_i\dots;t;p,q)
=
\prod_{\substack{1\le i\le n\\1\le j\le \lambda_1}}
\frac{\theta(p^{j}q t^{-\lambda'_j} z_i^{\pm 1}/u_0;q)}
     {\theta(p^{j}q z_i^{\pm 1}/u_0;q)}
\prod_{\substack{1\le i\le n\\1\le j\le \mu_1}}
\frac{\theta(pq^{j} t^{-\mu'_j} z_i^{\pm 1}/u_0;p)}
     {\theta(pq^{j} z_i^{\pm 1}/u_0;p)}
\]
It follows as in the proof of Lemma 6.3 of \cite{OkounkovA:1998a} that
\[
F^{(n)}_{\lambda\mu}(u_0:\dots z_i\dots)
\in
{\cal A}^{(n)}_{\lambda\mu}(u_0;t;p,q);
\]
on the other hand, we find that
\[
\lim_{\substack{z_i\to p^{-\lambda_i} q^{-\mu_i} t^{i-1} u_0\\i=1\dots n}}
\prod_{1\le i\le n} \theta(pq z_i^{\pm 1}/u_0;p,q)_{l,m}
F^{(n)}_{\lambda\mu}(u_0:\dots z_i\dots;t;p,q)
\]
is generically nonzero.
\end{proof}

\begin{rem}
The function $F^{(n)}_{\lambda\mu}$ is a special case of the interpolation
functions introduced below (Definition \ref{defn:interpolation}).  Indeed,
one can show that
\[
F^{(n)}_{\blambda}(u_0{:};t,p,q)
=
{\cal R}^{*(n)}_{\blambda}(;pqt^{-n}/u_0,u_0;t;p,q).
\]
The existence of such a factorizable special case of the interpolation
functions will turn out to be crucial to the arguments of \cite{bctheta}.
\end{rem}

The reason we have introduced this filtration is the following fact:

\begin{lem}
The difference operator ${\cal D}^{(n)}_q(u_0,t_0,t_1;t,p)$ is triangular
with respect to the above filtration; that is, for all
$\blambda\in\Lambda_n^2$,
\[
{\cal D}^{(n)}_q(u_0,t_0,t_1;t,p)
{\cal A}^{(n)}_{\blambda}(\sqrt{q}u_0;t;p,q)
\subset
{\cal A}^{(n)}_{\blambda}(u_0;t;p,q),
\]
with equality for generic values of the parameters.
\end{lem}

\begin{proof}
Let $\blambda=(\lambda,\mu)$.  Choose $l\ge \lambda_1$, $m\ge \mu_1$, and
consider a function
\[
f\in {\cal A}^{(n)}_{\lambda\mu}(\sqrt{q}u_0;t;p,q)
\]
For $\kappa\subset l^n$, $\nu\subset m^n$, define
\begin{align}
C_{\kappa\nu}(f)
&=
\lim_{\substack{z_i\to p^{-\kappa_i} q^{1/2-\nu_i} t^{i-1} u_0\\i=1\dots n}}
\prod_{1\le i\le n} \theta(pq^{1/2} z_i^{\pm 1}/u_0;p,q)_{l,m}
f(\dots z_i\dots)\\
C'_{\kappa\nu}(f)
&=
\lim_{\substack{z_i\to p^{-\kappa_i} q^{-\nu_i} t^{i-1} u_0\\i=1\dots n}}
\prod_{1\le i\le n} \theta(pq z_i^{\pm 1}/u_0;p,q)_{l,m}
({\cal D}^{(n)}_q(u_0,t_0,t_1;t,p)f)(\dots z_i\dots)
\end{align}
We claim that we can write
\[
C'_{\kappa\nu}
=
\sum_{\nu\subset\rho}
c_{\kappa\nu\rho} C_{\kappa\rho},
\]
where the coefficients $c_{\kappa\nu\rho}$ are meromorphic and
independent of the choice of $f$.  Indeed, this follows readily from the
definition of ${\cal D}$; compare the proof of Theorem 3.2 of \cite{bcpoly}.
More precisely, we see that a given term of the corresponding sum
involves the specialization
\[
\lim_{\substack{z_i\to p^{-\kappa_i} q^{\sigma_i/2-\nu_i} t^{i-1} u_0\\i=1\dots n}}
\prod_{1\le i\le n} (pq^{1/2} z_i^{\pm 1}/u_0;p,q)_{l,m}
f(\dots z_i\dots);
\]
if the sequence $\frac{1-\sigma_i}{2}+\nu_i$ does not induce a partition,
then the remaining factors vanish, while if it does give a partition, that
partition necessarily contains $\nu$.  We also find that the diagonal
coefficient $c_{\kappa\nu\nu}$ is generically nonzero; the result follows.
\end{proof}

Now, given a pair of spaces with corresponding tight filtrations, equipped
with a (sufficiently general) scalar product, there is a unique (up to
scalar multiples) orthogonal pair of bases compatible with the filtration.
In the case of the above scalar product, this suggests the following
definition.

\begin{sq_defn}
For all partition pairs $\blambda\in\Lambda_n^2$, the function
\[
R^{(n)}_{\blambda}(\dots z_i\dots;t_0,t_1,t_2,t_3;u_0,u_1;t;p,q)
\]
is defined to be the unique (up to scalar multiples) element of ${\cal
  A}^{(n)}_{\blambda}(u_0;t;p,q)$ such that
\[
\langle
R^{(n)}_{\blambda}(\dots z_i\dots;t_0,t_1,t_2,t_3;u_0,u_1;t;p,q),
g\rangle_{t_0,t_1,t_2,t_3;u_0;u_1;t;p,q}
=
0
\]
whenever $g\in {\cal A}^{(n)}_{\bkappa}(u_1;t;p,q)$ for some
$\bkappa\subsetneq \blambda$.
\end{sq_defn}

Since our adjoint difference operators preserve the filtrations, they would
necessarily be diagonal in the corresponding bases, {\em if} they were
well-defined.  Unfortunately, we have as yet no reason to believe that the
scalar product is nondegenerate relative to the filtration; that is, that
its restriction to ${\cal A}^{(n)}_{\blambda}(u_0;t;p,q)$ and ${\cal
  A}^{(n)}_{\blambda}(u_1;t;p,q)$ is nondegenerate for all partition pairs.
If this condition were to fail for a given pair
$\bkappa$, then the function $R^{(n)}_{\blambda}$ would not be
uniquely determined for $\blambda\supsetneq\bkappa$, and the argument
breaks down.

There is one special case in which we can prove the scalar product
generically nondegenerate.

\begin{prop}
For generic parameters satisfying $t^{2n-2}t_0t_1t_2t_3u_0u_1=pq$, and any
partition $\lambda\in\Lambda_n$, the scalar product
$\langle\rangle_{t_0,t_1,t_2,t_3;u_0,u_1;t;p,q}$ is nondegenerate between
${\cal A}^{(n)}_{0\lambda}(u_0;t;p,q)$ and ${\cal
  A}^{(n)}_{0\lambda}(u_1;t;p,q)$.
\end{prop}

\begin{proof}
To show a scalar product generically nondegenerate, it suffices to exhibit
a nondegenerate specialization.  Choose $l$ such that the spaces
\[
\theta(pq/u_0;p,q)_{0,l} {\cal A}^{(n)}_{0\lambda}(u_0;t;p,q)\quad
\text{and}\quad \theta(pq/u_1;p,q)_{0,l} {\cal A}^{(n)}_{0\lambda}(u_1;t;p,q)
\]
consist of holomorphic functions, and specialize the parameters so that
every parameter except $u_0$, $u_1$ is real, between $0$ and $1$, while
$u_0$ and $u_1$ are complex conjugates satisfying $0<|u_0|=|u_1|<q^l$.
(This is possible as long as $p<q^{2l-1}t^{2n-2}$.)  Then the contour in
the scalar product can be taken to be the unit torus, on which the weight
function is clearly strictly positive.  Moreover, the filtrations with
respect to $u_0$ and $u_1$ are conjugate to each other.  The scalar
product thus becomes a positive definite Hermitian inner product, and is
therefore nondegenerate.
\end{proof}

This in particular proves the existence and uniqueness of the above
biorthogonal functions, as long as one of the partitions is trivial.  In
general, however, it is unclear how to construct a manifestly nondegenerate
instance of the scalar product.  We will therefore give a more direct
construction of these functions, and by computing their scalar products
show that this problem generically does not arise.  (In addition, the above
construction gives functions that are only guaranteed to be orthogonal when
the corresponding pairs of partitions are distinct but comparable; it will
follow below (as one would expect) that comparability is not necessary.)

To do this, we need a different adjoint pair of difference operators.

First, define
\[
{\cal D}^{-(n)}_q(u_0;t,p)
:=
D^{(n)}_q(u_0,q u_0,p/u_0,\frac{1}{t^{n-1}u_0 q};t,p)
%
\]
Next, define
\begin{align}
({\cal D}^{+(n)}_q&(u_0{:}u_1{:}u_2,u_3,u_4;t,p)f)(\dots z_i\dots)\notag\\
&:=
\prod_{1\le i\le n}
\frac{\theta(pq t^{n-i} u_1/u_0;p)}
{\prod_{2\le r\le 5} \theta(u_r t^{n-i} u_1;p)}\\
&\phantom{{}:={}}
\sum_{\sigma\in \{\pm 1\}^n}
\prod_{1\le i\le n}
\frac{\prod_{1\le r\le 5} \theta(u_r z_i^{\sigma_i};p)}
     {\theta(pq z_i^{\sigma_i}/u_0;p)\theta(z_i^{2\sigma_i};p)}
\prod_{1\le i<j\le n}
\frac{\theta(t z_i^{\sigma_i} z_j^{\sigma_j};p)}
     {\theta(z_i^{\sigma_i} z_j^{\sigma_j};p)}
f(\dots q^{\sigma_i/2} z_i\dots),\notag
\end{align}
where $u_5=p^2q/t^{n-1}u_0u_1u_2u_3u_4$.  Note that aside from the
normalization factor, ${\cal D}^{+(n)}_q$ is symmetric in $u_1$ through $u_5$.

These act as lowering and raising operators with respect to the filtration:

\begin{lem}
For all $\blambda\in\Lambda_n^2$ with $(0,1)^n\subset\blambda$,
\[
{\cal D}^{-(n)}_q(u_0;t,p)
{\cal A}^{(n)}_{\blambda}(q^{3/2} u_0;t;p,q)
\subset
{\cal A}^{(n)}_{\blambda-(0,1)^n}(u_0;t;p,q)
\]
Similarly, for all $\blambda\in\Lambda_n^2$,
\[
{\cal D}^{+(n)}_q(u_0{:}u_1{:}u_2,u_3,u_4;t,p)
{\cal A}^{(n)}_{\blambda}(q^{-1/2} u_0;t;p,q)
\subset
{\cal A}^{(n)}_{\blambda+(0,1)^n}(u_0;t;p,q).
\]
Moreover, the restriction of ${\cal D}^{-(n)}$ is generically surjective,
while the restriction of ${\cal D}^{+(n)}$ is generically injective.
\end{lem}

\begin{proof}
As above.
\end{proof}

\begin{thm}
If $f\in {\cal A}^{(n)}(q^{-1/2} u_0;p,q)$,
   $g\in {\cal A}^{(n)}(u_1;p,q)$ and
$t^{2n-2}u_0u_1t_0t_1t_2t_3=pq$, then
\[
\langle {\cal D}^{+(n)}_{q}(u_0{:}t_0{:}t_1,t_2,t_3;t,p) f,
g\rangle_{t_0,t_1,t_2,t_3;u_0,u_1;t,p,q}
=
C
\langle f,{\cal D}^{-(n)}_{q}(u'_1;t,p) g
\rangle_{t'_0,t'_1,t'_2,t'_3;u'_0,u'_1;t,p,q},
\]
where
\[
(t'_0,t'_1,t'_2,t'_3,u'_0,u'_1)=
(q^{1/2}t_0,q^{1/2}t_1,q^{1/2}t_2,q^{1/2}t_3,q^{-1/2}u_0,q^{-3/2}u_1)
\]
and
\[
C = 
\prod_{1\le i\le n}
\frac{\theta(t^{n-i} t_1 t_2,t^{n-i} t_1 t_3,t^{n-i} t_2 t_3,pq t^{n-i}
t_0/u_0;p)}
{\theta(t^{n-i} t_0 u_1/q,t^{n-i} t_1 u_1/q,t^{n-i} t_2 u_1/q,t^{n-i} t_3 u_1/q,t^{n-i} u_0 u_1/q,t^{n-i} u_0 u_1/q^2,p t_0u_1 t^{2n-1-i};p)}
\]
\end{thm}

\section{Integral operators}

Just as our second proof of Theorem \ref{thm:normalization} is related to
an adjoint pair of difference operators, the argument of van Diejen and
Spiridonov is related to an adjoint pair of {\it integral} operators.
To understand these operators, we first need to understand what happens to
the $I^{(0)}_{BC_n}$ integral when the integrand is multiplied by an element
of ${\cal A}(t_0;p,q)$.  We define a corresponding integral operator as
follows.

\begin{defn}
If $f\in {\cal A}(u_0;p,q)$, then $I^{(n)}(u_0;p,q)f$ is the
function on the set $\prod_{0\le r\le 2n+3} u_r=pq$ defined by
\begin{align}
({\cal I}^{*(n)}(u_0;p,q)f)&(u_1,\dots u_{2n+3})\notag\\
&=
\frac{(p;p)^n(q;q)^n}{2^n n!\prod_{0\le r<s\le 2n+3} \Gamma(u_r u_s;p,q)}
\int_{C^n}
f(\dots z_i\dots)
\prod_{1\le i<j\le n} \frac{1}{\Gamma(z_i^{\pm 1}z_j^{\pm 1};p,q)}\\
&\phantom{
\frac{(p;p)^n(q;q)^n}{2^n n!\prod_{0\le r<s\le 2n+3} \Gamma(u_r u_s;p,q)}
\int_{C^n}\quad
}
\prod_{1\le i\le n}
\frac{
\prod_{0\le r\le 2n+3} \Gamma(u_r z_i^{\pm 1};p,q)}
{\Gamma(z_i^{\pm 2};p,q)}
\frac{dz_i}{2\pi\sqrt{-1}z_i},\notag
\end{align}
with the usual conventions about the choice of contour.
\end{defn}

In particular, by Corollary \ref{cor:big_jackson}, it follows that
\[
{\cal I}^{*(n)}(u_0;p,q)1 = 1.
\]

To determine the action of this integral operator in general, it suffices
to consider $f$ in a spanning set.  We may thus restrict our attention to
functions of the form
\[
f(\dots z_i\dots)
=
\prod_{1\le i\le n}
\frac{
\prod_{1\le j\le l} x_j^{-1}\theta(x_jz_i^{\pm 1};q)
\prod_{1\le j\le m} y_j^{-1}\theta(y_jz_i^{\pm 1};p)}
{\theta(pq z_i^{\pm 1}/u_0;p,q)_{l,m}}.
\]
If we write the theta functions in the numerator as a ratio of elliptic
$\Gamma$ functions, and similarly absorb the denominator factors into a
ratio of elliptic $\Gamma$ functions, we find that the resulting integral
is proportional to an integral of type $I^{(l+m)}_{BC_n}$ in which the extra
$2l+2m$ parameters have pairwise products $p^2q$ and $pq^2$.  If we then
apply Theorem \ref{thm:big_bailey}, we find that the right-hand side
becomes a sum via residue calculus.  We thus obtain the following result.

\begin{thm}\label{thm:int_op_main}
If
\[
f(\dots z_i\dots)
=
\prod_{1\le i\le n}
\frac{
\prod_{1\le j\le l} x_j^{-1}\theta(x_jz_i^{\pm 1};q)
\prod_{1\le j\le m} y_j^{-1}\theta(y_jz_i^{\pm 1};p)}
{\theta(pq z_i^{\pm 1}/u_0;p,q)_{l,m}},
\]
and $\prod_{0\le r\le 2n+3} u_r=pq$, then
\begin{align}
({\cal I}^{*(n)}(u_0;p,q)f)&(u_1\dots u_{2n+3})\notag\\
&=
\prod_{1\le r\le 2n+3} \frac{1}{\theta(pq/u_0u_r;p,q)_{l,m}}\label{eq:int_op_main}\\
&\phantom{{}={}}
\left(
\prod_{1\le i\le l} (1+R(x_i))
\frac{\theta(p^{-l} u_0 x_i;q)\prod_{1\le r\le 2n+3}
  \theta(u_r x_i;q)}{x_i^n \theta(x_i^2;q)}
\prod_{1\le i<j\le l}
\frac{\theta(p x_ix_j;q)}{\theta(x_ix_j;q)}
\right)\notag\\
&\phantom{{}={}}
\left(
\prod_{1\le i\le m} (1+R(y_i))
\frac{\theta(q^{-m} u_0 y_i;p)\prod_{1\le r\le 2n+3}
  \theta(u_r y_i;p)}{y_i^n \theta(y_i^2;p)}
\prod_{1\le i<j\le m}
\frac{\theta(q y_iy_j;p)}{\theta(y_iy_j;p)}
\right),
\notag
\end{align}
where $R(x_k)$ is an operator acting on $g(\dots x_i\dots)$ via the
substitution $x_k\mapsto x_k^{-1}$; thus the factors in parentheses are
sums of $2^l$ and $2^m$ terms respectively.
\end{thm}

Since the factors in parentheses are clearly holomorphic in $u_1\dots
u_{2n+3}$, and the given functions span ${\cal A}(u_0;p,q)$, we obtain
the following as an immediate consequence:

\begin{cor}
If $f\in {\cal A}(u_0;p,q)$ is such that
\[
\theta(pq z_i^{\pm 1}/u_0;p,q)_{l,m} f(\dots z_i\dots)
\]
is holomorphic, then
\[
\prod_{1\le r\le 2n+3} \theta(pq/u_0u_r;p,q)_{l,m}
({\cal I}^{*(n)}(u_0;p,q)f)(u_1\dots u_{2n+3})
\]
is holomorphic on the set $\prod_{0\le r\le 2n+3} u_r = pq$.
\end{cor}

\begin{rems}
Similarly, the left-hand side of \eqref{eq:int_op_main} is manifestly a
holomorphic $q$-theta function in the $x$'s, and a holomorphic $p$-theta
function in the $y$'s; that this is true of the right-hand side follows
from a symmetrization argument analogous to those we have just encountered
in studying difference operators.  And, indeed, the two sums are really
just minor variants of the difference operators we have already seen.
\end{rems}

\begin{rems}\label{rems:int_op_main_alt_proof}
As the above argument is based on Theorem \ref{thm:big_bailey}, it cannot
be directly applied in the limit $p\to 0$.  In fact, one can also derive
this result from Corollary \ref{cor:big_jackson}, for which direct,
non-elliptic, proofs are known in the $p\to 0$ limit \cite{GustafsonRA:1992}.
The basic observation is that if two of the parameters have product $q$,
i.e., if two of the $\Gamma$ factors combine to produce a factor of the form
\[
\prod_{1\le i\le n} \frac{1}{\theta(a z_i^{\pm 1};q)},
\]
then the integrand is essentially invariant under $a\mapsto 1/a$ (aside
from an overall constant).  However, the {\em integral} does not share this
invariance, because inverting $a$ changes the constraint on the contour.
The two contours differ only in whether they contain the points $z=a^{\pm
  1}$; as a result, the difference in the two integrals is (proportional
to) the $n-1$-dimensional integral of the residue at that point.  This
$n-1$-dimensional integral simplifies to the above form, with $l=1$, $m=0$;
the difference of the original $n$-dimensional integrals simplifies to the
desired right-hand side.  This argument can then be repeated as necessary
to prove the theorem for arbitrary values of $l,m\ge 0$.
\end{rems}

\begin{rems}\label{rems:int_op_main_rosengren}
The fact that we obtain an $l+m$-tuple sum is, of course, directly related
to the fact that we needed $2l+2m$ $\Gamma$ factors to represent the
numerator of $f$.  In general, if we took
\[
f(\dots z_i\dots)
=
\prod_{1\le i\le n}
\frac{\prod_{1\le j\le m} \theta(x_j z_i^{\pm 1};p,q)_{a_j,b_j}}
{\theta(pq z_i^{\pm 1}/u_0;p,q)_{\sum a_j,\sum b_j}},
\]
residue calculus would again give a sum, this time a $2m$-tuple sum (i.e.,
the product of an $m$-tuple sum for $p$ and an $m$-tuple sum for $q$).  On
the other hand, we could also compute ${\cal I}^{*(n)}(u_0;p,q)f$ by
specialization of Theorem \ref{thm:int_op_main}, which would give a
sum with $2^{\sum a_j+b_j}$ terms.  The fact that this sum simplifies
underlies Rosengren's arguments in section 7 of \cite{RosengrenH:2004}.
\end{rems}

\begin{rems}
It is particularly striking that the right-hand side factors as a product
of two sums, one involving only $q$-theta functions, and one involving
only $p$-theta functions.  This factorization phenomenon appears to hold
quite generally in the theory of elliptic hypergeometric integrals, but
only when the relevant balancing condition holds.
\end{rems}

Since ${\cal I}^{*(n)}(u_0;p,q)$ takes $BC_n$-symmetric functions to
$A_{2n+2}$-symmetric functions, it is not quite suitable for our purposes.
However, we can readily obtain $BC$-symmetric functions by suitable
specialization.

\begin{defn}
Define operators ${\cal I}^{+(n)}_t(u_0;p,q)$, ${\cal I}^{(n)}_t(u_0{:}u_1,u_2;p,q)$,
and ${\cal I}^{-(n)}_t(u_0{:}u_1,u_2,u_3,u_4;p,q)$ by:
\begin{align}
({\cal I}^{+(n)}_t(u_0;p,q)f)(z_1\dots z_{n+1})
&=
({\cal I}^{*(n)}(u_0;p,q)f)(\frac{t^{-n-1}pq}{u_0},\dots \sqrt{t} z_i^{\pm
    1}\dots)\\
({\cal I}^{(n)}_t(u_0{:}u_1,u_2;p,q)f)(z_1\dots z_n)
&=
({\cal I}^{*(n)}(u_0;p,q)f)(u_1,u_2,
\frac{t^{-n}pq}{u_0u_1u_2},\dots \sqrt{t} z_i^{\pm
    1}\dots)\\
({\cal I}^{-(n)}_t(u_0{:}u_1,u_2,u_3,u_4;p,q)f)(z_1\dots z_{n-1})
&=
({\cal I}^{*(n)}(u_0;p,q)f)(u_1,u_2,u_3,u_4,
\frac{t^{1-n}pq}{u_0u_1u_2u_3u_4},\dots \sqrt{t} z_i^{\pm
    1}\dots).
\end{align}
\end{defn}

\begin{thm}
The above operators are triangular with repsect to the filtration of
${\cal A}^{(n)}(u_0;p,q)$; to be precise,
\begin{align}
{\cal I}^{+(n)}_t(u_0;p,q){\cal A}^{(n)}_{\blambda}(u_0;t;p,q)
&\subseteq
{\cal A}^{(n+1)}_{\blambda}(t^{1/2} u_0;t;p,q)\\
{\cal I}^{(n)}_t(u_0{:}u_1,u_2;p,q){\cal A}^{(n)}_{\blambda}(u_0;t;p,q)
&\subseteq
{\cal A}^{(n)}_{\blambda}(t^{1/2} u_0;t;p,q),
\end{align}
and, if $\blambda_n=(0,0)$,
\[
{\cal I}^{-(n)}_t(u_0{:}u_1,u_2,u_3,u_4;p,q){\cal A}^{(n)}_{\blambda}(u_0;t;p,q)
\subseteq
{\cal A}^{(n-1)}_{\blambda}(t^{1/2} u_0;t;p,q).
\]
Moreover, ${\cal I}^{+(n)}_t(u_0;p,q)$ is generically injective, and ${\cal
  I}^{-(n)}_t(u_0{:}u_1,u_2,u_3,u_4;p,q)$ is generically surjective.
\end{thm}

\begin{proof}
It suffices to consider the action of the operators on the functions
\[
F^{(n)}_{\lambda\mu}(u_0:\dots z_i\dots;t;p,q)
=
\prod_{\substack{1\le i\le n\\1\le j\le \lambda_1}}
\frac{\theta(p^{j}q t^{-\lambda'_j} z_i^{\pm 1}/u_0;q)}
     {\theta(p^{j}q z_i^{\pm 1}/u_0;q)}
\prod_{\substack{1\le i\le n\\1\le j\le \mu_1}}
\frac{\theta(pq^{j} t^{-\mu'_j} z_i^{\pm 1}/u_0;p)}
     {\theta(pq^{j} z_i^{\pm 1}/u_0;p)}
\]
considered above.  Applying Theorem \ref{thm:int_op_main}, we find that
each term of the resulting sum is also of this form, with appropriately
constrained partitions.  The one exception is ${\cal I}^{-(n)}_t$ in the
case when $\lambda_n$ or $\mu_n>0$, which we will consider below.
\end{proof}

As promised, the integral operators indeed satisfy appropriate adjointness
relations.

\begin{thm}
If $f\in {\cal A}^{(n)}(u_0;p,q)$, $g\in {\cal A}^{(n)}(t^{-1/2} u_1;p,q)$ and
$t^{2n-2}u_0u_1t_0t_1t_2t_3=pq$, then
\[
\langle {\cal I}^{(n)}_{t}(u_0{:}t_0,t_1;p,q) f,
g\rangle_{t'_0,t'_1,t'_2,t'_3;u'_0,u'_1;t;p,q}
=
\langle f,
{\cal I}^{(n)}_{t}(u'_1{:}t'_2,t'_3;p,q)
g\rangle_{t_0,t_1,t_2,t_3;u_0,u_1;t;p,q},
\]
where
\[
(t'_0,t'_1,t'_2,t'_3,u'_0,u'_1)=(t^{1/2}t_0,t^{1/2}t_1,t^{-1/2}t_2,t^{-1/2}t_3,t^{1/2}u_0,t^{-1/2}u_1).
\]
Similarly, if $f\in {\cal A}^{(n)}(u_0;p,q)$, $g\in {\cal
  A}^{(n-1)}(u'_1;p,q)$ and $t^{2n-2}u_0u_1t_0t_1t_2t_3=pq$, then
\[
\langle {\cal I}^{-(n)}_{t}(u_0{:}t_0,t_1,t_2,t_3;p,q) f,
g\rangle_{t'_0,t'_1,t'_2,t'_3,u'_0,u'_1;t;p,q}
=
\langle f,{\cal I}^{+(n-1)}_{t}(u'_1;p,q) g
\rangle_{t_0,t_1,t_2,t_3,u_0,u_1;t;p,q},
\]
where
\[
(t'_0,t'_1,t'_2,t'_3,u'_0,u'_1)=(t^{1/2}t_0,t^{1/2}t_1,t^{1/2}t_2,t^{1/2}t_3,t^{1/2}u_0,t^{-1/2}u_1).
\]
\end{thm}

\begin{proof}
In each case, the definition of the integral operators allows us to
express the inner products as double integrals; the stated identities
correspond to changing the order of integration.
\end{proof}

Recall that for the operators ${\cal D}^-$ and ${\cal I}^-$, we were only
able to show triangularity with respect to a portion of the filtration; for
some functions, the methods we used were insufficient to understand the
images.  The key observation for dealing with those cases is that the
difficult case for ${\cal D}^-$ is precisely the (generic) image of ${\cal
  I}^+$, and similarly the difficult case for ${\cal I}^-$ is the image of
${\cal D}^+$.  Thus to complete our understanding of the action of these
operators on the filtration, it will suffice to prove the following result.

\begin{thm}\label{thm:mp_chain}
For any function $f\in {\cal A}(u_0;p,q)$,
\[
{\cal D}^{-(n)}_q(q^{-3/2} t^{1/2} u_0;t,p){\cal I}^{+(n-1)}_t(u_0;p,q)f=0.
\]
Similarly, for any function $f\in {\cal A}^{(n)}(q^{-1/2} u_0;p,q)$,
\[
{\cal I}^{-(n)}_t(u_0{:}t_0,t_1,t_2,t_3;p,q)
{\cal D}^{+(n)}_q(u_0{:}t_0,t_1,t_2,t_3;t,p)
f=0.\label{eq:mp_chain2}
\]
\end{thm}

\begin{proof}
For the first identity, take
\[
f(\dots z_i\dots)
=
\prod_{1\le i\le n}
\frac{
\prod_{1\le j\le l} x_j^{-1}\theta(x_jz_i^{\pm 1};q)
\prod_{1\le j\le m} y_j^{-1}\theta(y_jz_i^{\pm 1};p)}
{\theta(pq z_i^{\pm 1}/u_0;p,q)_{l,m}};
\]
we can thus compute its image via Theorem \ref{thm:int_op_main}
and the definition of ${\cal D}^{-(n)}$.  The vanishing of the resulting
sum follows as a special case of Lemma \ref{lem:vanishI} below.

For the second identity, we can argue as in the proof of adjointness of
the difference operators to express the image as the integral of $f(\dots
z_i\dots)$ with respect to an appropriate $BC_n$-symmetric density.  That
this density vanishes identically follows from Lemma \ref{lem:vanishII}
below.
\end{proof}

\begin{lem}\label{lem:vanishI}
For arbitrary parameters satisfying $vw\prod_{1\le i\le n} q_i^2=1$, and
generic $z_1$,\dots $z_n$,
\[
\prod_{1\le i\le n} (1+R(z_i))
\frac{\theta(v q_i z_i,w q_i z_i;p)}{\theta(z_i^2;p)}
\prod_{1\le i<j\le n}
  \frac{q_j^{-1}\theta(q_i q_j z_i z_j,q_j z_i/q_i z_j;p)}
       {\theta(z_i z_j,z_i/z_j;p)}
=
0.
\]
\end{lem}

\begin{proof}
For $n=1$, the summand is manifestly antisymmetric under $R(z_i)$, and thus
the lemma follows in that case.  Thus assume $n>1$, set $v=u/Q$,
$w=(uQ)^{-1}$ with $Q:=\prod_{1\le i\le n} q_i$, and consider the sum as a
function of $u$.  We readily verify that it is a $BC_1$-symmetric theta
function in $u$ of degree $n$; we thus need only show that it vanishes at
more than $n$ independent points.  If $u=Q z_n/q_n$, the terms involving
$R(z_n)$ vanish; moreover, if we pull out $BC_{n-1}$-symmetric factors, we
obtain a special case of the $n-1$-dimensional sum.  By symmetry, the
identity holds for any point of the form $u=Q z_i^{\pm 1}/q_i$; since $n>1$,
these $2n$ points are generically independent, and the result follows.
\end{proof}

We note the following related result in passing:

\begin{cor}
For arbitrary parameters satisfying $tuvw\prod_{1\le i\le n} q_i^2=1$, and
generic $z_1$,\dots $z_n$,
\[
t\theta(uv,uw,vw;p)
\prod_{1\le i\le n} (1+R(z_i))
\frac{\theta(t z_i/q_i,u q_i z_i,v q_i z_i,w q_i z_i;p)}{z_i\theta(z_i^2;p)}
\prod_{1\le i<j\le n}
  \frac{q_j^{-1}\theta(q_i q_j z_i z_j,q_j z_i/q_i z_j;p)}
       {\theta(z_i z_j,z_i/z_j;p)}
\]
is symmetric under permutations of $t$, $u$, $v$, $w$.
\end{cor}

\begin{proof}
The sum is manifestly symmetric in $u$, $v$, $w$, so it suffices to show
that it is invariant under the exchange of $t$ and $u$.  Thus take the
difference of the given sum and its image upon exchanging $t$ and $u$.  If
we then set $t=q_{n+1}y_{n+1}$, $u=q_{n+1}/y_{n+1}$, we obtain the
$n+1$-dimensional instance of the lemma.
\end{proof}

\begin{lem}\label{lem:vanishII}
For generic values of $y_1$,\dots $y_n$,
\[
\prod_{1\le i\le n} (1+R(y_i))
\prod_{1\le i<j\le n} \frac{\theta(uy_iy_j;p)}{\theta(y_iy_j;p)}
\prod_{1\le i\le n}
\frac{
\prod_{1\le r\le n-1} \theta(u^{-1/2} z_r^{\pm 1} y_i;p)
}
{y_i^{n-2}\theta(y_i^2;p)}
=
0.
\]
\end{lem}

\begin{proof}
When $n=1$, the summand is antisymmetric under $R(y_1)$, and the sum
therefore vanishes.  Now, consider the sum for general $n$ as a function
of $z_{n-1}$.  This is manifestly a $BC_1$-symmetric theta function of
degree $n$; it thus suffices to show that it vanishes at more than $n$
independent points.  If $z_{n-1}=u^{1/2}y_n$, the terms coming from
$R(y_n)$ vanish; we thus obtain an instance of the $n-1$-dimensional sum,
which vanishes by induction.  By symmetry, the sum vanishes at any point
of the form $z_{n-1}=u^{1/2}y_i^{\pm 1}$; this gives $2n$ independent
values at which the sum vanishes, proving the lemma.
\end{proof}

A similar argument applies to the following result, which can also be
obtained from Theorem \ref{thm:big_bailey} via residue calculus.

\begin{thm}
Choose integers $m\ge l\ge 0$, and suppose $q^{m-l}t_0t_1t_2t_3=q$.
Then we have the following identity.
\begin{align}
\prod_{1\le i\le m}
&(1+R(x_i))
\frac{\theta(t_0 x_i,t_1 x_i,t_2 x_i,t_3 x_i;p)
      \prod_{1\le r\le l} \theta(q^{-1/2} x_i y_r^{\pm 1};p)}
     {x_i^{l+1} \theta(x_i^2;p)}
\prod_{1\le i<j\le m}
\frac{\theta(q x_ix_j;p)}{\theta(x_ix_j;p)}\notag\\
&=
\prod_{0\le i<m-l}
\frac{\theta(q^i t_0t_1,q^it_0t_2,q^it_0t_3;p)}{(-q^{1/2})^{m-1}t_0}\\
&\phantom{{}={}}
\prod_{1\le i\le l}
(1+R(y_i))
\frac{\theta(\frac{q^{1/2}}{t_0} y_i,\frac{q^{1/2}}{t_1} y_i,\frac{q^{1/2}}{t_2} y_i,\frac{q^{1/2}}{t_3} y_i;p)
      \prod_{1\le r\le m} \theta(q^{-1/2} y_i x_r^{\pm 1};p)}
     {y_i^{m+1} \theta(y_i^2;p)}
\prod_{1\le i<j\le l}
\frac{\theta(q y_iy_j;p)}{\theta(y_iy_j;p)}
\notag
\end{align}
\end{thm}

\begin{proof}
By the usual symmetry argument, we find that both sides are
$BC_m$-symmetric theta functions of degree $l$ in $x$.  By induction, both
sides agree if $x_m$ is of the form $t_r$ or $q^{-1/2} y_i^{\pm 1}$; this
gives $2l+4$ independent points at which the functions agree, which shows
that they agree everywhere.
\end{proof}

This gives rise to some commutation relations between our difference and
integral operators.

\begin{cor}\label{cor:diff_int_comm}
For any function $f\in {\cal A}^{(n)}(q^{1/2}u_0;p,q)$,
\begin{align}
{\cal I}^{(n)}_t(u_0{:}t_0,t_1;p,q)
{\cal D}^{(n)}_q(u_0,t_0,t_1;t,p)f
&=
{\cal D}^{(n)}_q(t^{1/2}u_0,t^{1/2}t_0,t^{1/2}t_1;t,p)
{\cal I}^{(n)}_t(q^{1/2}u_0{:}q^{1/2}t_0,q^{1/2}t_1;p,q)f
\\
{\cal I}^{(n)}_t(u_0{:}t_0,t_1;p,q)
{\cal D}^{(n)}_q(u_0,t_0,t_2;t,p)f
&=
{\cal D}^{(n)}_q(t^{1/2}u_0,t^{1/2}t_0,t^{-1/2}t_2;t,p)
{\cal I}^{(n)}_t(q^{1/2}u_0{:}q^{1/2}t_0,q^{-1/2}t_1;p,q)f
\\
{\cal I}^{+(n)}_t(u_0;p,q)
{\cal D}^{(n)}_q(u_0,t_0,t_1;t,p)f
&=
{\cal D}^{(n+1)}_q(t^{1/2}u_0,t^{-1/2}t_0,t^{-1/2}t_1;t,p)
{\cal I}^{+(n)}_t(q^{1/2}u_0;p,q)f
\end{align}
while for any function $f\in {\cal A}^{(n)}(q^{-1/2}u_0;p,q)$,
\begin{align}
{\cal I}^{(n)}_t(u_0{:}t_0,t_1;p,q)&
{\cal D}^{+(n)}_q(u_0{:}t_0{:}t_1,t_2,t_3;t,p)f\\
&=
{\cal D}^{+(n)}_q(t^{1/2}u_0{:}t^{1/2}t_0{:}t^{1/2}t_1,t^{-1/2}t_2,t^{-1/2}t_3;t,p)
{\cal I}^{(n)}_t(q^{-1/2}u_0{:}q^{1/2}t_0,q^{1/2}t_1;p,q)f\notag
\end{align}
\end{cor}

\begin{proof}
In each case, arguing as in the proof of adjointness of difference
operators transforms the left-hand side into an integral of $f$ against a
$BC_n$-symmetric density which itself can be transformed via the theorem to
give the right-hand side.
\end{proof}

\section{Biorthogonal functions}

Now that we have suitable difference and integral operators, we are now in
a position to construct the desired biorthogonal functions.

\begin{defn}
For each integer $n\ge 0$, we define a family of functions
\[
\tilde{\cal R}^{(n)}_{\blambda}(z_1\dots z_n;t_0{:}t_1,t_2,t_3;u_0,u_1;t;p,q)
\in {\cal A}^{(n)}_{\blambda}(u_0;t;p,q)
\]
indexed by a partition pair $\blambda$ of length at most $n$ and with
parameters satisfying $t^{2n-2}t_0t_1t_2t_3u_0u_1=pq$, as follows.  For
$n=0$, we take
\[
\tilde{\cal R}^{(0)}(;t_0{:}t_1,t_2,t_3;u_0,u_1;t;p,q):=1.
\]
Otherwise, if $\blambda_n=(0,0)$, we set
\begin{align}
\tilde{\cal R}^{(n)}_{\blambda}&(;t_0{:}t_1,t_2,t_3;u_0,u_1;t;p,q)\\
&:=
{\cal I}^{+(n-1)}_t(t^{-1/2} u_0;p,q)
\tilde{\cal R}^{(n-1)}_{\blambda}(;t^{1/2} t_0{:}t^{1/2} t_1,t^{1/2}
t_2,t^{1/2} t_3;t^{-1/2} u_0,t^{1/2} u_1;t;p,q).\notag
\end{align}
If $(0,1)^n\subset \blambda$, set
\begin{align}
\tilde{\cal R}^{(n)}_{\blambda}&(;t_0{:}t_1,t_2,t_3;u_0,u_1;t;p,q)\\
&:=
{\cal D}^{+(n)}_q(u_0{:}t_0{:}t_1,t_2,t_3;t,p)
\tilde{\cal R}^{(n)}_{\blambda-(0,1)^n}(;q^{1/2} t_0{:}q^{1/2} t_1,q^{1/2} t_2,q^{1/2}
t_3;q^{-1/2} u_0,q^{-3/2} u_1;t;p,q).\notag
\end{align}
Finally, if $(1,0)^n\subset \blambda$, but $(0,1)^n\not\subset\blambda$, set
\begin{align}
\tilde{\cal R}^{(n)}_{\blambda}&(;t_0{:}t_1,t_2,t_3;u_0,u_1;t;p,q)\\
&:=
{\cal D}^{+(n)}_p(u_0{:}t_0{:}t_1,t_2,t_3;t,q)
\tilde{\cal R}^{(n)}_{\blambda-(1,0)^n}
(;p^{1/2} t_0{:}p^{1/2} t_1,p^{1/2} t_2,p^{1/2} t_3;p^{-1/2} u_0,p^{-3/2} u_1;t;p,q).\notag
\end{align}
\end{defn}

\begin{rem}
The above definition closely resembles, and indeed was inspired by,
Okounkov's integral representation for interpolation polynomials
\cite{OkounkovA:1998a}; in fact, in an appropriate limit, our ${\cal
  I}^{+(n)}$ becomes Okounkov's integral operator (which can thus be
expressed as a contour integral, rather than a $q$-integral).
\end{rem}

It is clear that this inductively defines a family of functions as
described; note also that the last relation still holds if
$(1,1)^n\subset\blambda$, since the corresponding $p$- and $q$-difference
operators ``commute''.  In addition, it is clear that these functions
should agree with the functions $R$ we attempted to define above, aside
from the fact that the scalar multiplication freedom has been eliminated:

\begin{prop}\label{prop:biorth_normalization}
The functions $\tilde{\cal R}$ satisfy the normalization condition
\[
\tilde{\cal R}^{(n)}_{\blambda}(\dots t^{n-i}
t_0\dots;t_0{:}t_1,t_2,t_3;u_0,u_1;t;p,q)
=
1.
\]
\end{prop}

Since the ``diagonal'' coefficients of the $+$ operators with respect to
the filtration are generically nonzero, we find that they form a section
of the filtration; that is:

\begin{prop}
For any partition pair $\blambda$, and for generic values of the
parameters, the functions
\[
\tilde{\cal R}^{(n)}_{\bkappa}(;t_0{:}t_1,t_2,t_3;u_0,u_1;t;p,q)
\]
for $\bkappa\subset\blambda$ form a basis of ${\cal
  A}^{(n)}_{\blambda}(t;p,q)$.
\end{prop}

Also, since each of the $+$ operators used above factors as a tensor
product, we find that the same holds for our family of functions.

\begin{lem}
Each function $\tilde{\cal R}^{(n)}_\blambda$ is a product of a $q$-abelian
and a $p$-abelian function; more precisely, we have
\begin{align}
\tilde{\cal R}^{(n)}_{\lambda\mu}(\dots z_i\dots;t_0{:}t_1,t_2,t_3;u_0,u_1;t;p,q)
&=
\tilde{\cal R}^{(n)}_{\lambda 0}(\dots z_i\dots;t_0{:}t_1,t_2,t_3;u_0,u_1;t;p,q)\\
&\phantom{{}={}}
\tilde{\cal R}^{(n)}_{0 \mu}(\dots z_i\dots;t_0{:}t_1,t_2,t_3;u_0,u_1;t;p,q)
.\notag
\end{align}
\end{lem}

Similarly, from adjointness and Theorem \ref{thm:mp_chain}, we
can conclude:

\begin{thm}\label{lem:biorthogonality}
The functions $\tilde{\cal R}^{(n)}_{\blambda}$ satisfy the biorthogonality
relation
\[
\langle
\tilde{\cal R}^{(n)}_{\blambda}(;t_0{:}t_1,t_2,t_3;u_0,u_1;t;p,q),
\tilde{\cal R}^{(n)}_{\bkappa}(;t_0{:}t_1,t_2,t_3;u_1,u_0;t;p,q)
\rangle_{t_0,t_1,t_2,t_3,u_0,u_1;t,p,q}
=0
\]
whenever $\bkappa\ne\blambda$.  In particular, $\tilde{\cal
  R}^{(n)}_{\blambda}$ is orthogonal to the space ${\cal
  A}^{(n)}_{\bkappa}(t;p,q)$ whenever
$\blambda\not\subset\bkappa$.
\end{thm}

\begin{rem}
In particular, it follows that our functions agree with the univariate
biorthogonal functions considered in \cite[Appendix A]{SpiridonovVP:2003}.
Note that in the univariate case, the definition involves only the raising
difference operators; the integral operators are unnecessary.  This gives
rise to a generalized Rodriguez-type formula; compare
\cite{RahmanM/SuslovSK:1993}.
\end{rem}

\begin{thm}\label{thm:biorth_diffeqs}
The functions $\tilde{\cal R}^{(n)}_{\blambda}$ satisfy the difference equations:
\begin{align}
{\cal D}^{(n)}_p(u_0,t_0,t_1;t,q)
\tilde{\cal R}^{(n)}_{\blambda}(;p^{1/2}t_0{:}p^{1/2}t_1,p^{-1/2}t_2,p^{-1/2}t_3;p^{1/2}u_0,p^{-1/2}u_1&;t;p,q)\\
&=
\tilde{\cal R}^{(n)}_{\blambda}(;t_0{:}t_1,t_2,t_3;u_0,u_1;t;p,q)\notag\\
{\cal D}^{(n)}_q(u_0,t_0,t_1;t,p)
\tilde{\cal R}^{(n)}_{\blambda}(;q^{1/2}t_0{:}q^{1/2}t_1,q^{-1/2}t_2,q^{-1/2}t_3;q^{1/2}u_0,q^{-1/2}u_1&;t;p,q)\\
&=
\tilde{\cal R}^{(n)}_{\blambda}(;t_0{:}t_1,t_2,t_3;u_0,u_1;t;p,q)\notag
\end{align}
and the integral equation
\[
{\cal I}^{(n)}_t(u_0{:}t_0,t_1;p,q)
\tilde{\cal R}^{(n)}_{\blambda}(;t_0{:}t_1,t_2,t_3;u_0,u_1;t;p,q)
=
\tilde{\cal R}^{(n)}_{\blambda}(;t^{1/2}t_0{:}t^{1/2}t_1,t^{-1/2}t_2,t^{-1/2}t_3;t^{1/2}u_0,t^{-1/2}u_1;t;p,q).
\]
\end{thm}

\begin{proof}
Since each of the operators respects the factorization of $\tilde{\cal
  R}^{(n)}_{\lambda\mu}$, it suffices to consider the cases $\lambda=0$ or
$\mu=0$, which are clearly equivalent.  In particular, the inner product is
now generically nondegenerate, and thus $\tilde{\cal R}^{(n)}_{\lambda0}$
and $\tilde{\cal R}^{(n)}_{0\mu}$ are uniquely determined by
biorthogonality and the normalization condition.  Since each of the three
operators we are considering has triangular adjoint, the left-hand sides
satisfy biorthogonality; on the other hand, we readily compute that each
operator preserves the normalization condition.
\end{proof}

\begin{rem}
This gives rise to an alternate proof of the commutation relations of
Corollary \ref{cor:diff_int_comm}, by comparing the actions of the two
sides on the appropriate basis of biorthogonal functions.  Similarly, one
obtains the commutation relations:
\begin{align}
{\cal I}^{(n+1)}_t(u_0{:}t_0,t_1;p,q)
{\cal I}^{+(n)}_t(t^{-1/2} u_0;p,q)
&=
{\cal I}^{+(n)}_t(u_0;p,q)
{\cal I}^{(n)}_t(t^{-1/2}u_0{:}t^{1/2}t_0,t^{1/2}t_1;p,q)\\
{\cal I}^{(n)}_t(u_0{:}t_0,t_1;p,q)
{\cal I}^{(n)}_t(t^{-1/2}u_0{:}t^{-1/2}t_0,t^{1/2}t_2;p,q)
&=
{\cal I}^{(n)}_t(u_0{:}t_0,t_2;p,q)
{\cal I}^{(n)}_t(t^{-1/2}u_0{:}t^{-1/2}t_0,t^{1/2}t_1;p,q)\\
{\cal D}^{(n)}_q(u_0,t_0,t_1;t,p)
{\cal D}^{(n)}_q(q^{1/2}u_0,q^{1/2}t_0,q^{-1/2}t_2;t,p)
&=
{\cal D}^{(n)}_q(u_0,t_0,t_2;t,p)
{\cal D}^{(n)}_q(q^{1/2}u_0,q^{1/2}t_0,q^{-1/2}t_1;t,p)\\
{\cal D}^{+(n)}_q(u_0{:}t_0{:}t_1,t_2,t_3;t,p)
{\cal D}^{(n)}_q(q^{-1/2}u_0,q^{1/2}t_0,{}&q^{1/2}t_1;t,p)\notag\\
=
{\cal D}^{(n)}_q(u_0,t_0,t_1;t,p)
&{\cal D}^{+(n)}_q(q^{1/2}u_0{:}q^{1/2}t_0{:}q^{1/2}t_1,q^{-1/2}t_2,q^{-1/2}t_3;t,p)
\end{align}
In contrast to Corollary \ref{cor:diff_int_comm}, it is unclear how to
prove these commutation relations directly.
\end{rem}

\begin{cor}
For any partition $\lambda$,
\[
T^{(n)}_{\omega,p}
\tilde{\cal R}^{(n)}_{0\lambda}(;p^{1/2}t_0{:}p^{1/2}t_1,p^{-1/2}t_2,p^{-1/2}t_3;p^{1/2}u_0,p^{-1/2}u_1;t;p,q)
=
\tilde{\cal R}^{(n)}_{0\lambda}(;t_0{:}t_1,t_2,t_3;u_0,u_1;t;p,q).
\]
Moreover,
\[
\tilde{\cal R}^{(n)}_{0\lambda}(;p^{k_0} t_0{:}p^{k_1} t_1,p^{k_2} t_2,p^{k_3}
t_3;p^{l_0} u_0,p^{l_1} u_1;t;p,q)
=
\tilde{\cal R}^{(n)}_{0\lambda}(;t_0{:}t_1,t_2,t_3;u_0,u_1;t;p,q)
\]
for all choices of integers $k_?$, $l_?$ such that $k_0+k_1+k_2+k_3+l_0+l_1=0$.
\end{cor}

\begin{proof}
The first claim follows from the fact that
\[
{\cal D}^{(n)}_p(u_0,t_0,t_1;t,q)f=T^{(n)}_{\omega,p}f
\]
for any $p$-abelian function $f$.  Now, when $l_0=0$, the second claim
follows from the definition of ${\cal R}^{(n)}_{0\lambda}$ and the fact that
${\cal D}^{+(n)}_q(u_0{:}t_0{:}t_1,t_2,t_3;t,p)$ is a $p$-abelian function
of the $u_?$ and $t_?$ parameters.  Iterating the first claim and using
the fact that $\tilde{\cal R}^{(n)}_{0\lambda}$ is $p$-abelian gives an
instance of the second claim with $l_0=1$, and thus the claim holds in
general.
\end{proof}

To see how the operators act when $t_0$ is not among the parameters of the
operator, we need to determine how $\tilde{\cal R}^{(n)}$ changes when we
switch $t_0$ and $t_1$.  This leaves the biorthogonality relation unchanged,
so multiplies the function by a constant; to determine that constant,
it suffices to compute the following evaluation.

\begin{prop}
For generic values of the parameters,
\[
\tilde{\cal R}^{(n)}_{\blambda}(\dots
t^{n-i}t_1\dots;t_0{:}t_1,t_2,t_3;u_0,u_1;t;p,q)
=
\frac{
\cC^0_{\blambda}(t^{n-1}t_1t_2,t^{n-1}t_1t_3,pqt^{n-1}t_0/u_0,t^{1-n}/t_1u_1;t;p,q)}
{
\cC^0_{\blambda}(t^{n-1}t_0t_2,t^{n-1}t_0t_3,pqt^{n-1}t_1/u_0,t^{1-n}/t_0u_1;t;p,q)}.
\]
\end{prop}

\begin{proof}
This follows by comparing the actions of ${\cal
  D}^{+(n)}_q(u_0{:}t_0{:}t_1,t_2,t_3;t,p)$ and ${\cal
  D}^{+(n)}_q(u_0{:}t_1{:}t_0,t_2,t_3;t,p)$.  This gives a recurrence for
the desired specialization, having the right-hand side as unique solution.
\end{proof}

It will be convenient at this point to introduce ``hatted'' parameters.
These are defined as follows.  First, we have:
\[
\hat{t}_0 = \sqrt{t_0t_1t_2t_3/pq}=\frac{t^{1-n}}{\sqrt{u_0u_1}}.
\]
The remaining parameters are then defined by giving invariants of the
transformation.  To be precise, we define $\hat{t}_1$, $\hat{t}_2$,
$\hat{t}_3$, $\hat{u}_0$, and $\hat{u}_1$ by insisting that
\[
\hat{t}_0\hat{t}_1=t_0t_1\quad
\hat{t}_0\hat{t}_2=t_0t_2\quad
\hat{t}_0\hat{t}_3=t_0t_3\quad
\frac{\hat{u}_0}{\hat{t}_0} = \frac{u_0}{t_0}\quad
\frac{\hat{u}_1}{\hat{t}_0} = \frac{u_1}{t_0}.
\]
Note in particular that
\[
t^{2n-2}\hat{t}_0\hat{t}_1\hat{t}_2\hat{t}_3\hat{u}_0\hat{u}_1=pq.
\]
The action of the hat transformation on the $t$ parameters is, of course,
quite familiar from the theory of Koornwinder polynomials
\cite{KoornwinderTH:1992,SahiS:1999} (aside from the factor of $p$ required
to preserve symmetry); the action on the $u$ parameters is then essentially
forced by the balancing condition.  We furthermore define
$z_i(\blambda;\hat{t}_0):=(p,q)^{\blambda_i}t^{n-i}\hat{t}_0$.

In the following formulas, the ratios of $\Gamma$ functions that appear are
sometimes ill-defined, in that some of the factors vanish.  These should be
interpreted by multiplying the argument of each $\Gamma$ function by the
same scale factor, then taking the limit as that scale factor approaches 1.
Alternatively, it turns out in each case that the ratio can be formally
expressed in terms of theta functions alone, and that upon doing so, the
resulting formula is well-defined.  Similar comments apply to ratios of
$\theta$ functions.  In particular, we note that
\[
\prod_{1\le i\le n}
\frac{\theta(v z_i(\blambda;w)^{\pm 1};p)}{\theta(v z_i(0,0;w)^{\pm 1};p)}
=
\prod_{1\le i\le n}
\frac{\Gamma(q v z_i(\blambda;w)^{\pm 1},v z_i(0,0;w)^{\pm 1};p,q)}{\Gamma(v
  z_i(\blambda;w)^{\pm 1},q v z_i(0,0;w)^{\pm 1};p,q)}
\propto
\frac{{\cal C}^0_{\blambda}(t^{n-1} q v w,t^{n-1} pqw/v;p,q)}
     {{\cal C}^0_{\blambda}(t^{n-1} v w,t^{n-1} pw/v;p,q)},
\]
where the constant of proportionality is independent of $v$.

\begin{cor}\label{cor:biorth:diff_int_eqs}
We have the difference equations
\begin{align}
{\cal D}^{(n)}_q(u_0,t_0,t_1;t,p)
{\cal D}^{(n)}_q(q^{1/2} u_0,q^{-1/2} t_2,q^{-1/2} t_3;t,p)
&\tilde{\cal R}^{(n)}_{\blambda}(;t_0{:}t_1,t_2,t_3;q u_0,q^{-1} u_1;t;p,q)\\
&=
\frac{{\cal E}^{\cal D}_{\blambda}(\hat{t}_1{:}\hat{t}_0;t;p{:}q)}
     {{\cal E}^{\cal D}_{\blambda}(\hat{u}_0{:}\hat{t}_0;t;p{:}q)}
\tilde{\cal R}^{(n)}_{\blambda}(;t_0{:}t_1,t_2,t_3;u_0,u_1;t;p,q),\notag\\
{\cal D}^{(n)}_q(u_0,t_2,t_3;t,p)
{\cal D}^{(n)}_q(q^{1/2} u_0,q^{-1/2} t_0,q^{-1/2} t_1;t,p)
&\tilde{\cal R}^{(n)}_{\blambda}(;t_0{:}t_1,t_2,t_3;q u_0,q^{-1} u_1;t;p,q)\\
&=
\frac{{\cal E}^{\cal D}_{\blambda}(\hat{t}_1/q{:}\hat{t}_0;t;p{:}q)}
     {{\cal E}^{\cal D}_{\blambda}(\hat{u}_0{:}\hat{t}_0;t;p{:}q)}
\tilde{\cal R}^{(n)}_{\blambda}(;t_0{:}t_1,t_2,t_3;u_0,u_1;t;p,q),\notag
\end{align}
where
\[
{\cal E}^{\cal D}_{\blambda}(v{:}w;t;p{:}q)
:=
\prod_{1\le i\le n}
\frac{\theta(v z_i(\blambda;w)^{\pm 1};p)}
     {\theta(v z_i(0,0;w)^{\pm 1};p)}.
\]
Similarly,
\begin{align}
{\cal I}^{(n)}_t(t^{1/2}u_0{:}t^{-1/2}t_2,t^{-1/2}t_3;p,q)
{\cal I}^{(n)}_t(u_0{:}t_0,t_1;p,q)&
\tilde{\cal R}^{(n)}_{\blambda}(;t_0{:}t_1,t_2,t_3;u_0,u_1;p,q)\\
&=
\frac{{\cal E}^{\cal I}_{\blambda}(\hat{t}_1{:}\hat{t}_0;t;p,q)}
     {{\cal E}^{\cal I}_{\blambda}(\hat{u}_0{:}\hat{t}_0;t;p,q)}
\tilde{\cal R}^{(n)}_{\blambda}(;t_0{:}t_1,t_2,t_3;t u_0,u_1/t;p,q)\notag
\\
{\cal I}^{(n)}_t(t^{1/2}u_0{:}t^{-1/2}t_0,t^{-1/2}t_1;p,q)
{\cal I}^{(n)}_t(u_0{:}t_2,t_3;p,q)&
\tilde{\cal R}^{(n)}_{\blambda}(;t_0{:}t_1,t_2,t_3;u_0,u_1;p,q)\\
&=
\frac{{\cal E}^{\cal I}_{\blambda}(\hat{t}_1/t{:}\hat{t}_0;t;p,q)}
     {{\cal E}^{\cal I}_{\blambda}(\hat{u}_0{:}\hat{t}_0;t;p,q)}
\tilde{\cal R}^{(n)}_{\blambda}(;t_0{:}t_1,t_2,t_3;t u_0,u_1/t;p,q)\notag
\end{align}
where
\[
{\cal E}^{\cal I}_{\blambda}(v{:}w;t;p,q)
:=
\prod_{1\le i\le n}
\frac{\Gamma(v z_i(\blambda;w)^{\pm 1},t v z_i(0,0;w)^{\pm
    1};p,q)}{\Gamma(t v z_i(\blambda;w)^{\pm 1},v z_i(0,0;w)^{\pm 1};p,q)}
\]
\end{cor}

The $-$ and $+$ operators give similar equations:

\begin{thm}
\begin{align}
{\cal D}^{-(n)}_q(u_0;t,p)
{\cal D}^{+(n)}_q(q^{3/2} u_0{:}\frac{t_0}{q^{1/2}}{:}\frac{t_1}{q^{1/2}},\frac{t_2}{q^{1/2}},\frac{t_3}{q^{1/2}};t,p)
&\tilde{\cal R}^{(n)}_{\blambda}(;t_0{:}t_1,t_2,t_3;q u_0,q^{-1} u_1;t;p,q)\\
&=
C
\frac{{\cal E}^{\cal D}_{\blambda}(\hat{t}_0/q{:}\hat{t}_0;p)}
     {{\cal E}^{\cal D}_{\blambda}(\hat{u}_0{:}\hat{t}_0;p)}
\tilde{\cal R}^{(n)}_{\blambda}(;t_0{:}t_1,t_2,t_3;u_0,u_1;t;p,q),\notag
\end{align}
where
\[
C =
\prod_{1\le i\le n}
\frac{
\theta(t^{n-i} u_0t_0,t^{n-i} u_0t_1,t^{n-i} u_0t_2,t^{n-i} u_0t_3,p
t^{n-i} t_0/qu_0,pq t^{n-i},t^{n-i}t_0 t_1 t_2 t_3/pq^2;p)}
{
\theta(p t^{2n-i-1} t_0 u_1,t^{n-i} t_0 t_1/q,t^{n-i} t_0 t_2/q,t^{n-i} t_0 t_3/q;p)
}.
\]
Similarly,
\begin{align}
{\cal I}^{-(n+1)}_t(t^{1/2}u_0{:}\frac{t_0}{t^{1/2}},\frac{t_1}{t^{1/2}},\frac{t_2}{t^{1/2}},\frac{t_3}{t^{1/2}};p,q)
{\cal I}^{+(n)}_t(u_0;p,q)&
\tilde{\cal R}^{(n)}_{\blambda}(;t_0{:}t_1,t_2,t_3;u_0,u_1;p,q)\\
&=
\frac{{\cal E}^{\cal I}_{\blambda}(t/\hat{t}_0{:}\hat{t}_0;t;p,q)}
     {{\cal E}^{\cal I}_{\blambda}(\hat{u}_0{:}\hat{t}_0;t;p,q)}
\tilde{\cal R}^{(n)}_{\blambda}(;t_0{:}t_1,t_2,t_3;t u_0,u_1/t;p,q)\notag
\end{align}
\end{thm}

\begin{proof}
In each case, by adjointness, both sides satisfy biorthogonality, and must
therefore be proportional.  To determine the constant of proportionality,
we can compare to one of the corresponding equations from Corollary
\ref{cor:biorth:diff_int_eqs}.  Indeed, the fact of proportionality shows
that the relevant products of difference (or integral) operators differ in
their action only by a diagonal transformation; as a result, we can
compute the ratio of their constants of proportionality using {\em any}
section of the filtration.  In particular, it is straightforward to compute
diagonal coefficients using the sections with which
we proved triangularity in the first place, thus giving the desired result.
\end{proof}

\begin{rem}
We thus find that for $v\in
\{\hat{t}_1,\hat{t}_2,\hat{t}_3,\hat{t}_0/q,\hat{t}_1/q,\hat{t}_2/q,\hat{t}_3/q\}$,
we have a difference operator ${\cal D}(v)$ (of ``order'' 2) such that
\[
{\cal D}(v)\tilde{\cal R}^{(n)}_{\blambda}(;t_0{:}t_1,t_2,t_3;q u_0,q^{-1}
u_1;t;p,q)
=
\frac{{\cal E}^{\cal D}_{\blambda}(v{:}\hat{t}_0;t;p{:}q)}
     {{\cal E}^{\cal D}_{\blambda}(\hat{u}_0{:}\hat{t}_0;t;p{:}q)}
\tilde{\cal R}^{(n)}_{\blambda}(;t_0{:}t_1,t_2,t_3;u_0,u_1;t;p,q);
\]
moreover
\[
{\cal D}^{+(n)}_q(u_0{:}t_0{:}t_1,t_2,t_3;t,p)
{\cal D}^{-(n)}(q^{-1/2}u_0;t,p)
\]
essentially gives us such an operator for $v=\hat{t}_0$.  We conjecture
that such an operator exists for all $v$; since the ``eigenvalue'' is
effectively just a $BC_1$-symmetric theta function of degree n in $v$, this
conjecture certainly holds for $n\le 7$.  Such a collection of difference
operators, together with the various spaces of higher-degree difference
operators obtained by composing them, would seem to give the analogue of
the center of the affine Hecke algebra applicable to our biorthogonal
functions.  Indeed, in the Koornwinder limit, the conjecture certainly
holds, and the resulting space of operators is precisely the subspace of
the center of the affine Hecke algebra having degree at most 1.
\end{rem}

In particular, this gives us a recurrence for the nonzero values of the
inner product.  Define
\[
\Delta^{(n)}(z;t_0,t_1,t_2,t_3,t_4,t_5;t;p,q)
=
\prod_{1\le i<j\le n}
\frac{\Gamma(t z_i^{\pm 1}z_j^{\pm 1};p,q)}
     {\Gamma(z_i^{\pm 1}z_j^{\pm 1};p,q)}
\prod_{1\le i\le n}
\frac{
\prod_{0\le r\le 5} \Gamma(t_r z_i^{\pm 1};p,q)
}
{\Gamma(z_i^{\pm 2};p,q)};
\]
in other words, this is simply the density with respect to which our
functions are biorthogonal.

\begin{thm}\label{thm:biorth_ip}
For any partition pair $\blambda$ of length at most $n$, and for generic
values of the parameters,
\begin{align}
\langle
\tilde{\cal R}^{(n)}_{\blambda}(;t_0{:}t_1,t_2,t_3;u_0,u_1&;t;p,q),
\tilde{\cal R}^{(n)}_{\blambda}(;t_0{:}t_1,t_2,t_3;u_1,u_0;t;p,q)
\rangle_{t_0,t_1,t_2,t_3,u_0,u_1;t,p,q}\\
&=
\frac{\Delta^{(n)}(\dots z_i(0,0;\hat{t}_0)\dots;\hat{t}_0,\hat{t}_1,\hat{t}_2,\hat{t}_3,\hat{u}_0,\hat{u}_1;t;p,q)}
{\Delta^{(n)}(\dots
  z_i(\blambda;\hat{t}_0)\dots;\hat{t}_0,\hat{t}_1,\hat{t}_2,\hat{t}_3,\hat{u}_0,\hat{u}_1;t;p,q)}\notag\\
&=
\Delta_{\blambda}(t^{2n-2} \hat{t}_0^2|t^n,t^{n-1}\hat{t}_0\hat{t}_1,t^{n-1}\hat{t}_0\hat{t}_2,t^{n-1}\hat{t}_0\hat{t}_3,t^{n-1}\hat{t}_0\hat{u}_0,t^{n-1}\hat{t}_0\hat{u}_1;t;p,q)^{-1}
.\\
&=
\Delta_{\blambda}(\frac{1}{u_0u_1}|t^n,t^{n-1}t_0t_1,t^{n-1}t_0t_2,t^{n-1}t_0t_3,\frac{t^{1-n}}{t_0u_0},\frac{t^{1-n}}{t_0u_1};t;p,q)^{-1}
.
\end{align}
\end{thm}

This of course, is the direct analogue of the formula for the inner
products of Koornwinder polynomials.

If $t_0t_1=p^{-l}q^{-m}t^{1-n}$, then the integral converts via residue
calculus to a sum, and we thus obtain the following discrete
biorthogonality property.

\begin{thm}
For any partition pairs $\blambda,\bkappa\subset (l,m)^n$, and
for otherwise generic parameters satisfying $t_0t_1=p^{-l}q^{-m}t^{1-n}$,
\begin{align}
\sum_{\bmu\subset (l,m)^n}
\tilde{\cal R}^{(n)}_{\blambda}(\dots
z_i(\bmu;t_0)\dots;t_0{:}t_1,t_2,t_3;u_0,u_1;t;p,q)
&\tilde{\cal R}^{(n)}_{\bkappa}(\dots
z_i(\bmu;t_0)\dots;t_0{:}t_1,t_2,t_3;u_1,u_0;t;p,q)\\
&\frac{\Delta^{(n)}(\dots z_i(\bmu;t_0)\dots;t_0,t_1,t_2,t_3,u_0,u_1;t;p,q)}
{\Delta^{(n)}(\dots z_i(0,0;t_0)\dots;t_0,t_1,t_2,t_3,u_0,u_1;t;p,q)}
=
0\notag
\end{align}
unless $\blambda=\bkappa$.
\end{thm}

\begin{rem}
Note that when $t_0t_1=p^{-l}q^{-m}t^{1-n}$, we have
\[
z_i(\bmu;t_0)=z_{n+1-i}((l,m)^n-\bmu;t_1)^{-1},
\]
and thus summing over $z_i(\bmu;t_1)$ gives the same result.
\end{rem}

This result leads to a very important special case of the $\tilde{\cal R}$
functions.

\begin{cor}
If $t_1u_1=t^{1-n}$, then
\[
\tilde{\cal R}^{(n)}_{\blambda}(\dots
z_i(\bkappa;t_1)\dots;t_0{:}t_1,t_2,t_3;u_0,u_1;t;p,q)=0
\]
unless $\blambda\subseteq\bkappa$.  Moreover, in this case
$\tilde{\cal R}^{(n)}_{\blambda}$ is independent of $t_2$ and $t_3$, and
up to scalar multiplication, is independent of $t_0$.
\end{cor}

\begin{proof}
First suppose that we have $t_0t_1=p^{-l}q^{-m}t^{1-n}$ for $l,m$ such that
$\blambda,\bkappa\subset (l,m)^n$, and consider the discrete
biorthogonality relation.  We observe that for $f\in {\cal
  A}^{(n)}(u_1;p,q)$ such that
\[
\prod_{1\le i\le n} \theta(pq z_i^{\pm 1}/u_1;p,q)_{l,m} f(\dots z_i\dots)
\]
is holomorphic, and for partitions $\blambda\subset (l,m)^n$,
\begin{align}
\lim_{v\to t^{1-n}/u_1} 
\prod_{1\le i\le n} \theta(pq z_i(\blambda;v)^{\pm 1}/u_1;p,q)_{l,m}
&f(\dots z_i(\blambda;v)\dots)\\
&=
\lim_{\substack{z_i\to (p,q)^{\bkappa_i} t^{i-1} u_1\\i=1\dots n}}
\prod_{1\le i\le n} \theta(pq z_i^{\pm 1}/u_1;p,q)_{l,m}
f(\dots z_i\dots).\notag
\end{align}
In other words, if $f\in {\cal A}^{(n)}_{\bkappa}(u_1;t;p,q)$, then the
inner product of our function with $f$ can be expressed as a sum over
partition pairs contained in $\bkappa$, by the very definition of the
filtration.  The desired vanishing property follows immediately.  Moreover,
this orthogonality is independent of the specific values for $t_2$, $t_3$,
and thus changing $t_2$ or $t_3$ can at most multiply our function by a
scalar; this scalar must then be 1 by the normalization formula.

We thus find that the result holds whenever $t_0$ is of the above form.
Since the given quantity is a product of abelian functions of $t_0$ for any
choice of $\blambda$, $\bkappa$, the fact that it holds for $t_0$ of the
form $p^{-l}q^{-m}t^{1-n}/t_1$ implies that it holds in general.  Symmetry
in $t_0$, $t_2$, $t_3$ then shows that the dependence on $t_0$ is only via
the normalization.
\end{proof}

With this in mind, we consider the following alternate normalization in
this case.

\begin{defn}\label{defn:interpolation}
The {\em interpolation functions} ${\cal
  R}^{*(n)}_{\blambda}(;t_0,u_0;t;p,q)$
are defined by
\[
{\cal R}^{*(n)}_{\blambda}(;t_0,u_0;t;p,q)
=
\Delta^0_{\blambda}(t^{n-1}t_0/u_0|t^{n-1}t_0t_1,t_0/t_1;t;p,q)
\tilde{\cal R}^{(n)}_{\blambda}(;t_1{:}t_0,t_2,t_3;u_0,t^{1-n}/t_0;t;p,q),
\]
where the right-hand side is independent of the choice of $t_1$, $t_2$,
$t_3$, as long as $t^{n-1} t_1t_2t_3u_0 = pq$.
The {\em multivariate elliptic binomial coefficient}
$\binomE{\blambda}{\bkappa}_{[a,b];t;p,q}$ is defined by
\[
\binomE{\blambda}{\bkappa}_{[a,b];t;p,q}
:=
\Delta_\bkappa(\frac{a}{b}|t^n,1/b;t;p,q)
{\cal R}^{*(n)}_{\bkappa}(\dots
z_i(\blambda;t^{1-n}a^{1/2})\dots;t^{1-n}a^{1/2},b/a^{1/2};t;p,q),
\]
for $n\ge \ell(\blambda),\ell(\bkappa)$.
\end{defn}

\begin{rem}
An alternate definition uses the fact that
\[
{\cal R}^{*(n)}_{l^n+\lambda,m^n+\mu}(;t_0,u_0;t;p,q)
=
\frac{(pq/t_0u_0)^{-2l|\mu|-2m|\lambda|}\prod_{1\le i\le n} \theta(t_0 x_i^{\pm 1};p,q)_{l,m}}
     {\prod_{1\le i\le n} \theta((pq/u_0) x_i^{\pm 1};p,q)_{l,m}}
{\cal R}^{*(n)}_{\lambda,\mu}(;p^lq^m t_0,u_0/p^lq^m;t;p,q)
\]
(which follows by two applications of equation \eqref{eq:Rscomp} below)
together with the action of ${\cal I}^{+(n)}$ to obtain an integral
representation generalizing that of \cite{OkounkovA:1998a}.
\end{rem}

We note in particular that
\[
{\cal R}^{*(n)}_{\bkappa}(\dots z_i\dots,a,ta,...t^{m-1}a;a,b;t;p,q),
=
\frac{\Delta_{\bkappa}(\frac{t^{n-1}a}{b}|t^{n-m},1/b;t;p,q)}
     {\Delta_{\bkappa}(\frac{t^{n-1}a}{b}|t^n,1/b;t;p,q)}
{\cal R}^{*(n-m)}_{\bkappa}(\dots z_i\dots;t^m a,b;t;p,q),
\]
and thus the multivariate elliptic binomial coefficient is independent of
$n$ (as long as $\ell(\blambda),\ell(\bkappa)\le n$, that is).

The significance of these interpolation functions is that one can express
connection coefficients for the biorthogonal functions in terms of
multivariate elliptic binomial coefficients.  (The proof requires a more
thorough study of these binomial coefficients, and will thus be deferred to
\cite{bctheta}.)

\begin{thm} \cite{bctheta} \label{thm:biorth_conn}
If we define connection coefficients $c_{\blambda\bmu}$ by
\[
\tilde{\cal R}^{(n)}_{\blambda}(;t_0{:}t_1,t_2,t_3;u_0,u_1;t;p,q)\\
=
\sum_{\bmu\subset\blambda}
c_{\blambda\bmu}
\tilde{\cal R}^{(n)}_{\bmu}(;t_0{:}t_1v,t_2,t_3;u_0,u_1/v;t;p,q),
\]
then
\[
c_{\blambda\bmu}
=
\frac{
\Delta^{0}_{\blambda}(1/u_0u_1|1/v,t^{n-1}t_2t_3,pqt^{n-1}t_0/u_0,t_1v/u_1;t;p,q)}
{
\Delta^{0}_{\bmu}(v/u_0u_1|v,t^{n-1}t_2t_3,pqt^{n-1}t_0/u_0,t_1v/u_1;t;p,q)}
\binomE{\blambda}{\bmu}_{[1/u_0u_1,1/v];t;p,q}
\]
\end{thm}

If the biorthogonal function on the right is specialized to an
interpolation function, we obtain the following generalization of
Okounkov's binomial formula for Koornwinder polynomials \cite{OkounkovA:1998a}:

\begin{cor}
\[
\tilde{\cal R}^{(n)}_{\blambda}(;t_0{:}t_1,t_2,t_3;u_0,u_1;t;p,q)
=
\sum_{\bkappa\subseteq\blambda}
c_{\bkappa}
{\cal R}^{*(n)}_{\bkappa}(\dots z_i(\blambda;\hat{t}_0)\dots;\hat{t}_0,\hat{u}_0;t;p,q)
{\cal R}^{*(n)}_{\bkappa}(;t_0,u_0;t;p,q),
\]
where
\[
c_{\bkappa}=
\Delta_{\bkappa}(t^{n-1}t_0/u_0|t^n,pq/u_0t_1,pq/u_0t_2,pq/u_0t_3;t;p,q)
\]
\end{cor}

Since $c_{\bkappa}$ above remains the same when the parameters are replaced
by their hatted analogues, we obtain the following corollary, the analogue
of the ``evaluation symmetry'' property of Koornwinder polynomials.

\begin{cor}\label{cor:ev_symm}
For any partition pairs $\blambda$, $\bkappa$ of length at most $n$, and
for generic values of the parameters,
\[
\tilde{\cal
  R}^{(n)}_{\blambda}(\dots z_i(\bkappa;t_0)\dots;t_0{:}t_1,t_2,t_3;u_0,u_1;t;p,q)
=
\tilde{\cal
  R}^{(n)}_{\bkappa}(\dots z_i(\blambda;\hat{t}_0)\dots;\hat{t}_0{:}\hat{t}_1,\hat{t}_2,\hat{t}_3;\hat{u}_0,\hat{u}_1;t;p,q).
\]
\end{cor}

Before leaving the topic of biorthogonal functions, it remains to justify
our assertions that these are a generalization of Koornwinder polynomials.
The inner product clearly can be degenerated into the Koornwinder inner
product; the difficulty is the filtration.  Indeed, in order to degenerate
the inner product, we must take $p\to 0$, $u_0\to \{0,\infty\}$, at which
point the definition of the filtration breaks.  It turns out that the
filtration actually does have a well-defined limit; however, we have been
unable to find an argument for this other than as a corollary of the
following result.

\begin{thm}
Fix otherwise generic parameters $t_0$, $t_1$, $t_2$, $t_3$.  Then the
limits
\begin{align}
\lim_{u_0\to 0} \lim_{p\to 0}&
\tilde{\cal R}^{(n)}_{0\lambda}(\dots
z_i\dots;t_0{:}t_1,t_2,t_3;u_0,\frac{pq}{t^{2n-2}u_0t_0t_1t_2t_3};t;p,q)\\
\lim_{u_0\to \infty} \lim_{p\to 0}&
\tilde{\cal R}^{(n)}_{0\lambda}(\dots
z_i\dots;t_0{:}t_1,t_2,t_3;u_0,\frac{pq}{t^{2n-2}u_0t_0t_1t_2t_3};t;p,q)
\end{align}
agree, and give a family of $BC_n$-symmetric Laurent polynomials
orthogonal with respect to the Koornwinder inner product.  Moreover, these
polynomials are diagonal with respect to dominance of monomials, and thus
are precisely the Koornwinder polynomials (normalized to have principal
specialization 1).
\end{thm}

\begin{proof}
The key observation is that, although the definition of the {\it
  filtration} blows up in the limit, the raising difference and integral
operators have perfectly fine limits.  Consequently, the above limits are
indeed well-defined; as the choice $u_0\to 0$ or $u_0\to\infty$ has no
effect on the limiting operators, it can have no effect on the limiting
functions.  Since the space of $BC_n$-symmetric $p$-theta functions of
degree $m$ tends in the limit $p\to 0$ to the space of $BC_n$-symmetric
Laurent polynomials of degree at most $m$ in each variable, our functions
become rational functions in that limit.  Taking the limit $u_0\to
0,\infty$ causes the poles of the rational functions to move to $0$ and
$\infty$, thus giving Laurent polynomials.  Finally, we observe that
because the above limits agree, biorthogonality becomes orthogonality in
the limit.  (Recall that ${\cal R}^{(n)}_{0\lambda}$ is $p$-abelian in its
parameters, so the factor of $p$ in $u_1$ can be moved around arbitrarily
before taking the limit.)  We have thus proved the first claim.

To see that these agree with Koornwinder polynomials, we observe that the
operator ${\cal D}^{(n)}_q(u_0,t_0,t_1;t)$ also has a well-defined limit;
standard arguments (\cite[Theorem 3.2]{bcpoly}) show that the limit is
triangular with respect to dominance of monomials.  It thus follows from
Theorem \ref{thm:biorth_diffeqs} that the limiting polynomials are
eigenfunctions of a pair of triangular difference operators, and thus must
themselves be triangular.  The normalization then follows from Proposition
\ref{prop:biorth_normalization}.
\end{proof}

\begin{rem}
In order to determine the constant of proportionality, i.e., determine the
leading coefficient of the limiting polynomial, we need simply determine
how the raising operators affect the leading coefficient.  For the
difference operator, this is straightforward; for the integral operator,
we can appeal to \ref{thm:int_op_main} and, by using the fact
\[
\sum_{\lambda\subset m^n}
(-1)^{mn-|\lambda|}
m_\lambda(y_1,y_2,\dots y_n)
e_{m^n-\lambda}(z_1,z_2,\dots z_m,1/z_1,1/z_2,\dots 1/z_m)
=
\prod_{\substack{1\le i\le n\\1\le j\le m}}
(y_i+1/y_i-z_j-1/z_j)
\]
(where $m_\lambda$ is a $BC_n$-symmetric monomial), reduce to the
difference operator case.  The result, of course, is simply Macdonald's
``evaluation'' conjecture; Theorem \ref{thm:biorth_ip} then gives the
nonzero values of the inner product.  (For more details, see \cite{ICMS}.)
The remaining (``symmetry'') conjecture does not follow from the methods
given above, however (although there are at least two different arguments
for deducing it from evaluation: \cite{vanDiejenJF:1996},
\cite{OkounkovA:1998a}).  The argument we will give in \cite{bctheta} does
descend to the Koornwinder case; indeed, the result is precisely the proof
given in \cite{bcpoly}.
\end{rem}

It follows from \cite[Theorem 7.25]{bcpoly} that the filtration has the
following limit.

\begin{cor}
Choose an integer $n\ge 0$, and a partition $\lambda$ of at most $n$ parts.
Then the limits $u_0\to 0$ or $u_0\to\infty$ of the space ${\cal
A}^{(n)}_{0\lambda}(u_0;t;0,q)$ agree, and are given by the span
\[
\langle P_\mu(x_1^{\pm 1},x_2^{\pm 1},\dots x_n^{\pm 1};q,t)\rangle_{\mu\subset\lambda},
\]
where $P_\mu$ is an ordinary Macdonald polynomial.
\end{cor}

\begin{rem}
It would be nice to have a direct proof of this Corollary, or the
corresponding result for a refined partial order; in particular, for the
dominance partial order, the limiting filtration should correspond to
dominance of monomials.
\end{rem}

\section{Type II transformations}

The connection coefficient formula for our biorthogonal functions, Theorem
\ref{thm:biorth_conn}, has a number of nice consequences for the
multivariate elliptic binomial coefficients.  For instance, by taking
$v=1$, we obtain the limiting case
\[
\lim_{b\to 1}
\frac{\Delta^0_\blambda(a|b;t;p,q)}
     {\Delta^0_\bkappa(a/b|1/b;t;p,q)}
\binomE{\blambda}{\bkappa}_{[a,b];t;p,q}=\delta_{\blambda\bkappa}.
\]
Also, if we perform the change of basis corresponding to $t_1\to t_1 v$,
then the change of basis corresponding to $t_1v\to t_1 vw$, the result
should be the same as if we directly changed $t_1\to t_1vw$.  We thus
obtain the following sum:

\begin{thm}\label{thm:binom_jackson}\cite{bctheta}
For otherwise generic parameters satisfying $bcde=pqa$,
\[
\binomE{\blambda}{\bkappa}_{[a,c];t;p,q}
=
\frac{\Delta^0_\bkappa(a/c|1/c,bd,be,pqa/b;t;p,q)}
     {\Delta^0_\blambda(a|c,bd,be,pqa/b;t;p,q)}
\sum_{\bkappa\subset\bmu\subset\blambda}
\Delta^0_\bmu(a/b|c/b,pqa,d,e;t;p,q)
\binomE{\blambda}{\bmu}_{[a,b];t;p,q}
\binomE{\bmu}{\bkappa}_{[a/b,c/b];t;p,q}.
\]
In particular,
\[
\sum_{\bkappa\subset\bmu\subset\blambda}
\binomE{\blambda}{\bmu}_{[a,b];t;p,q}
\binomE{\bmu}{\bkappa}_{[a/b,1/b];t;p,q}
=\delta_{\blambda\bkappa}.
\]
\end{thm}

\begin{rem}
Although this identity, along with the other sums mentioned in this
section, does indeed follow from Theorem \ref{thm:biorth_conn}, we should
mention that the argument in \cite{bctheta} proceeds in the opposite
direction, using these identities (and others) to prove the binomial
formula, and from this Theorem \ref{thm:biorth_conn}.  On the other hand,
the above argument provides a more straightforward {\em interpretation} of the
identity than that given in \cite{bctheta}.
\end{rem}

If we take $\blambda=(l,m)^n$, $\bkappa=0$ above, the above identity turns
out to be a product of two general instances of Warnaar's Jackson-type
summation (conjectured in \cite{WarnaarSO:2002}, and proved by Rosengren
\cite{RosengrenH:2001}).  Warnaar's Schlosser-type summation is also a special
case; see \cite{bctheta}.

Our reason for discussing this here is that there is an integral analogue
of the above sum, generalizing Theorem \ref{thm:normalization}.

\begin{thm}\label{thm:Rsint_spec1}
For otherwise generic parameters satisfying $t^{2n-2}t_0t_1t_2t_3u_0u_1=pq$,
\begin{align}
\langle
{\cal R}^{*(n)}_{\blambda}&(;t_0,u_0;t;p,q),
{\cal R}^{*(n)}_{\bkappa}(;t_1,u_1;t;p,q)
\rangle_{t_0,t_1,t_2,t_3,u_0,u_1;t;p,q}\notag\\
&=
\Delta^0_{\blambda}(t^{n-1}t_0/u_0|t^{n-1}t_0t_2,t^{n-1}t_0t_3;t;p,q)
\Delta^0_{\bkappa}(t^{n-1}t_1/u_1|t^{n-1}t_1t_2,t^{n-1}t_1t_3,t^{n-1}t_1t_0,t^{n-1}t_1u_0;t;p,q)\notag\\
&\phantom{{}={}}
{\cal R}^{*(n)}_{\blambda}(
\dots z_i(\bkappa;t_1/\sqrt{t^{n-1}t_1u_1})\dots
;t_0\sqrt{t^{n-1}t_1u_1},u_0\sqrt{t^{n-1}t_1u_1};t;p,q)
\end{align}
\end{thm}

\begin{proof}
Using the connection coefficient identity, we may express both interpolation
functions as linear combinations of biorthogonal functions.  Substituting
in the known values for the inner products of the biorthogonal functions,
we thus obtain a sum over partition pairs $\bmu\subset\blambda,\bkappa$.
That this sum gives the desired right-hand side is itself a special case
of the connection coefficient identity.

Alternatively, we can mimic the proof of Theorem \ref{thm:normalization},
using the fact that ${\cal D}^{(n)}_q(u_0,t_0,t_3;t,p)$ acts nicely on
${\cal R}^{*(n)}_{\blambda}(;t_0,u_0;t;p,q)$.  If we define
\begin{align}
F^{(n)}_{\blambda\bkappa}&(t_0,t_1,t_2,t_3,u_0,u_1;t;p,q)\\
&:=
\frac{
\langle
{\cal R}^{*(n)}_{\blambda}(;t_0,u_0;t;p,q),
{\cal R}^{*(n)}_{\bkappa}(;t_1,u_1;t;p,q)
\rangle_{t_0,t_1,t_2,t_3,u_0,u_1;t;p,q}}
{\Delta^0_{\blambda}(t^{n-1}t_0/u_0|t^{n-1}t_0t_2,t^{n-1}t_0t_3;t;p,q)
\Delta^0_{\bkappa}(t^{n-1}t_1/u_1|t^{n-1}t_1t_2,t^{n-1}t_1t_3;t;p,q)}
\notag
\end{align}
then adjointness gives
\begin{align}
F^{(n)}_{\blambda\bkappa}(t_0,t_1,t_2,t_3,u_0,u_1;t;p,q)
&=
F^{(n)}_{\blambda\bkappa}(q^{1/2}t_0,q^{-1/2}t_1,q^{1/2}t_2,q^{-1/2}t_3,q^{1/2}u_0,q^{-1/2}u_1;t;p,q)\\
&=
F^{(n)}_{\blambda\bkappa}(p^{1/2}t_0,p^{-1/2}t_1,p^{1/2}t_2,p^{-1/2}t_3,p^{1/2}u_0,p^{-1/2}u_1;t;p,q).
\end{align}
and thus
\[
F^{(n)}_{\blambda\bkappa}(t_0,t_1,t_2,t_3,u_0,u_1;t;p,q)
=
F^{(n)}_{\blambda\bkappa}(wt_0,t_1/w,t_2w,t_3/w,u_0w,u_1/w;t;p,q)
\]
for any $w\in \C^*$.  Taking the limit $w\to \sqrt{t^{n-1}t_1u_1}$ and
expanding via residue calculus, we obtain a sum over partition pairs
contained in $\bkappa_1^n$, in which only the term associated to $\bkappa$
survives.  (Recall that the contour must be deformed around the poles of
${\cal R}^{*(n)}_{\bkappa}$.)  We thus find
\[
F^{(n)}_{\blambda\bkappa}(t_0,t_1,t_2,t_3,u_0,u_1;t;p,q)
\propto
{\cal R}^{*(n)}_{\blambda}(
\dots z_i(\bkappa;t_1/\sqrt{t^{n-1}t_1u_1})\dots
;t_0\sqrt{t^{n-1}t_1u_1},u_0\sqrt{t^{n-1}t_1u_1};t;p,q).
\]
where the constant of proportionality is independent of $\blambda$.
This constant can be resolved by taking the limit
$w\to\sqrt{t^{n-1}t_0u_0}$ in the case $\blambda=0$.
\end{proof}

\begin{rems}
The left-hand side above is invariant under exchanging $(\blambda,t_0,u_0)$
and $(\bkappa,t_1,u_1)$.  That the right-hand side is invariant is a
special case of evaluation symmetry (Corollary \ref{cor:ev_symm}).  We can
also use that same special case of evaluation symmetry to see that this
generalizes Theorem \ref{thm:binom_jackson}.  Indeed, if we specialize so that
$t_0t_1=p^{-l}q^{-m}t^{1-n}$ with $\blambda,\bkappa\subset (l,m)^n$, then
the above left-hand side becomes a sum over $\bmu\subset (l,m)^n$.  Using
evaluation symmetry, the factor
\[
{\cal R}^{*(n)}_{\bkappa}(\dots z_i(\bmu;t_0)\dots;t_1,u_1)
=
{\cal R}^{*(n)}_{\bkappa}(\dots z_i((l,m)^n-\bmu;t_1)\dots;t_1,u_1)
\]
can be rewritten in terms of
\[
{\cal R}^{*(n)}_{\bmu}(\dots z_i((l,m)^n-\bkappa;x)\dots;x,y)
\]
for suitable $x$ and $y$.  Replacing $\bkappa$ by $(l,m)^n-\bkappa$ gives
Theorem \ref{thm:binom_jackson}.

Similarly, replacing $\bkappa$ by $(l,m)^n-\bkappa$ in the general version
and comparing the results, we find
\[
{\cal R}^{*(n)}_{l^n-\lambda,m^n-\mu}(;t_1,u_1;t;p,q)
=
\frac{(pq/t_1u_1)^{2l|\mu|+2m|\lambda|}\prod_{1\le i\le n} \theta(t_1 x_i^{\pm 1};p,q)_{l,m}}
     {\prod_{1\le i\le n} \theta((pq/u_1) x_i^{\pm 1};p,q)_{l,m}}
{\cal R}^{*(n)}_{\lambda,\mu}(;u_1/p^lq^m,p^lq^m t_1;t;p,q)
\label{eq:Rscomp}
\]
as both sides have the same inner product with ${\cal
  R}^{*(n)}_{\blambda}(;t_0,u_0;t;p,q)$.
\end{rems}

\begin{rems}
Note that the second proof of the theorem did not use the connection
coefficient identity, and is thus independent of \cite{bctheta}.
\end{rems}

If we take $\bkappa=0$ above, we obtain the following identity,
generalizing Kadell's lemma (see, for instance Corollary 5.14 of
\cite{bcpoly}).

\begin{cor}
For otherwise generic parameters satisfying $t^{2n-2}t_0t_1t_2t_3t_4t_5=pq$,
\[
\langle
{\cal R}^{*(n)}_{\blambda}(;t_0,t_1;t;p,q)
\rangle_{t_0,t_1,t_2,t_3,t_4,t_5;t;p,q}
=
\Delta^0_{\blambda}(t^{n-1}t_0/t_1|t^{n-1}t_0t_2,t^{n-1}t_0t_3,t^{n-1}t_0t_4,t^{n-1}t_0t_5;t;p,q).
\]
\end{cor}

The connection coefficient argument gives the following identity as well.

\begin{thm}\label{thm:Rsint_spec2}
For otherwise generic parameters satisfying
$t^{2n-2}t_0t_1t_2t_3u_0u_1=pq$,
\begin{align}
\langle
\cR^{*(n)}_\blambda(;t_0,u_0;t;p,q),
\cR^{*(n)}_\bkappa(;t_0,u_1;t;p,q)&
\rangle_{t_0,t_1,t_2,t_3,u_0,u_1;t;p,q}\notag\\
&=
\Delta^{0}_\blambda(\frac{t^{n-1}t_0}{u_0}|t^{n-1}t_0t_1,t^{n-1}t_0t_2,t^{n-1}t_0t_3,t^{n-1}t_0u_1;t;p,q)\notag\\
&\phantom{{}={}}
\Delta^{0}_\bkappa(\frac{t^{n-1}t_0}{u_1}|t^{n-1}t_0t_1,t^{n-1}t_0t_2,t^{n-1}t_0t_3,t^{n-1}t_0u_0;t;p,q)\notag\\
&\phantom{{}={}}
\tilde{R}^{(n)}_{\blambda}(\dots z_i(\bkappa,t'_0)\dots;t'_0{:}t'_1,t'_2,t'_3;u'_0,u'_1;t;p,q),
\end{align}
where the primed parameters are determined by
\begin{align}
t^{n-1}t'_0t'_1 &= t^{n-1}t_0t_1& t^{n-1}t'_0t'_2 &= t^{n-1}t_0t_2&
t^{n-1}t'_0t'_3 &= t^{n-1}t_0t_3\notag\\
t^{n-1}t'_0u'_0 &= t^{n-1}t_0u_0&
t^{n-1}t'_0u'_1 &= \frac{1}{t^{n-1}t_0u_1}& t^{\prime 2}_0 &= \frac{t_0}{t^{n-1}u_1}\notag
\end{align}
\end{thm}

\begin{rem}
The above transformation of the parameters is involutive, and conjugates
the exchange $u_0\leftrightarrow u_1$ to the ``hat'' transformation.
\end{rem}

A further application of connection coefficients gives the following
result, containing both Theorems \ref{thm:Rsint_spec1} and
\ref{thm:Rsint_spec2} as special cases.

\begin{thm}\label{thm:Rsint_spec}
For otherwise generic parameters satisfying
$t^{2n-2}t_0t_1t_2t_3u_0u_1=pq$,
\begin{align}
\langle
\cR^{*(n)}_\blambda(;t_0v,u_0;t;p,q),
\cR^{*(n)}_\bkappa(;t_0,u_1;t;p,q)&
\rangle_{t_0,t_1,t_2,t_3,u_0,u_1;t;p,q}\notag\\
&=
\Delta^0_\blambda(t^{n-1}t_0v/u_0|t^{n-1}vt_0t_1,t^{n-1}vt_0t_2,t^{n-1}vt_0t_3,t^{n-1}t_0u_1;t;p,q)\notag\\
&\phantom{{}={}}
\Delta^0_\bkappa(t^{n-1}t_0/u_1|t^{n-1}t_0t_1,t^{n-1}t_0t_2,t^{n-1}t_0t_3,t^{n-1}t_0u_0;t;p,q)\notag\\
&\phantom{{}={}}
\tilde{R}^{(n)}_{\blambda}(\dots
z_i(\bkappa,t'_0)\dots;t'_0v{:}t'_1,t'_2,t'_3;u'_0,u'_1/v;t;p,q),
\end{align}
with primed parameters as above.
\end{thm}

Now, Theorem \ref{thm:binom_jackson} is sufficiently general that the
univariate argument for deriving Bailey-type transformations from
Jackson-type summations applies, giving the following identity.

\begin{thm}\label{thm:binom_bailey}\cite{bctheta}
The sum
\[
\frac{\Delta^0_\blambda(a|b,apq/bf;t;p,q)}
     {\Delta^0_\bkappa(a/c|b/c,apq/bd;t;p,q)}
\sum_{\bkappa\subset\bmu\subset\blambda}
\frac{\Delta^0_\bmu(a/b|c/b,f,g;t;p,q)}
     {\Delta^0_\bmu(a/b|1/b,d,e;t;p,q)}
\binomE{\blambda}{\bmu}_{[a,b];t;p,q}
\binomE{\bmu}{\bkappa}_{[a/b,c/b];t;p,q}
\]
is symmetric in $b$ and $b'$, where $bb'de=capq$, $bb'fg=apq$.
\end{thm}

\begin{rems}
Again, taking $\blambda=(l,m)^n$, $\bkappa=0$ gives an identity
conjectured by Warnaar, in this case his conjectured multivariate
Frenkel-Turaev transformation \cite{WarnaarSO:2002}.
\end{rems}

\begin{rems}
This identity can also be obtained by comparing various ways of computing
connection coefficients for biorthogonal functions in which $t_0$, $t_3$,
$u_0$ are left fixed, but $t_1$, $t_2$, $u_1$ change.
\end{rems}

It turn out that this identity also has an integral analogue.  For each
integer $n\ge 0$, and partition pairs $\blambda$, $\bmu$ of length at most
$n$, we define a (meromorphic) function
\begin{align}
\II^{(n)}_{\blambda,\bmu}(t_0,t_1{:}t_2,t_3{:}t_4,t_5,t_6,t_7;t;p,q)\qquad&\\
:=
\frac{(p;p)^n(q;q)^n\Gamma(t;p,q)^n}{2^n n!}
\int_{C^n}&
{\cal R}^{*(n)}_{\blambda}(\dots x_i\dots;t_0,t_1;t;p,q)
{\cal R}^{*(n)}_{\bmu}(\dots x_i\dots;t_2,t_3;t;p,q)\notag\\
&
\prod_{1\le i<j\le n} 
\frac{\Gamma(t x_i^{\pm 1}x_j^{\pm 1};p,q)}
     {\Gamma(  x_i^{\pm 1}x_j^{\pm 1};p,q)}
\prod_{1\le i\le n}
 \frac{\prod_{0\le r\le 7} \Gamma(t_r x_i^{\pm 1};p,q)}
      {\Gamma(x_i^{\pm 2};p,q)}
\frac{dx_i}{2\pi\sqrt{-1}x_i},\notag
\end{align}
on the domain $t^{2n-2}t_0t_1t_2t_3t_4t_5t_6t_7=p^2q^2$,
where the contour $C^n$ is constrained in the usual way by the poles of the
integrand.

\begin{thm}\label{thm:int_bailey_lm}
If $t^{2n-2}t_0t_1t_2t_3t_4t_5t_6t_7=p^2q^2$ for some nonnegative integer
$n$, then
\begin{align}
\II^{(n)}_{\blambda,\bmu}(t_0,t_1{:}t_2,t_3{:}t_4,t_5,t_6,t_7;t;p,q)
&=
\Delta^0_{\blambda}(t^{n-1}t_0/t_1|t^{n-1} t_0t_4,t^{n-1} t_0t_5)
\prod_{1\le j\le n}
\prod_{\substack{r,s\in\{0,1,4,5\}\\r<s}}
\Gamma(t^{n-j} t_r t_s;p,q)
\notag\\
&\phantom{{}={}}
\Delta^0_{\bmu}(t^{n-1}t_2/t_3|t^{n-1} t_2t_6,t^{n-1} t_2t_7)
\prod_{1\le j\le n}
\prod_{\substack{r,s\in\{2,3,6,7\}\\r<s}}
\Gamma(t^{n-j} t_r t_s;p,q)
\notag\\
&\phantom{{}={}}
\II^{(n)}_{\blambda,\bmu}(t_0/u,t_1/u{:}ut_2,ut_3{:}t_4/u,t_5/u,ut_6,ut_7;t;p,q),
\end{align}
where $u$ is chosen so that
\[
u^2 = \sqrt{\frac{t_0t_1t_4t_5}{t_2t_3t_6t_7}} =
\frac{pqt^{1-n}}{t_2t_3t_6t_7}=\frac{t_0t_1t_4t_5}{pqt^{1-n}}.
\]
\end{thm}

\begin{proof}
If, following the second proof of Theorem \ref{thm:Rsint_spec1}, we attempt
to mimic the difference operator proof of Theorem \ref{thm:normalization},
we immediately encounter the difficulty that we no longer have adjointness
between two instances of ${\cal D}^{(n)}_q$, but rather between an instance
of ${\cal D}^{(n)}_q$ and an instance of ${\cal D}^{+(n)}_q$.  The one
exception is when $t^{n-1}t_0t_1t_4t_5=p$, in which case
\[
{\cal D}^{(n)}_q(t_0,t_1,t_4;t,p)
\]
and
\[
{\cal D}^{(n)}_q(\sqrt{q}t_2,\sqrt{q}t_3,\sqrt{q}t_6;t,p)
\]
are adjoint; the resulting transformation is precisely the case $u^2=q$ of
the theorem.

To extend this argument, we will thus need to extend the difference
operators.  With this in mind, we define a difference operator
$D^{(n)}_{l,m}(u_0,u_1,u_2;t;p,q)$ for $l,m\ge 0$ as follows.
\begin{align}
D^{(n)}_{0,1}(u_0,u_1,u_2;t;p,q)
&=
D^{(n)}_q(u_0,u_1,u_2,p/t^{n-1}u_0u_1u_2;t,p)\\
D^{(n)}_{1,0}(u_0,u_1,u_2;t;p,q)
&=
D^{(n)}_p(u_0,u_1,u_2,q/t^{n-1}u_0u_1u_2;t,q)\\
D^{(n)}_{l+l',m+m'}(u_0,u_1,u_2;t;p,q)
&=
D^{(n)}_{l,m}(u_0,u_1,u_2;t;p,q)
D^{(n)}_{l',m'}(S_{l,m}^{1/2}u_0,S_{l,m}^{1/2}u_1,S_{l,m}^{1/2}u_2;t;p,q),
\end{align}
where $S_{l,m}=(p,q)^{l,m}$. Since $D^{(n)}_{l,m}$ is a composition of $p$- and
$q$-difference operators, it itself is a difference operator; that it is
well-defined follows by verifying that the two ways of computing
$D^{(n)}_{1,1}$ agree.

\begin{lem}
Let $u_0$, $u_1$, $u_2$, $u_3$ be such that $t^{n-1}u_0u_1u_2u_3=pq/S_{l,m}$.
Then
\begin{align}
D^{(n)}_{l,m}(u_0,u_1,u_2;t;p,q)&
{\cal R}^{*(n)}_{\blambda}(;S_{l,m}^{1/2}u_0,S_{l,m}^{1/2}u_1;t;p,q)\\
&=
\frac{
\Delta^0_{\blambda}(t^{n-1}u_0/u_1|pq/u_1u_2,pq/u_1u_3;t;p,q)}
{\prod_{1\le i\le n}\prod_{0\le r<s\le 3} \Gamma(t^{n-i} u_ru_s;p,q)}
{\cal R}^{*(n)}_{\blambda}(;u_0,u_1;t;p,q)
\notag
\end{align}
In particular,
\[
D^{(n)}_{l,m}(u_0,u_1,u_2;t;p,q)
=
D^{(n)}_{l,m}(u_0,u_1,u_3;t;p,q).
\]
\end{lem}

\begin{proof}
The first claim holds when $(l,m)\in\{(0,1),(1,0)\}$; an easy induction
gives it in general.

Since $D^{(n)}_{l,m}(u_0,u_1,u_2;t;p,q)$ is a
difference operator, it is uniquely determined by this action; since the
given formula is symmetric between $u_2$ and $u_3$, the operator itself is
symmetric.
\end{proof}

\begin{lem}
The different instances of $D^{(n)}_{l,m}$ are related by
\[
D^{(n)}_{l,m}(u_0,u_1,u_2;t;p,q)
=
\prod_{\substack{1\le i\le n\\0\le r\le 3}}
\frac{\Gamma(u'_r x_i^{\pm 1};p,q)}
     {\Gamma(u_r x_i^{\pm 1};p,q)}
D^{(n)}_{l,m}(u'_0,u'_1,u'_2;t;p,q)
\prod_{\substack{1\le i\le n\\0\le r\le 3}}
\frac{\Gamma(S_{l,m}^{1/2}u_r x_i^{\pm 1};p,q)}
     {\Gamma(S_{l,m}^{1/2}u'_r x_i^{\pm 1};p,q)}
\]
\end{lem}

\begin{proof}
If $u'_3=u_3$, the result follows by a simple induction; the general case
then follows by combining that case with the symmetry between $u_0$, $u_1$,
$u_2$ and $u_3$.
\end{proof}

In particular, we can define the operator
\[
D^{(n)}_{l,m}(t;p,q)
:=
\prod_{\substack{1\le i\le n\\0\le r\le 3}}
  \Gamma(u_r x_i^{\pm 1};p,q)
D^{(n)}_{l,m}(u_0,u_1,u_2;t;p,q)
\prod_{\substack{1\le i\le n\\0\le r\le 3}}
  \frac{1}{\Gamma(S_{l,m}^{1/2}u_r x_i^{\pm 1};p,q)},
\]
which is independent of $u_0$, $u_1$, $u_2$.

This operator is self-adjoint with respect to the cross-terms in the $\II$
density; that is, with respect to 
\[
\prod_{1\le i<j\le n}
\frac{\Gamma(t x_i^{\pm 1} x_j^{\pm 1};p,q)}
     {\Gamma(x_i^{\pm 1} x_j^{\pm 1};p,q)}
\prod_{1\le i\le n}
\frac{dx_i}{2\pi\sqrt{-1}x_i}.
\]
This follows from the fact that
\[
\int
f D^{(n)}_{l,m}(t_0,t_1,t_2;t;p,q) g
\Delta^{(n)}(t_0,\dots,t_5;t;p,q)
=
\int
g D^{(n)}_{l,m}(t'_3,t'_4,t'_5;t;p,q) f
\Delta^{(n)}(t'_0,\dots,t'_5;t;p,q),
\]
where
\[
(t'_0,t'_1,t'_2,t'_3,t'_4,t'_5)=(S_{l,m}^{1/2}t_0,S_{l,m}^{1/2}t_1,S_{l,m}^{1/2}t_2,S_{l,m}^{-1/2}t_3,S_{l,m}^{-1/2}t_4,S_{l,m}^{-1/2}t_5),
\]
which in turn follows by induction from the cases $(l,m)\in
\{(0,1),(1,0)\}$.

We also have a sort of commutation relation satisfied by
$D^{(n)}_{l,m}(t;p,q)$.

\begin{lem}
If $l,m,l',m'$ are nonnegative integers and
$u_0u_1u_2u_3 = S_{l,m}S_{l',m'}p^2q^2$, then we have the following
identity of difference operators.
\begin{align}
D^{(n)}_{l,m}(t;p,q)&
\prod_{\substack{1\le i\le n\\0\le r\le 3}}
  \Gamma(u_r x_i^{\pm 1};p,q)
D^{(n)}_{l',m'}(t;p,q)
\prod_{\substack{1\le i\le n\\0\le r\le 3}}
  \frac{1}{\Gamma(S_{l',m'}^{-1/2} u_r x_i^{\pm 1};p,q)}
\\
&=
\prod_{\substack{1\le i\le n\\0\le r\le 3}}
  \Gamma(S_{l,m}^{-1/2}u_r x_i^{\pm 1};p,q)
D^{(n)}_{l',m'}(t;p,q)
\prod_{\substack{1\le i\le n\\0\le r\le 3}}
  \frac{1}{\Gamma((S_{l,m}S_{l',m'})^{-1/2} u_r x_i^{\pm 1};p,q)}
D^{(n)}_{l,m}(t;p,q)
\notag
\end{align}
\end{lem}

\begin{proof}
By comparing actions on ${\cal
  R}^{*(n)}_{\blambda}((S_{l,m}S_{l',m'})^{1/2}u_0,(S_{l,m}S_{l',m'})^{1/2}t_0;t;p,q)$,
we find that
\begin{align}
D^{(n)}_{l,m}(u_0,t_0,t_1;t;p,q)
&D^{(n)}_{l',m'}(S_{l,m}^{1/2}u_0,S_{l,m}^{1/2}t_0,S_{l,m}^{-1/2}t_2;t;p,q)\\
&=
D^{(n)}_{l',m'}(u_0,t_0,t_2;t;p,q)
D^{(n)}_{l,m}(S_{l',m'}^{1/2}u_0,S_{l',m'}^{1/2}t_0,S_{l',m'}^{-1/2}t_1;t;p,q)
\notag
\end{align}
Expressing this in terms of $D^{(n)}_{l,m}(t;p,q)$ and simplifying gives
the desired result.
\end{proof}

Now, consider an integral of the form
\[
\int
[D^{(n)}_{l,m}(t_0,t_1,t_4;t;p,q)f]
[D^{(n)}_{l',m'}(t_2,t_3,t_6;t;p,q)g]
\Delta^{(n)}(t_0,t_1,t_2,t_3,t_4,t_5,t_6,t_7;t;p,q)
\]
where $f\in {\cal A}^{(n)}(S_{l,m}^{1/2}t_1;p,q)$, $g\in {\cal
    A}^{(n)}(S_{l',m'}^{1/2} t_3;p,q)$, and the parameters satisfy the
    relations
\[
t^{n-1}t_0t_1t_4t_5 = pq \frac{S_{l',m'}}{S_{l,m}},\qquad
t^{n-1}t_2t_3t_6t_7 = pq \frac{S_{l,m}}{S_{l',m'}}.
\]
If we rewrite this integral in terms of $D^{(n)}(t;p,q)$ and
$\Delta^{(n)}(;t;p,q)$ and apply self-adjointness of
$D^{(n)}_{l,m}(t;p,q)$, the resulting composition of difference operators
can be transformed by the commutation relation.  The result is of the same
form, and we thus obtain the following identity.
\begin{align}
\int&
(D^{(n)}_{l,m}(t_0,t_1,t_4;t;p,q)f)
(D^{(n)}_{l',m'}(t_2,t_3,t_6;t;p,q)g)
\Delta^{(n)}(t_0,t_1,t_2,t_3,t_4,t_5,t_6,t_7;t;p,q)\\
&=
\int
(D^{(n)}_{l',m'}(t'_0,t'_1,t'_4;t;p,q)f)
(D^{(n)}_{l,m}(t'_2,t'_3,t'_6;t;p,q)g)
\Delta^{(n)}(t'_0,t'_1,t'_2,t'_3,t'_4,t'_5,t'_6,t'_7;t;p,q),
\notag
\end{align}
where\allowbreak
\[
t'_r =
\begin{cases}
(S_{l,m}/S_{l',m'})^{1/2}t_r & r\in \{0,1,4,5\}\\
(S_{l',m'}/S_{l,m})^{1/2}t_r & r\in \{2,3,6,7\}.
\end{cases}
\]
If we set
\begin{align}
f &= {\cal R}^{*(n)}_\blambda(;S_{l,m}^{1/2}t_0,S_{l,m}^{1/2}t_1;t;p,q),\\
g &= {\cal R}^{*(n)}_\bmu(;S_{l',m'}^{1/2}t_2,S_{l',m'}^{1/2}t_3;t;p,q),
\end{align}
we obtain the special case $u^2=S_{l',m'}/S_{l,m}$ of the theorem.  Since
this set is dense, the theorem holds in general.
\end{proof}

\begin{rems}
Note that the above proof did not use Theorem \ref{thm:biorth_conn}, and is
thus independent of the results of \cite{bctheta}.  In fact, one can use
this result to prove Theorem \ref{thm:biorth_conn}, as follows.  Connection
coefficients for interpolation functions can be obtained from the special
case $t_2t_7=pq$ (essentially Theorem \ref{thm:Rsint_spec}), by comparing
the result to that of Theorem \ref{thm:Rsint_spec1}.  One can then reverse
the first proof of Theorem \ref{thm:Rsint_spec1} to show that the functions
given by the binomial formula are indeed biorthogonal; Theorem
\ref{thm:biorth_conn} then follows via Theorem \ref{thm:binom_jackson}.
\end{rems}

\begin{rems}
In the special case $t^{n-1}t_0t_2=1/S_{l,m}$, we recover Theorem
\ref{thm:binom_bailey}.  Also, the univariate case $n=1$ is precisely the
case $n=m=1$ of the $A_n$ transformation.
\end{rems}

\begin{rems}
Similarly, using our integral operators, one can give a direct proof for
the case $u^2=t$, which presumably only extends to an argument valid for
$u^2\in t^\Z$.  This is, however, probably the simplest proof in the
univariate case (since then the integral is independent of $t$).
\end{rems}

We can simplify this transformation somewhat by adding an appropriate
normalization factor.  Define a meromorphic function
\[
\tilde\II^{(n)}_{\blambda,\bmu}(t_0,t_1{:}t_2,t_3{:}t_4,t_5,t_6,t_7;t;p,q)
:=
Z_{\blambda\bmu}
\II^{(n)}_{\blambda,\bmu}(t^{1/2}t_0,t^{1/2}t_1{:}t^{1/2}t_2,t^{1/2}t_3{:}t^{1/2}t_4,t^{1/2}t_5,t^{1/2}t_6,t^{1/2}t_7;t;p,q),
\]
where
\begin{align}
Z_{\blambda\bmu}&=
\prod_{0\le r<s\le 7} \Gamma^+(t t_r t_s;t,p,q) Z_\blambda Z_\bmu
,\\
Z_\blambda&=
\cC^0_{\blambda}(t^n,pq/tt_1t_2,pq/tt_1t_3;t;p,q)
\prod_{4\le r\le 7} \cC^0_{\blambda}(pq/tt_1t_r;t;p,q)\\
Z_\bmu&=
\cC^0_{\bmu}(t^n,pq/tt_0t_3,pq/tt_1t_3;t;p,q)
\prod_{4\le r\le 7} \cC^0_{\bmu}(pq/tt_3t_r;t;p,q)
\end{align}
and the condition on the parameters is now
$t^{2n+2}t_0t_1t_2t_3t_4t_5t_6t_7=p^2q^2$.  Here $\Gamma^+(x;t,p,q)$ is defined by
\[
\Gamma^+(x;t,p,q) :=
\prod_{i,j,k\ge 0} (1-t^i p^j q^k x)(1-t^{i+1} p^{j+1} q^{k+1}/x),
\]
so for instance
\[
\Gamma^+(tx;t,p,q) = \Gamma^+(x;t,p,q)\Gamma(x;p,q).
\]
Note that for generic $p,q,t$, the integer $n$ can be deduced from the
balancing condition on the parameters, and thus could in principle be
omitted from the notation for $\tilde\II$.  Note that it follows from the
residue calculus of the appendix that $\tilde\II_{00}$ is holomorphic for
each $n$; it may very well be holomorphic for $\blambda,\bmu\ne 0$, but
this would require a deeper understanding of the singularities of $\cR^*$
as a function of the parameters.

\begin{cor}\label{cor:int_bailey_lm}
We have
\[
\tilde\II^{(n)}_{\blambda,\bmu}(t_0,t_1{:}t_2,t_3{:}t_4,t_5,t_6,t_7;t;p,q)
=
\tilde\II^{(n)}_{\blambda,\bmu}(t_0/u,t_1/u{:}ut_2,ut_3{:}t_4/u,t_5/u,ut_6,ut_7;t;p,q),
\]
where $u$ is chosen so that
\[
u^2 = \sqrt{\frac{t_0t_1t_4t_5}{t_2t_3t_6t_7}} =
\frac{pqt^{-n-1}}{t_2t_3t_6t_7}=\frac{t_0t_1t_4t_5}{pqt^{-n-1}}.
\]
\end{cor}

Since $\tilde\II$ is also invariant under permutations of $t_4$, $t_5$,
$t_6$, $t_7$, it is in fact invariant under an action of the Weyl group
$D_4$.  Since there are three double cosets $S_4\backslash D_4/S_4$, there
is one other type of nontrivial transformation, namely:
\[
\tilde\II^{(n)}_{\blambda,\bmu}(t_0,t_1{:}t_2,t_3{:}t_4,t_5,t_6,t_7;t;p,q)
=
\tilde\II^{(n)}_{\blambda,\bmu}(u/t_1,u/t_0{:}u/t_3,u/t_2{:}v/t_4,v/t_5,v/t_6,v/t_7;t;p,q),
\]
where $u^2=t_0t_1t_2t_3$, $v^2=t_4t_5t_6t_7$, and $t^{n+1}uv=pq$.
In terms of the unnormalized integral, this reads
\begin{align}
\II^{(n)}_{\bmu,\bnu}(t_0,t_1{:}t_2,t_3{:}t_4,t_5,t_6,t_7;t;p,q)
={}&
\prod_{1\le j\le n}\prod_{0\le r\le 3,4\le s\le 7} \Gamma^+(t^{n-j} t_r t_s;t,p,q)\notag\\*
&
\Delta^0_{\bmu}(t^{n-1}t_0/t_1|t^nt_0t_4,t^nt_0t_5,t^nt_0t_6,t^nt_0t_7;t;p,q)\notag\\*
&\Delta^0_{\bnu}(t^{n-1}t_2/t_3|t^nt_2t_4,t^nt_2t_5,t^nt_2t_6,t^nt_2t_7;t;p,q)\notag\\*
&\II^{(n)}_{\bmu,\bnu}(u/t_1,u/t_0{:}u/t_3,u/t_2{:}v/t_4,v/t_5,v/t_6,v/t_7;t;p,q).
\end{align}

The reason for the factors of $t^{1/2}$ in the definition of
$\tilde\II^{(n)}$ is that the integral satisfies a further identity.

\begin{thm}
Let $n\ge m\ge 0$ be nonnegative integers such that
$\ell(\lambda),\ell(\bmu)\le m$, and suppose the parameters satisfy
$t^{n-m}t_0t_2=1$.  Then
\[
\tilde\II^{(n)}_{\blambda,\bmu}(t_0,t_1{:}t_2,t_3{:}t_4,t_5,t_6,t_7;t;p,q)
=
\tilde\II^{(m)}_{\blambda,\bmu}(1/t_2,t_1{:}1/t_0,t_3{:}t_4,t_5,t_6,t_7;t;p,q).
\]
\end{thm}

\begin{proof}
To compute the left-hand side, we must take a limit (as the condition on
the contour cannot be satisfied); what we find is that we must take
residues in $n-m$ of the variables, effectively setting those variables to
$t^{1/2} t_0,\dots,t^{n-m-1/2}t_0$, or equivalently (taking reciprocals) to
to $t^{n-m-1/2}t_2,\dots,t^{1/2}t_2$.  The result is the desired
$m$-dimensional instance of $\tilde\II$.
\end{proof}

\begin{rem} 
Note that the requirement that $\ell(\blambda),\ell(\bmu)\le m$ and $n\ge
m\ge 0$ with $n,m\in \Z$ is equivalent to a requirement that both sides be
well-defined.
\end{rem}

We thus find that we have a {\em formal} symmetry under a larger group,
isomorphic to the Weyl group $A_1D_4$.

If one of the partition pairs is trivial, the effective symmetry group
becomes larger.  To be precise, define
\[
\II^{(n)}_{\blambda}(t_0,t_1{:}t_2,t_3,t_4,t_5,t_6,t_7;t;p,q)
:=
\II^{(n)}_{\blambda,0}(t_0,t_1{:}t_2,t_3{:}t_4,t_5,t_6,t_7;t;p,q),
\]
and similarly for $\tilde\II^{(n)}_{\blambda}$.  This function is now
manifestly symmetric under permutations of $t_2$ through $t_7$; together
with the symmetry of Theorem \ref{thm:int_bailey_lm}, this gives rise to
the Weyl group $D_6$.  Since $S_6\backslash D_6/S_6$ has four double
cosets, we thus obtain a further transformation.

\begin{cor}
We have
\begin{align}
\II^{(n)}_{\blambda}(t_0,t_1{:}t_2,t_3,t_4,t_5,t_6,t_7;t;p,q)
&=
\Delta^0_{\blambda}(t^{n-1}t_0/t_1|t^{n-1}t_0t_2,\dots,t^{n-1}t_0t_7;t;p,q)
\prod_{1\le i\le n} \prod_{0\le r<s\le 7} \Gamma(t^{n-i} t_rt_s;p,q)\notag\\
&\phantom{{}={}}\II^{(n)}_{\blambda}(u/t_1,u/t_0{:}u/t_2,u/t_3,u/t_4,u/t_5,u/t_6,u/t_7;t;p,q),
\end{align}
where $u^2 = \sqrt{t_0t_1t_2t_3t_4t_5t_6t_7} = pq/t^{n-1}$.
\end{cor}

\begin{rem}
In the limit $t_6t_7=p^{l+1}q^{m+1}$, the right-hand side becomes a sum; taking
$\blambda=0$ and reparametrizing, we obtain the following integral
representation for a Warnaar-type sum:
\begin{align}
\II^{(n)}(pq/u_0,t_1,t_2,t_3,t_4,t_5,t_6,t_7;t;p,q)
={}&
\prod_{1\le i\le n}
  \frac{\Gamma(t^{n-i}t,t^{n-i}u_0^2;p,q)
        \prod_{1\le r<s\le 7} \Gamma(t^{n-i} t_rt_s;p,q)}
       {\prod_{1\le r\le 7} \Gamma(t^{n-i} u_0t_r;p,q)}\notag\\*
&\sum_{\bmu\subset (l,m)^n}
\Delta_{\bmu}(t^{n-1}u_0^2/pq|t^n,u_0/t_1,u_0/t_2,\dots,u_0/t_7;t;p,q),
\end{align}
assuming $t_1=p^lq^mu_0$ and
\[
t^{2n-2}t_2t_3t_4t_5t_6t_7=p^{1-l}q^{1-m}.
\]
Of course, other, less-symmetric, integral representations can be obtained
from transformations of the left-hand side.
\end{rem}

The formal group (adding in the dimension-changing transformation) now
becomes the Weyl group $D_7$; it is, however, unclear what significance
this has, since we cannot in general compose such transformations.  Thus
rather than obtaining the full Weyl group, we only obtain a union of two
$D_6\backslash D_7/D_6$ double cosets (out of three).  This gives rise to
several new dimension altering transformations, some of which correspond to
well-defined integrals.  Thus for instance, we find that
\[
\tilde\II^{(n)}_{\blambda}(t_0,t_1{:}t_2,t_3,t_4,t_5,t_6,t_7;t;p,q)
=
\tilde\II^{(n+m)}_{\blambda}(t_0/u,t_1u{:}t_2/u,t_3/u,t_4/u,t_5/u,t_6/u,t_7u;t;p,q),
\]
where
\[
u^2 = \sqrt{t_0t_2t_3t_4t_5t_6/t_1t_7} = \frac{pqt^{-n-1}}{t_1t_7} =
\frac{t_0t_2t_3t_4t_5t_6}{pqt^{-n-1}},
\]
such that $u^2=t^m$ with $m\in \Z$, $n,n+m\ge \ell(\blambda)$.  (In all,
there are essentially 9 distinct dimension altering transformations, coming
to the 12 legal $S_6\backslash D_7/S_6$ double cosets not in $D_6$ (modulo
inverses))

If $\blambda=0$, the group enlarges even further; in that case, the main
group is the Weyl group $E_7$, while the ``formal'' group is the Weyl group
$E_8$.  Moreover, the action of $E_8$ comes from the usual root system,
with roots of the form
\[
(\pm \frac{1}{2},\pm \frac{1}{2},\pm \frac{1}{2},\pm \frac{1}{2},\pm \frac{1}{2},\pm \frac{1}{2},\pm \frac{1}{2},\pm \frac{1}{2})
\]
(with an even number of $-$ signs) and permutations of
\[
(1,1,0,0,0,0,0,0),
(1,-1,0,0,0,0,0,0)
\]
(Thus, for instance, Corollary \ref{cor:int_bailey_lm} corresponds to the
reflection in the root
$(\frac{1}{2},\frac{1}{2},-\frac{1}{2},-\frac{1}{2},\frac{1}{2},\frac{1}{2},-\frac{1}{2},-\frac{1}{2})$.)
The subgroup $E_7$ is then the stabilizer of the root
$(\frac{1}{2},\frac{1}{2},\frac{1}{2},\frac{1}{2},\frac{1}{2},\frac{1}{2},\frac{1}{2},\frac{1}{2})$,
corresponding to $\sqrt{t_0t_1t_2t_3t_4t_5t_6t_7}=pq/t^{n+1}$.  Since there
are again four double cosets $S_8\backslash E_7/S_8$, we do not obtain any
new forms of the main transformation (and similarly for the dimension
altering transformation).  The various subgroups considered above are
related to $E_8$ as follows:
\begin{align}
E_7 &= Stab_{E_8}(\sqrt{t_0t_1t_2t_3t_4t_5t_6t_7})\\
D_7 &= Stab_{E_8}(\sqrt{t_1^3t_2t_3t_4t_5t_6t_7/t_0})\\
D_6 &= Stab_{E_8}(t_1/t_0,\sqrt{t_0t_1t_2t_3t_4t_5t_6t_7})\\
A_1D_4 &= Stab_{E_8}(
\sqrt{t_1^3t_2t_3t_4t_5t_6t_7/t_0},
\sqrt{t_0t_1t_3^3t_4t_5t_6t_7/t_2},
t_1t_3)\\
D_4 &= Stab_{E_8}(
t_1/t_0,t_3/t_2,t_1t_3,\sqrt{t_0t_1t_2t_3t_4t_5t_6t_7})
\end{align}

If $t=q$, the integrand can be expressed as a product of two determinants,
and is itself expressible as a determinant.  This gives rise to a system of
three-term quadratic recurrences, which turn out to be a form of Sakai's
elliptic Painlev\'e equation \cite{SakaiH:2001}; this generalizes the
result of \cite{KajiwaraK/MasudaT/NoumiM/OhtaY/YamadaY:2003} for the
univariate case (the authors of which also observed the existence of an
$E_7$ symmetry in that case).  This also generalizes results of
\cite{ForresterPJ/WitteNS:2002} at the Selberg level (showing that certain
Selberg-type integrals give solutions of the ordinary Painlev\'e
equations).  Also of interest are the cases $t=q^2$, $t=\sqrt{q}$, when the
integral can be expressed as a pfaffian, and thus satisfies a system of
four-term quadratic recurrences.  See \cite{recur} for more details.

%
%

Finally, to obtain a reasonable degeneration of the integral in the limit
$p\to 0$, we would need two ``upper'' parameters, of order $O(p)$, while
the remaining parameters would have order $O(1)$; we would then use the
fact that $\Gamma(pq/x;p,q)=\Gamma(x;p,q)^{-1}$ to move the upper
parameters to the denominator.  This property is in fact not invariant
under $E_8$, or even under the above $E_7$; instead we obtain a {\em
  different} instance of $E_7$ (as the stabilizer of the root
$(0,0,0,0,0,0,1,1)$, assuming the upper parameters are $t_6$ and $t_7$)
from the $E_8$ action, while the $E_7$ action reduces to $E_6$.  (The
one-dimensional instance of the resulting integral identity is a trivial
consequence of the hypergeometric series representation of Rahman
\cite{RahmanM:1986}).  If we further degenerate the integral to the
multivariate Askey-Wilson case (a.k.a. the Koornwinder density), the
symmetry group reduces to $D_5$, and the corresponding identity was proved
in \cite{bcpoly}.

\section{Appendix: Meromorphy of integrals}

In the above work, we have made heavy use of the fact that the various
contour integrals we consider are meromorphic functions of the parameters.
This does not quite follow from the meromorphy of the integrands, as can be
seen from the following two examples:
\begin{align}
\int_{|z|=1} \frac{1}{1+t(z+1/z)} \frac{dz}{2\pi\sqrt{-1}z}
&=
(1-4t^2)^{-1/2},\ |t|<1/4\\
\int_{|z|=1}
e^{1/(z-2)} \frac{dz}{2\pi\sqrt{-1}(z-t)}
&=
e^{1/(t-2)},\ |t|<1
\end{align}
In both cases, the integrand is meromorphic in a neighborhood of the
contour, but there are obstructions to meromorphically continuing the
integral.  (The second integrand, of course, has an essential singularity,
but so do the integrands of interest to us.)  It turns out, however, that
these are typical of the only two such obstructions: an initial contour
that separates branches of a component of the polar divisor of the
integrand, or such a component that leaves the domain of meromorphy.

We will prove this fact in Theorem \ref{thm:meromorphy} below, but first
need a lemma about meromorphy of residues.  Note that with $g$ as described
in the hypotheses of the lemma, the polar divisor of $g$ is a codimension 1
analytic subvariety of $D\times P$, and thus each component is either of
the form $D\times \chi_0$, or can be viewed as a family of point sets in
$D$ parametrized by $P$.

\begin{lem}
Let $D$ be a nonempty open subset of $\CP^1$, and let $P$ be an irreducible
normal holomorphic variety.  Let $g$ be a meromorphic function on $D\times
P$, and let $\chi$ be a component of the polar divisor of $g$ which is
closed in $\overline{D}\times P$.  Let $P'$ be the subset of $P$ on which
the fibers of $\chi$ are disjoint from the other polar divisors of
$g$ (the complement of a codimension 1 subvariety), and define a function
$f(p)$ on $P'$ by
\[
f(p) = \sum_{d:(d,p)\in \chi} \Res_{z=d} g(z,p).
\]
Then $f$ extends uniquely to a meromorphic function on all of $P$.
\end{lem}

\begin{proof}
Note that on any compact subset of $P$, $\chi$ is bounded away from the
complement of $D$, and thus its fibers lie in a compact subset of $D$, so
are finite in number.  In particular, the above sum is thus well-defined,
and gives a holomorphic function on $P'$.

Now, by Levi's theorem, we can freely remove any codimension 2 subvariety
$W$ of $P$ without affecting the extension of $f$; in particular, we may
assume that $P$ is regular (since its singular locus is codimension 2 by
normality).  We can then further restrict to a neighborhood of a general
regular point, to reduce to the case $P\subset \C^n$ an open polydisc; we
can also then write $g(d,p)=g_1(d,p)/g_2(d,p)$ for $g_1$, $g_2$
holomorphic.  Let $Z\subset P$ be the locus for which $g_2(d,p)$ is
identically 0 as a function of $d$; then by multiplying $g$ by a suitable
function of $p$ alone, we can remove all codimension 1 components of $Z$,
leaving a codimension 2 locus which can be removed by another application
of Levi's theorem.

Now, consider a point $p_0\in P$.  Reducing $P$ as necessary, we can assume
that the fiber of $\chi$ over $p_0$ consists of a single point $d_0$; we
can then reduce $D$ to assure that $g_2(d,p_0)$ also vanishes only at
$d_0$.  But then by the Weierstrass preparation theorem, $g_2(d,p)$ is the
product of a monic polynomial in $d$ with holomorphic coefficients and a
holomorphic function nowhere vanishing on $D$, which can be absorbed into
$g_1$.  Moreover, $g_2(d,p)$ factors as $h_1(d,p)h_2(d,p)$ where the monic
polynomial $h_1(d,p)$ vanishes precisely along $\chi$, and the monic
polynomial $h_2(d,p)$ is relatively prime to $h_1(d,p)$.  We can thus write
\[
g(d,p) = \frac{i_1(d,p)}{h_1(d,p)}+\frac{i_2(d,p)}{h_2(d,p)}+i_0(d,p),
\]
where $i_0$ is holomorphic in $d$, and $i_1$, $i_2$ are polynomials with
meromorphic coefficients of degree less than $\deg(h_1)$, $\deg(h_2)$
respectively.  But the above sum of residues is then precisely the leading
coefficient of $i_1(d,p)$, and is thus meromorphic in a neighborhood of
$p_0$.
\end{proof}

\begin{rem}
In fact, $i_1(d,p)\Delta(p)$ is holomorphic, where $\Delta(p)$ is the
resultant of the polynomials $h_1(d,p)$ and $h_2(d,p)$.
\end{rem}

Given a closed contour $C$ in $\CP^1$, every point not in $C$ of course
has an associated winding number; we extend this by linearity to formal
linear combinations of contours.

\begin{thm}\label{thm:meromorphy}
Let $D$ be a nonempty connected open subset of $\CP^1$, let $C$ be a finite
complex linear combination of contours in $D$, and let $P$ be an
irreducible normal holomorphic variety.  Let $g$ be a function meromorphic
on $D\times P$, and suppose the function $f$ is defined on an open subset
$U$ of $P$ by
\[
f(p) = \int_C g(z,p) dz
\]
(Thus, in particular, we assume that the polar divisor of $g$ in $D\times
U$ is disjoint from $C\times U$).

Suppose that each irreducible component $\chi$ of the polar divisor of $g$
is either of the form $D\times \chi_0$ or satisfies the assumptions:
\begin{itemize}
\item[1.] For every point $u\in U$, every point $d\in D$ such that
$(d,u)\in \chi$ has the same winding number with respect to $C$; call this
the winding number of $\chi$.
\item[2.] If $(d,p)$ is a limit point of $\chi$ in $\overline{D}\times P$
outside $D\times P$, then the winding number of $d$ with respect to $C$ is
the same as that of $\chi$ itself.
\end{itemize}
Then $f(p)$ extends uniquely to a meromorphic function on all of $P$.
\end{thm}

\begin{proof}
Again, we may as well assume that $P$ is an open polydisc in $\C^n$ for
some $n$.  We may then assume that the polar divisor of $g$ contains no
components of the form $D\times \chi_0$, since we can in that case simply
multiply $g$ by a holomorphic function to remove that pole.

Now, let $U'$ be an open subset of $P$, and consider a component $\chi$ of
the polar divisor of $g$ on $D\times U'$.  This is contained in a unique
component of the full polar divisor, with winding number $w_0$, say; on the
other hand, if $U'$ is not contained in $U$, $\chi$ can easily intersect
$C$ or have well-defined winding number different from the ``true'' winding
number $w_0$, in which case we call it ``problematical''.  We claim that
every point $p\in P$ has a neighborhood with only finitely many
problematical polar components.  Indeed, by condition $(2)$ above, we can
choose a bounded neighborhood $U'$ of $p$ such that the problematical
components of $g$ on $D\times U'$ are bounded away from the complement of
$D$, and are thus contained in $D'\times U'$ for some compact subset of
$D'$, in which $g$ can support only finitely many poles.

If $p$ is such that we can choose $U'$ so that the problematical components are
disjoint from all components with different ``true'' winding number (which
will hold for $p$ away from a codimension 1 subvariety), then we can obtain a
new contour $C'$ by deforming $C$ and adding small circles in such a way that
the integralp
\[
\int_C' g(z,p) dz
\]
is well-defined on $U'$, and such that on the intersection of two such
components, the functions agree.  Indeed, we can clearly deform $C$ in such a
way that the problematical components have well-defined winding numbers
w.r.to the new contour; by adding small circles around the problematical
components (shrinking $U'$ as necessary to allow these circles to be fixed)
we can make these winding numbers equal to the ``true'' winding numbers.  Any
two such contours will give the same integral, by Cauchy's theorem, and
thus these functions agree on intersections.

At a general point, we can still deform $C$ to give well-defined winding
numbers to the problematical components, but now have the difficulty that
they might intersect components with different winding numbers. Here, we
can observe that the above analytic continuation can be written as
\[
\int_C' g(z,p) dz + \text{finite sum of residues}
\]
where instead of adding a small circle around the problematical components,
we simply add the corresponding residue.  The first term is certainly
meromorphic (in fact, holomorphic near $p$); that the residue terms are
meromorphic (and thus that the theorem holds) results from the lemma.
\end{proof}

We need only the following special case.  Here
\[
(x;p,q)_\infty:=\prod_{0\le i,j} (1-p^i q^j x).
\]

\begin{cor}\label{cor:meromorphy}
Let $F(z;t_0,\dots,t_{m-1};u_0,\dots,u_{n-1};p,q)$ be a function
holomorphic on the domain
\[
z,t_0,\dots,t_{m-1},u_0,\dots,u_{n-1}\in \C^*,\ 0<|p|,|q|<1.
\]
Then the function defined for $|t_r|,|u_r|<1$ by
\[
G(t_0,\dots,t_{m-1};u_0,\dots,u_{n-1};p,q)
=
\prod_{0\le r<m,0\le s<n} (t_r u_s;p,q)_\infty
\int_{|z|=1}
\Delta(z;t_0,\dots,t_{m-1};p,q)
\frac{dz}{2\pi\sqrt{-1}z}
\]
where
\[
\Delta(z;t_0,\dots,t_{m-1};p,q)
=
\frac{F(z;t_0,\dots,t_{m-1};p,q)}
     {\prod_{0\le r<m} (t_r z;p,q)_\infty\prod_{0\le r<n} (u_r/z;p,q)_\infty},
\]
extends uniquely to a holomorphic function on the domain
\[
t_0,\dots,t_{m-1},u_0,\dots,u_{n-1}\in \C^*,\ 0<|p|,|q|<1.
\]
Away from the divisor of $\prod_{0\le r<m,0\le s<n} (t_r u_s;p,q)_\infty$,
this extension can be obtained by replacing the unit circle by any
(possibly disconnected) contour that contains the points $p^i q^j u_r$ and
excludes the points $1/p^i q^j t_r$, for $0\le i,j$.
\end{cor}

In particular, our multidimensional integrals can all be expressed as
iterated contour integrals of this form (in general restricted to a
subvariety of parameter space), so are meromorphic by straightforward
induction.  This does, however, tend to grossly overestimate the polar
divisor.  This overestimation can easily occur even in the one-dimensional
case, in the presence of symmetry.

In the case of the $BC_1$ integral, one role of the balancing condition, as
we have seen, is to make the summation limits factor into $p$-abelian and
$q$-abelian factors, which occurs because the density satisfies the
relation
\[
\Delta(p^i q^j z)\Delta(z) = \Delta(p^i z)\Delta(q^j z)
\]
for $i,j\in \Z$.  As observed by Spiridonov (personal communication), this
only determines the balancing condition up to a sign.  However,
one special case of this relation is the identity
\[
\Delta(\pm p^{i/2} q^{j/2})
\Delta(\pm p^{-i/2} q^{-j/2})
=
\Delta(\pm p^{i/2} q^{-j/2})
\Delta(\pm p^{-i/2} q^{j/2}),
\]
assuming both sides are defined; using the fact that
$\Delta(z)=\Delta(1/z)$, we conclude that
\[
\Delta(\pm p^{i/2} q^{j/2})^2
=
\Delta(\pm p^{i/2} q^{-j/2})^2
\]
so that
\[
\Delta(p^{i/2} q^{j/2}) = \pm \Delta(p^{i/2} q^{-j/2}),\qquad
\Delta(-p^{i/2} q^{j/2}) = \pm \Delta(-p^{i/2} q^{-j/2}).
\]
The balancing condition for the $BC_1$ integral then has the effect of
choosing the sign in this identity:
\[
\Delta(\pm p^{i/2} q^{j/2}) = -\Delta(\pm p^{i/2} q^{-j/2}).
\]
This motivates the hypotheses for the following result.

\begin{lem}
Let $\Delta(z;{\bf p})$ be a $BC_1$-symmetric function on $\C^*\times P$,
with $P$ an irreducible normal subvariety of
$\{t_0,t_1,\dots,t_{d-1},p,q\in \C^*:|p|,|q|<1\}$.  Suppose furthermore
that the following conditions are satisfied.
\begin{itemize}
\item[1.] The function
\[
\prod_{0\le r<d} (t_r z^{\pm 1};p,q) \Delta(z;{\bf p})
\]
is holomorphic.
\item[2.] At a generic point of $P$, the denominator has only simple zeros;
  it has triple zeros only in codimension 2.
\item[3.] For any integers $i,j$,
\[
\Delta(\pm p^{i/2} q^{j/2};{\bf p})
=
-\Delta(\pm p^{i/2} q^{-j/2};{\bf p}),
\]
as an identity of meromorphic functions on $P$.
\end{itemize}
Then the function on $P$ defined for $|t_r|<1$ by
\[
\prod_{0\le r<s<d} (t_r t_s;p,q)_\infty
\int_{|z|=1} \Delta(z) \frac{dz}{2\pi\sqrt{-1}z}
\]
extends to a holomorphic function on $P$.
\end{lem}

\begin{proof}
The integral extends meromorphically to this domain by Corollary
\ref{cor:meromorphy}; that it has at most simple poles along the subvarieties
$t_r t_s=p^{-i} q^{-j}$, $i,j\ge 0$ follows immediately from the fact that
at a generic point of such a subvariety, there are no higher-order
collisions of poles.  However, these considerations still leave open the
possibility that the given function might have poles along the subvarieties
$t_r^2 = p^{-i} q^{-j}$, $i,j\ge 0$.

We thus need, without loss of generality, to show that the above function
is holomorphic at a generic point ${\bf p}_0$ such that $t_0 = \pm p^{-l/2}
q^{-m/2}$, $i,j\ge 0$.  Now, in a neighborhood of such a point, the
analytic continuation is given by
\[
\int_C \Delta(z) \frac{dz}{2\pi\sqrt{-1}z},
\]
where $C=C^{-1}$ is a contour containing $p^i q^j t_r$ for $i,j\ge 0$,
$0\le r<d$.  Now, let $C'$ be a modified symmetric contour that still
contains $p^i q^j t_r$ for $r>0$ and $p^i q^j t_0$ for $i\ge l$ or $j\ge
m$, but excludes $p^i q^j t_0$ for $0\le i\le l$, $0\le j\le m$.  Then we
claim that
\[
\int_C \Delta(z) \frac{dz}{2\pi\sqrt{-1}z}
+
\int_{C'} \Delta(z) \frac{dz}{2\pi\sqrt{-1}z}
\]
is holomorphic on a neighborhood of ${\bf p}_0$.  Indeed, anywhere that two
poles coalesce, the poles have the same overall winding number with respect
to the two contours.  Thus to show the first term is holomorphic, it
suffices to prove that the difference of the two terms is holomorphic.  But
this is just a sum of residues; it is therefore sufficient to prove that
\[
\sum_{0\le i\le l}
\sum_{0\le j\le m}
\Res_{z=t_0 p^i q^j}
\Delta(z;{\bf p}),
\]
is holomorphic near ${\bf p}_0$, or in other words that
\[
\lim_{{\bf p}\to {\bf p}_0}
\sum_{0\le i\le l}
\sum_{0\le j\le m}
\Res_{z=t_0 p^i q^j}
\Delta(z;{\bf p})
\]
is well-defined.  We claim in fact that
\[
\lim_{{\bf p}\to {\bf p}_0}
(1-\pm p^{l/2}q^{m/2}t_0)
[
\Res_{z=t_0 p^i q^j}
\Delta(z;{\bf p})
+
\Res_{z=t_0 p^i q^{m-j}}
\Delta(z;{\bf p})
]
=
0,
\]
which then makes the poles of the summands cancel pairwise, giving the
desired result.  Now,
\[
\lim_{{\bf p}\to {\bf p}_0}
(1-\pm p^{l/2}q^{m/2}t_0)
\Res_{z=t_0 p^i q^j}
\Delta(z;{\bf p})
\propto
\lim_{{\bf p}\to {\bf p}_0}
\lim_{z\to t_0 p^i q^j}
(1-\pm p^{l/2}q^{m/2}t_0)
(1-t_0 p^i q^j/z)
\Delta(z;{\bf p}),
\]
and this limit is well-defined, again because at most two poles coalesce at
any given point.  Now, if we pull out the denominator factors $(t_0 z^{\pm
  1};p,q)_\infty$ of $\Delta(z;{\bf p})$, we can explicitly compute their
contributions to the limit, and use the fact that limits of holomorphic
functions can be exchanged to conclude that
\[
\lim_{{\bf p}\to {\bf p}_0}
\lim_{z\to t_0 p^i q^j}
(1-\pm p^{l/2}q^{m/2}t_0)
(1-t_0 p^i q^j/z)
\Delta(z;{\bf p})
=
\frac{1}{2}
\lim_{{\bf p}\to {\bf p}_0}
(1-\pm p^{l/2}q^{m/2}t_0)^2
\Delta(\pm p^{i-l/2}q^{j-m/2};{\bf p}).
\]
The claim follows.
\end{proof}

Similarly, for higher dimensional integrals, the simple inductive argument
leads to predictions of extremely high order poles along the divisors $t_r
t_s = p^{-l}q^{-m}$, $l,m\ge 0$.  That this does not occur for our
integrals follows via a similar argument from the fact that
\[
\Delta(p^a q^b z_0,p^c q^d z_0,z_3,\dots,z_n;{\bf p})
=
-\Delta(p^a q^d z_0,p^c q^b z_0,z_3,\dots,z_n;{\bf p})
\]
for our integrands; the consequence is that when moving the contour over a
given collection of poles of the form $p^i q^j t_0$, $0\le i\le l$, $0\le
j\le m$, the residues of residues that arise all cancel pairwise.  Somewhat
more generally, we have the following.

\begin{lem}
Let $\Delta(z_1,\dots,z_n;{\bf p})$ be a symmetric meromorphic function on
$(\C^*)^n\times P$, where $P$ is an irreducible normal subvariety of the
domain $\{t_0,t_1,\dots,t_{d-1},u_0,\dots,u_{d-1},p,q\in \C^*:|p|,|q|<1\}$.
Suppose furthermore that the following conditions are satisfied.
\begin{itemize}
\item[1.] The function
\[
\prod_{1\le i\le n}\prod_{0\le r<d} (t_r z_i,u_r/z_i;p,q) \Delta(z_1,\dots,
z_n;{\bf p})
\]
is holomorphic.
\item[2.] At a generic point of $P$, the denominator has only simple zeros.
\item[3.] For any integers $a,b,c,d\ge 0$,
\[
\Delta(p^a q^b z,p^c q^d z,z_3,\dots,z_n;{\bf p})
=
-\Delta(p^a q^d z,p^c q^b z,z_3,\dots,z_n;{\bf p}),
\]
as an identity of meromorphic functions on $(\C^*)^{n-1}\times P$.
\end{itemize}
For generic ${\bf p}\in P$, choose a contour $C_{\bf p}$ containing all
points of the form $u_r p^i q^j$ with $i,j\ge 0$, $0\le r<d$, and excluding
all points of the form $(t_r p^i q^j)^{-1}$ with $i,j\ge 0$, $0\le r<d$.  Let
$C'_{\bf p}$ be a similar contour that differs from $C_{\bf p}$ by
excluding the points $u_0 p^i q^j$ with $0\le i\le l$, $0\le j\le m$.  Then
\begin{align}
\int_{C_{\bf p}^n} \Delta(z;{\bf p}) \prod_{1\le i\le n}
\frac{dz_i}{2\pi\sqrt{-1}z_i}&{}
-
\int_{C_{\bf p}^{\prime n}} \Delta(z;{\bf p}) \prod_{1\le i\le n}
\frac{dz_i}{2\pi\sqrt{-1}z_i}\notag\\
&=
n
\int_{C_{\bf p}^{\prime n-1}}
\sum_{0\le (a,b)\le (l,m)}
\lim_{z_n\to p^a q^b t_0}
(1-p^a q^b t_0/z_n)
\Delta(z;{\bf p})
\prod_{1\le i<n}
\frac{dz_i}{2\pi\sqrt{-1}z_i}.
\end{align}
\end{lem}

\begin{proof}
A straightforward induction using the symmetry of the integrand tells us
that
\begin{align}
\int_{C_{\bf p}^n} \Delta(z;{\bf p}) \prod_{1\le i\le n}
\frac{dz_i}{2\pi\sqrt{-1}z_i}
&{}-
\int_{C_{\bf p}^{\prime n}} \Delta(z;{\bf p}) \prod_{1\le i\le n}
\frac{dz_i}{2\pi\sqrt{-1}z_i}\notag\\
&=
\sum_{1\le j\le n}
\int_{C_{\bf p}^{j-1}\times C_{\bf p}^{\prime n-j}}
\sum_{0\le (a,b)\le (l,m)}
\lim_{z_n\to p^a q^b t_0}
(1-p^a q^b t_0/z_n)
\Delta(z;{\bf p})
\prod_{1\le i<n}
\frac{dz_i}{2\pi\sqrt{-1}z_i};
\end{align}
indeed, the sum is simply the contributions from residues as we
deform the contours in $z_n$, $z_{n-1}$,\dots, $z_1$.  It thus suffices to
show that these $n$ terms all agree.  But the difference between the $j$th
term and the $j+1$st term is an $n-2$-dimensional integral of a sum of
double residues:
\[
\sum_{0\le (a,b),(c,d)\le (l,m)}
\lim_{z_n\to p^a q^b t_0}
\lim_{z_{n-1}\to p^c q^d t_0}
(1-p^a q^b t_0/z_n)
(1-p^a q^d t_0/z_{n-1})
\Delta(z_1,\dots,z_n;{\bf p})
\]
But again we can pull out the known pole factors and interchange limits of
the resulting holomorphic function; we conclude that the $(a,b),(c,d)$ and
$(a,d),(c,b)$ terms cancel.
\end{proof}

Applying this to the type I integral gives the following result; similar
results apply to the integral operators.

\begin{thm}
The function
\[
\prod_{0\le r<s\le 2m+2n+3} (t_r t_s;p,q)_\infty
I^{(m)}_{BC_n}(t_0,\dots,t_{2m+2n+3};p,q)
\]
extends to a holomorphic function on the domain $\prod_r t_r = (pq)^m$,
$|p|$, $|q|<1$.
\end{thm}

\begin{proof}
Indeed, the integrand satisfies the hypotheses of the two lemmas; the
second lemma readily shows that the integral has a simple pole along each
subvariety $t_r t_s p^a q^b=1$ (with residue equal to a sum of
$n-1$-dimensional integrals), while an induction using the first lemma
shows that the potential singularites for $t_r^2 p^a q^b=1$ are not
present.
\end{proof}

\begin{rem}
The situation for the $A_n$ integral is much more complicated, as we must
integrate against a test function in $Z$ to allow the use of a product
contour (which is legal since an inductive argument shows the $A_n$
integral to be meromorphic); but this extra integration, while preserving
meromorphy, can easily remove singularities.
\end{rem}

For the type II integral, a similar argument applies to contour
deformations; the additional poles coming from the cross terms are
sufficiently generic that the multidimensional lemma still holds.  There is
an important difference in that we have the additional constraint that the
contour $C=C^{-1}$ should contain the contours $t p^i q^j C$, $i,j\ge 0$. 
Thus when deforming through the collection of points $t^i q^j p^k t_0$,
$0\le (i,j,k)\le (a,b,c)$, it is necessary to first deform through the
points with $i=0$, then those with $i=1$, and so forth; otherwise the
contour constraint will be broken.  With that caveat, however, the results
still apply, and we obtain the following result.

\begin{thm}
Let $\II^{(m)}_n(t_0,t_1,\dots, t_{2m+3};t;p,q)$ be the $2m+4$-parameter
analogue of the type II integral, $m>0$.  The function
\[
\prod_{0\le i<n} \prod_{0\le r<s<2m} (t^i t_r t_s;p,q)_\infty
\II^{(m)}_n(t_0,\dots,t_{2m+3};t;p,q)
\]
extends to a holomorphic function on the domain $t^{2n-2}\prod_r
t_r=(pq)^m$, $|t|,|p|,|q|<1$.
\end{thm}

\begin{proof}
At a generic point with $t^i p^j q^k t_r t_s=1$, $r<s$, we can simply
deform the contour through the points $t^a p^b q^c t_r$, $0\le a\le i$
$0\le b\le j$, $0\le c\le k$ to obtain a holomorphic integral.  We thus
find that the desired integral is a sum of integrals over the new contour,
with integrands given by multiple residues at a sequence of points with
$a=0$, $a=1$,\dots The only such integrals that are singular at $t^i p^j
q^k t_r t_s=1$ are those involving $i+1$-tuple residues, and those have
simple poles.  Since we can obtain at most $n$-tuple residues from an
$n$-tuple integral, we conclude that we have at most simple poles, and
those only when $i<n$.

For $r=s$, if we first deform through the points $t^a p^b q^c t_r$ with
$0\le b\le j$, $0\le c\le k$, $0\le a<i/2$, we find that the only integrals
with possible singularities are those of $i/2$-tuple residues; these are
then generically holomorphic when $t^i p^j q^k t_r^2=1$ by induction.
\end{proof}

\begin{rem}
A similar result applies to $\II^{(n)}_{\blambda\bmu}$, with the caveat
that the interpolation functions may have poles independent of $z_1$,\dots
$z_n$; these poles would then in general survive as poles of the integral.
\end{rem}

\bibliographystyle{plain}

\begin{thebibliography}{10}

\bibitem{AndersonGW:1991}
G.~W. Anderson.
\newblock A short proof of {S}elberg's generalized beta formula.
\newblock {\em Forum Math.}, 3(4):415--417, 1991.

\bibitem{BarskyD/CarpentierM:1996}
D.~Barsky and M.~Carpentier.
\newblock Polyn\^omes de {J}acobi g\'en\'eralis\'es et int\'egrales de
  {S}elberg.
\newblock {\em Electron. J. Combin.}, 3(2):R1, 1996.

\bibitem{BaryshnikovY:2001}
Yu. Baryshnikov.
\newblock G{UE}s and queues.
\newblock {\em Probab. Theory Related Fields}, 119(2):256--274, 2001.

\bibitem{vanDiejenJF:1996}
J.~F. {\vphantom{Diejen}}van~Diejen.
\newblock Self-dual {K}oornwinder-{M}acdonald polynomials.
\newblock {\em Invent. Math.}, 126(2):319--339, 1996.

\bibitem{vanDiejenJF/SpiridonovVP:2000}
J.~F. {\vphantom{Diejen}}van~Diejen and V.~P. Spiridonov.
\newblock An elliptic {M}acdonald-{M}orris conjecture and multiple modular
  hypergeometric sums.
\newblock {\em Math. Res. Lett.}, 7(5-6):729--746, 2000.

\bibitem{vanDiejenJF/SpiridonovVP:2001}
J.~F. {\vphantom{Diejen}}van~Diejen and V.~P. Spiridonov.
\newblock Elliptic {S}elberg integrals.
\newblock {\em Internat. Math. Res. Notices}, (20):1083--1110, 2001.

\bibitem{DixonAL:1905}
A.~L. Dixon.
\newblock On a generalisation of {L}egendre's formula
  {$KE'-(K-E)K'=\frac{1}{2}\pi$}.
\newblock {\em Proc. London Math. Soc. (2)}, 3:206--224, 1905.

\bibitem{ForresterPJ/RainsEM:2005}
P.~J. Forrester and E.~M. Rains.
\newblock Interpretations of some parameter dependent generalizations of
  classical matrix ensembles.
\newblock {\em Probab. Theory Related Fields}, 131(1):1--61, 2005.

\bibitem{ForresterPJ/WitteNS:2002}
P.~J. Forrester and N.~S. Witte.
\newblock Application of the {$\tau$}-function theory of {P}ainlev\'e equations
  to random matrices: {${\rm P_{VI}}$}, the {JUE}, {C}y{UE}, c{JUE} and scaled
  limits.
\newblock {\em Nagoya Math. J.}, 174:29--114, 2004.

\bibitem{FrobeniusG:1882}
G.~Frobenius.
\newblock {\"U}ber die elliptischen {F}unctionen zweiter {A}rt.
\newblock {\em J. f{\"u}r die reine und angew. Math.}, 93:53--68, 1882.

\bibitem{GustafsonRA:1992}
R.~A. Gustafson.
\newblock Some {$q$}-beta and {M}ellin-{B}arnes integrals with many parameters
  associated to the classical groups.
\newblock {\em SIAM J. Math. Anal.}, 23(2):525--551, 1992.

\bibitem{KajiharaY/NoumiM:2003}
Y.~Kajihara and M.~Noumi.
\newblock Multiple elliptic hypergeometric series: an approach from the
  {C}auchy determinant.
\newblock {\em Indag. Math.}, 14:395--421, 2003.

\bibitem{KajiwaraK/MasudaT/NoumiM/OhtaY/YamadaY:2003}
K.~Kajiwara, T.~Masuda, M.~Noumi, Y.~Ohta, and Y.~Yamada.
\newblock {${}_{10}E_9$} solution to the elliptic {P}ainlev{\'e} equation.
\newblock {\em J. Phys. A}, 36(17):L263--L272, 2003.

\bibitem{KoornwinderTH:1992}
T.~H. Koornwinder.
\newblock {A}skey-{W}ilson polynomials for root systems of type {$BC$}.
\newblock In Donald St.~P. Richards, editor, {\em Hypergeometric functions on
  domains of positivity, Jack polynomials, and applications (Tampa, FL, 1991)},
  Contemp. Math. 138, pages 189--204. Amer. Math. Soc., Providence, RI, 1992.

\bibitem{OkounkovA:1998a}
A.~Okounkov.
\newblock {$BC$}-type interpolation {M}acdonald polynomials and binomial
  formula for {K}oornwinder polynomials.
\newblock {\em Transform. Groups}, 3(2):181--207, 1998.

\bibitem{OkounkovA:1998b}
A.~Okounkov.
\newblock ({S}hifted) {M}acdonald polynomials: {$q$}-integral representation
  and combinatorial formula.
\newblock {\em Compositio Math.}, 112(2):147--182, 1998.

\bibitem{RahmanM:1986}
M.~Rahman.
\newblock An integral representation of a {${}_{10}\phi_9$} and continuous
  bi-orthogonal {${}_{10}\phi_9$} rational functions.
\newblock {\em Can. J. Math.}, 38:605--618, 1986.

\bibitem{RahmanM/SuslovSK:1993}
M.~Rahman and S.~K. Suslov.
\newblock Classical biorthogonal rational functions.
\newblock In {\em Methods of approximation theory in complex analysis and
  mathematical physics (Leningrad, 1991)}, volume 1550 of {\em Lecture Notes in
  Math.}, pages 131--146. Springer, Berlin, 1993.

\bibitem{bctheta}
E.~M. Rains.
\newblock {$BC_n$}-symmetric abelian functions.
\newblock arXiv:math.CO/0402113.

\bibitem{recur}
E.~M. Rains.
\newblock Recurrences of elliptic hypergeometric integrals.
\newblock arXiv:math.CA/0504285.

\bibitem{bcpoly}
E.~M. Rains.
\newblock {$BC_n$}-symmetric polynomials.
\newblock {\em Transform. Groups}, to appear.
\newblock arXiv:math.QA/0112035.

\bibitem{ICMS}
E.~M. Rains.
\newblock A difference-integral representation of {K}oornwinder polynomials.
\newblock In V.~Kuznetsov and S.~Sahi, editors, {\em {J}ack,
  {H}all-{L}ittlewood and {M}acdonald polynomials}, {C}ontemporary
  {M}athematics. AMS, to appear.
\newblock arXiv:math.CA/0409437.

\bibitem{RichardsD/ZhengQ:2002}
D.~Richards and Q.~Zheng.
\newblock Determinants of period matrices and an application to {S}elberg's
  multidimensional beta integral.
\newblock {\em Adv. in Appl. Math.}, 28(3-4):602--633, 2002.

\bibitem{RosengrenH:2001}
H.~Rosengren.
\newblock A proof of a multivariable elliptic summation formula conjectured by
  {W}arnaar.
\newblock In {\em $q$-series with applications to combinatorics, number theory,
  and physics (Urbana, IL, 2000)}, volume 291 of {\em Contemp. Math.}, pages
  193--202. Amer. Math. Soc., Providence, RI, 2001.

\bibitem{RosengrenH:2004}
H.~Rosengren.
\newblock Elliptic hypergeometric series on root systems.
\newblock {\em Adv. Math.}, 181(3):417--447, 2004.

\bibitem{RosengrenH:2003}
H.~Rosengren.
\newblock New transformations for elliptic hypergeometric series on the root
  system {$A_n$}.
\newblock {\em Ramanujan J.}, to appear.
\newblock arXiv:math.CA/0305379.

\bibitem{RuijsenaarsSNM:1997}
S.~N.~M. Ruijsenaars.
\newblock First order analytic difference equations and integrable quantum
  systems.
\newblock {\em J. Math. Phys.}, 38:1069--1146, 1997.

\bibitem{SahiS:1999}
S.~Sahi.
\newblock Nonsymmetric {K}oornwinder polynomials and duality.
\newblock {\em Ann. of Math. (2)}, 150(1):267--282, 1999.

\bibitem{SakaiH:2001}
H.~Sakai.
\newblock Rational surfaces associated with affine root systems and geometry of
  the {P}ainlev{\'e} equations.
\newblock {\em Comm. Math. Phys.}, 220(1):165--229, 2001.

\bibitem{SpiridonovVP:2004}
V.~P. Spiridonov.
\newblock Short proofs of the elliptic beta integrals.
\newblock arXiv:math.CA/0408369.

\bibitem{SpiridonovVP:2001}
V.~P. Spiridonov.
\newblock On the elliptic beta function.
\newblock {\em Uspekhi Mat. Nauk}, 56(1(337)):181--182, 2001.

\bibitem{SpiridonovVP:2003}
V.~P. Spiridonov.
\newblock Theta hypergeometric integrals.
\newblock {\em Algebra i Analiz (St. Petersburg Math. J.)}, 15:161--215, 2003.

\bibitem{SpiridonovVP/ZhedanovAS:2000a}
V.~P. Spiridonov and A.~S. Zhedanov.
\newblock Classical biorthogonal rational functions on elliptic grids.
\newblock {\em C. R. Math. Rep. Acad. Sci. Canada}, 22(2):70--76, 2000.

\bibitem{SpiridonovVP/ZhedanovAS:2000b}
V.~P. Spiridonov and A.~S. Zhedanov.
\newblock Spectral transformation chains and some new biorthogonal rational
  functions.
\newblock {\em Commun. Math. Phys.}, 210:49--83, 2000.

\bibitem{VarchenkoAN:1989}
A.~N. Varchenko.
\newblock The determinant of a matrix of multidimensional hypergeometric
  integrals.
\newblock {\em Dokl. Akad. Nauk SSSR}, 308(4):777--780, 1989.

\bibitem{WarnaarSO:2002}
S.~O. Warnaar.
\newblock Summation and transformation formulas for elliptic hypergeometric
  series.
\newblock {\em Constr. Approx.}, 18(4):479--502, 2002.

\bibitem{WilsonJA:1991}
J.~A. Wilson.
\newblock Orthogonal functions from {G}ram determinants.
\newblock {\em SIAM J. Math. Anal.}, 22(4):1147--1155, 1991.

\end{thebibliography}

\end{document}